\documentclass[twoside,11pt]{article}
\usepackage{jmlr2e}
\usepackage{comment}
\newcommand{\papertitle}{Adaptation to the Range in $K$--Armed Bandits}
\newcommand{\authorsshort}{Hadiji and Stoltz}

\newcommand{\authorslong}{\name H{\'e}di Hadiji \email hedi.hadiji@math.u-psud.fr \\
\name Gilles Stoltz \email gilles.stoltz@math.u-psud.fr \\
\addr Universit{\'e} Paris-Saclay, CNRS, Laboratoire de math{\'e}matiques d'Orsay, 91405, Orsay, France
}

\jmlrheading{1}{2022}{1--23}{02/21}{$x$/22}{hadiji22a}{\authorsshort}

\ShortHeadings{\papertitle}{\authorsshort}
\firstpageno{1}

\newcommand{\bpc}[1]{\begin{proof}{{\textbf{#1}}}}
\newcommand{\bp}{\begin{proof}}

\newtheorem{reminder}{Reminder}
\newcommand{\stretchasneed}{}

\usepackage{fancyhdr}
\makeatletter
\renewcommand*{\@seccntformat}[1]{\csname the#1\endcsname.\quad}
\makeatother

\usepackage{natbib}

\usepackage{hyperref}       
\usepackage{url}            
\usepackage{booktabs}       

\usepackage[paperwidth=210mm,
                paperheight=297mm,
                top=2.5cm,
                textheight=247mm,
                footnotesep=2\baselineskip,
                left=2cm,
                right=2cm,
                marginparwidth=1.7cm,
                headsep=20pt,
                dvips=false]{geometry}
\usepackage{dsfont}
\usepackage{enumitem}
\usepackage{color}
\usepackage{indentfirst}
\usepackage{algorithm}
\usepackage{algorithmic}
\usepackage{subcaption}
\usepackage[font=footnotesize,labelfont=bf]{caption}
\usepackage{upgreek}
\usepackage{mathtools}
\usepackage{bigints}

\usepackage{color}
\usepackage{xcolor}
\definecolor{Red}{rgb}{1.00, 0.00, 0.00}
\definecolor{Blue}{rgb}{0.00, 0.00, 1.00}
\definecolor{Green}{rgb}{0.10, 0.7, 0}

\newcommand{\E}{\mathbb E}
\newcommand{\N}{\mathbb{N}}
\renewcommand{\P}{\mathbb P}
\newcommand{\R}{\mathbb R}

\newcommand{\nup}{\underline \nu}

\newcommand{\what}{\widehat}
\newcommand{\eps}{\varepsilon}
\renewcommand{\epsilon}{\varepsilon}

\newcommand{\AH}{\mbox{\rm \tiny AHB}}
\newcommand{\dep}{\mbox{\rm \tiny dep}}
\newcommand{\free}{\mbox{\rm \tiny free}}
\newcommand{\adv}{\mbox{\rm \tiny adv}}
\newcommand{\negent}{\mbox{\rm \tiny neg}}
\newcommand{\ts}{\mbox{\rm \tiny Ts}}

\newcommand{\Kinf}{\mathcal K_{\inf}}
\newcommand{\dint}{\, \text{d}}

\renewcommand{\d}{\mathrm{d}}
\renewcommand{\log}{\ln}
\newcommand{\wh}{\widehat}
\newcommand{\e}{\mathrm{e}}
\newcommand{\eqdef}{\stackrel{\mbox{\tiny def}}{=}}

\newcommand{\Dmod}{\mathcal D}
\newcommand{\cD}{\mathcal{D}}
\newcommand{\cF}{\mathcal{F}}
\newcommand{\cO}{\mathcal{O}}
\newcommand{\cS}{\mathcal{S}}
\newcommand{\cA}{\mathcal{A}}

\newcommand{\C}{C}

\newcommand{\dotptr}[2]{{#1^{\!\top}} #2}
\newcommand{\dotp}[2]{ \langle #1, \,  #2 \rangle }
\newcommand{\dotpb}[2]{ \bigl\langle #1, \,  #2 \bigr\rangle }
\newcommand{\dotpB}[2]{ \Bigl\langle #1, \,  #2 \Bigr\rangle }

\newcommand{\vertiii}[1]{{\left\vert\kern-0.25ex\left\vert\kern-0.25ex\left\vert #1
		\right\vert\kern-0.25ex\right\vert\kern-0.25ex\right\vert}}

\newcommand{\1}[1]{\mathds{1}_{\!\{  #1  \}}}
\renewcommand{\leq}{\leqslant}
\renewcommand{\geq}{\geqslant}

\DeclareMathOperator{\Ed}{E}
\newcommand{\argmax}{\mathop{\mathrm{argmax}}}
\newcommand{\argmin}{\mathop{\mathrm{argmin}}}

\DeclareMathOperator{\KL}{KL}
\DeclareMathOperator{\Var}{Var}
\DeclareMathOperator{\kl}{kl}

\DeclareMathOperator{\Tr}{Tr}
\DeclareMathOperator{\dom}{Dom}
\DeclareMathOperator{\interior}{Int}

\DeclareMathOperator{\segm}{Seg}

\newcommand{\changes}[1]{{{#1}}}

\newcommand{\secrefsuppl}[2]{{#2}}

\begin{document}

\title{\papertitle}
\author{\authorslong}
\editor{Ambuj Tewari}
\maketitle

\begin{abstract}%
We consider stochastic bandit problems with $K$ arms, each associated
with a distribution supported on a given finite range $[m,M]$.
We do not assume that the range $[m,M]$ is known and
show that there is a cost for learning this range. Indeed,
a new trade-off between distribution-dependent and distribution-free
regret bounds arises, which prevents from simultaneously
achieving the typical $\ln T$ and \smash{$\sqrt{T}$} bounds.
For instance, a \smash{$\sqrt{T}$} distribution-free regret bound
may only be achieved if the distribution-dependent regret bounds are
at least of order \smash{$\sqrt{T}$}. We exhibit a strategy achieving the rates
for regret imposed by the new trade-off.
\end{abstract}

\begin{keywords}
multiarmed bandits, adversarial learning, cumulative regret, information-theoretic proof techniques
\end{keywords}

\section{Introduction}
Stochastic multi-armed bandits form a standard setting
to deal with sequential decision-making problems
like the design of clinical trials---one of the first applications mentioned---or
online advertisement and online revenue management.

Except for notable exceptions discussed below,
virtually all articles on stochastic $K$--armed bandits
either assume that distributions of the arms belong to some parametric family---often,
one-dimensional exponential families---or are sub-Gaussian with a known
parameter $\sigma^2$. Among the latter category, the case of the non-parametric family of distributions
supported on a known range $[m,M]$ is of particular interest to us.

We show that the knowledge of the range $[m,M]$ is a crucial information
and that facing bounded bandit problems but ignoring the bounds $m$ and $M$
is much harder. We do so by studying what may be achieved and what cannot be achieved
anymore when the range $[m,M]$ is unknown and the strategies need to learn it.
We call this problem adaptation to the range, or scale-free regret minimization.
Why this problem is important and why we considered it is explained
in Section~\ref{sec:why}. \medskip

More precisely, we prove that adaptation to the range is actually possible but that it has a cost:
our most striking result (in Section~\ref{sec:main-simult})
is a severe trade-off between the scale-free distribution-dependent and distribution-free
regret bounds that may be achieved. For instance, no strategy adaptive to the range
can simultaneously achieve distribution-dependent regret bounds of order $\ln T$ and
distribution-free regret bounds of order \smash{$\sqrt{T}$} up to polynomial factors;
this is in contrast with the case of a known range where simple strategies like UCB strategies (by \citealp{UCB-orig})
do so.
Our general trade-off shows, for instance, that if one wants to keep the same \smash{$\sqrt{T}$}
order of magnitude for the scale-free distribution-free regret bounds, then the best scale-free distribution-dependent
rate that may be achieved is \smash{$\sqrt{T}$}.

We also provide (in Section~\ref{sec:AHB}) a strategy, based on exponential weights, that
adapts to the range and obtains optimal
distribution-dependent and distribution-free regret bounds in the eyes of the exhibited trade-off:
these are of respective orders $T^{1-\alpha}$ and $T^{\alpha}$, where
$\alpha \in [1/2,1)$ is a parameter of the strategy.

\subsection{Literature Review}
Optimal scale-free regret minimization under full monitoring for adversarial sequences is offered by the AdaHedge strategy
by~\citet{de2014follow}, which we will use as a building block in Section~\ref{sec:mainadahedge}.

For stochastic bandits, the main difficulty in adaptation to the range
is the adaptation to the upper end~$M$ (see~Remark~\ref{rk:th1OKDm+});
this is why \citet{honda15anon-asymptotic} could provide optimal $\ln T$ distribution-dependent regret bounds
for payoffs lying in ranges of the form $(-\infty,M]$, with a known~$M$.
\citet{lattimore2017scale} considers models of distributions with a known bound on their kurtosis,
which is a scale-free measure of the skewness of the distributions; he provides a scale-free algorithm
based on the median-of-means estimators, with $\ln T$ distribution-dependent regret bounds.
However, bounded bandits can have an arbitrarily high kurtosis, so our settings are not directly comparable.
\citet{cowan2015asymptotically} study adaptation to the range but in the restricted case of uniform distributions
over unknown intervals. They provide optimal $\ln T$ distribution-dependent regret bounds for that specific model:
in their model, the cost for adaptation is mild and lies only in the multiplicative constant before the $\ln T$.
In the setting of bounded bandits, we show that distribution-dependent regret bounds must be larger than $\ln T$,
but the argument of \citet[Remark~8]{lattimore2017scale} entails that any regret rate larger than $\ln T$, e.g., $(\ln T) \ln \ln T$,
may be achieved.
Similar results by~\citet{cowan2017normal} for Gaussian distributions with unknown means and variances were also obtained.

Finally, on the front of adversarial bandits, no prior work discussed adaptation to the range, to the best of our knowledge.

Additional important references performing adaptation in some other sense for stochastic and adversarial
$K$--armed bandits are discussed now, \changes{including some follow-up work
to this article}. \smallskip

\paragraph{Adaptation to the effective range or to unbounded ranges in adversarial bandits.}
\citet{gerchinovitz2016refined} show that it is impossible to adapt to the so-called effective range in adversarial bandits.
A sequence of rewards has effective range smaller than $b$ if for all rounds $t$, rewards $y_{t,a}$ at this round all
lie in an interval of the form $[m_t,M_t]$ with $M_t - m_t \leq b$.
The lower bound they exhibit relies on a sequence of changing intervals of fixed size. This problem is thus different from our setting.
See also positive results---regret upper bounds under additional assumptions---by \citet{cesa2017bandit} and \citet{thune2018adaptation} for adaptation to the effective range.

\citet{All06} deal with unbounded ranges $[m_t,M_t]$ in adversarial bandits and other
partial monitoring settings, where, e.g., $M_t = - m_t = t^\beta$ for some $\beta > 0$. They provide
regret upper bounds scaling with $t^{\beta/2}$ when $\beta$ is known, but do not detail the price to pay
for not knowing $\beta$---though they suggest to resort to a doubling trick in that case.
\smallskip

\paragraph{Adaptation to the variance.}
\citet{audibert2009exploration} consider a variant of UCB called UCB-V, which adapts to the unknown variance.
Its analysis assumes that rewards lie in a known range $[0, M]$. The results crucially use Bernstein's inequality,
which we state as Reminder~\ref{rm:bernstein_mart} in Appendix~\ref{sec:bernsrem}.
As Bernstein's inequality holds for random variables with supports in $(-\infty, M]$, the analysis of UCB-V might perhaps
be extended to this case as well.
Deviation bounds in Bernstein's inequality contain two terms, a main term scaling with the standard deviation,
and a remainder term, scaling with $M$. This remainder term, which seems harmless, is actually a true issue when $M$ is not known,
as shown by the results of the present article. \smallskip

\paragraph{Adaptation to other criteria.}
\citet{wei2018more}, \citet{zimmert2018optimal}, \citet{bubeck2017sparsity}, and many more,
provide strategies for adversarial bandits with rewards in a known range, say $[0,1]$, and
adapting to additional regularity in the data, like small variations or stochasticity of the data---but never
to the range itself.

\changes{
\paragraph{Follow-up works.}
Following an earlier version of this work, further research on range-adaptive bandit algorithms has been conducted. \citet{baudry2021from} bypass our lower bound by imposing minimal extra conditions on the reward distributions, avoiding the heavy-tail construction from Theorem~\ref{thm:adaptive_lower_bound}; in this context, they provide a fully range-adaptive algorithm. For adversarial multi-armed bandits,
\citet{rajaputtaagrawal} recover some small-loss bounds while being agnostic to the range, at the cost of degraded worst-case guarantees, and \citet{huang2021scale} obtain similar results under delayed feedback.
}

\subsection{Why Studying Adaptation to the Range for Finite-Armed Bandits}
\label{sec:why}

We encountered the problem of learning the range $[m,M]$ of bandits problems when
designing bandit algorithms for continuum-armed problems, see \citet{hedi-neurips}.
Therein, arms are indexed by some bounded interval $\mathcal{I}$, and the mean-payoff
function $f : \mathcal{I} \to \R$ is assumed to be smooth enough, e.g.,
H{\"o}lderian-smooth, with unknown regularity parameters $L$ and $\beta$.
The mean-payoff function $f$ has a bounded range, as it is continuous over a bounded
interval.
To optimally learn these smoothness parameters, histogram-reductions of
continuum-armed bandit problems to finite-armed bandits problems, of proper bandwidth,
are performed ({\`a} la \citealp{kleinberg}), by zooming out.
Any reasonable $K$--armed bandit algorithm may be used in this algorithmic scheme.
However, given this reduction, it had to be assumed that the range $[m,M]$ of $f$ is known,
as all $K$--armed bandit algorithms we were aware of assumed that the range of
the distributions over the arms was known. To get more complete adaptivity results in the continuum-armed case
and be able to ignore the range of the mean-payoff function $f$,
it was necessary and sufficient to deal with the similar issue of range adaptivity
in the case of finitely many arms---which this article provides.

We also believe that exhibiting impossibility results, like the existence of a severe trade-off
between distribution-dependent and distribution-free regret bounds in the case of the model
of bounded distributions with an unknown range, has consequences beyond that model.
This impossibility result holds in particular for all larger models, like
non-parametric models containing all distributions over the entire real line $\mathbb{R}$
satisfying certain assumptions on their tails to make sure
that they are not too large. We therefore provide some intrinsic limitation to learning in
$K$--armed stochastic bandits.

\changes{The techniques introduced extend to more complex settings, like linear bandits;
see \secrefsuppl{Appendix~D}{Appendix~\ref{sec:linear-app}}.}

\section{Setting and Main Results}
\label{sec:setting-stochastic}

We consider finitely-armed stochastic bandits with bounded and possibly signed rewards.
More precisely, $K \geq 2$ arms are available; we denote by $[K]$ the set $\{1,\ldots,K\}$ of these arms.
With each arm $a$ is associated a Borel probability distribution $\nu_a$ lying in
some known model $\Dmod$; a model is a set of Borel probability distributions over $\R$ with a first moment.
The models of interest in this article are discussed below; in the sequel,
we only consider Borel distributions even though we will omit this specification.

A bandit problem in $\Dmod$ is
a $K$--vector of probability distributions in $\Dmod$, denoted by
$\nup = (\nu_a)_{a \in [K]}$. The player knows $\Dmod$ but not $\nup$.
As is standard in this setting, we denote by $\mu_a = \Ed(\nu_a)$ the mean payoff provided by
an arm $a$. An optimal arm and the optimal mean payoff are respectively given by
$a^\star \in \argmax_{a \in [K]} \mu_a$ and $\mu^\star = \max_{a \in [K]} \mu_a$.
Finally, $\Delta_a = \mu^\star - \mu_a$ denotes the gap of an arm~$a$.

The online learning game goes as follows: at round $t \geq 1$, the player picks
an arm $A_t \in [K]$, possibly at random
according to a probability distribution $p_t = (p_{t,a})_{a \in [K]}$ based on an auxiliary randomization $U_{t-1}$,
e.g., uniformly distributed over $[0,1]$,
and then receives and observes a reward $Z_t$ drawn independently at random
according to the distribution $\nu_{A_t}$, given $A_t$.
More formally, a strategy of the player is a sequence of measurable mappings from the observations to the action set,
$(U_0, \, Z_1, U_1, \, \ldots, \, Z_{t-1}, U_{t-1}) \mapsto A_t$.
At each given time $T \geq 1$, we measure the performance of a strategy through its expected regret:
\begin{equation}
\label{eq:defregrstoch}
R_T(\nup) = T \mu^\star - \E \! \left[\sum_{t = 1}^{T} Z_t\right]
= T \mu^\star - \E \! \left[\sum_{t = 1}^{T} \mu_{A_t} \right] = \sum_{a = 1}^K \Delta_a \, \E\bigl[N_a(T)\bigr]\,,
\end{equation}
where we used the tower rule for the first equality and defined
$N_a(T)$ as the number of times arm $a$ was pulled between time rounds~$1$ and~$T$.

Doob's optional skipping (see \citealp[Chapter III, Theorem 5.2, page~145]{doob1953} for the original reference, see also
\citealp[Section~5.3]{CT88} for a more recent reference) shows that we may assume
that i.i.d.\ sequences of rewards $(Y_{t,a})_{t \geq 1}$ are drawn beforehand, independently at random, for each arm $a$
and that the obtained payoff at round $t \geq 1$ given the choice $A_t$ equals $Z_t = Y_{t,A_t}$.
We will use this second formulation in the rest of the paper as it is the closest
to the one of oblivious individual sequences described later in Section~\ref{sec:setting-oblivious}.
\changes{We may then assume that the auxiliary randomizations $U_0,U_1,\ldots$ are i.i.d.\ random variables independent from
the $(Y_{t,a})_{t \geq 1}$ and distributed according to a uniform distribution over $[0,1]$.}

\paragraph{Model: bounded signed rewards with unknown range.}
For a given range $[m,M]$, where $m < M$ are two real numbers, not necessarily nonnegative,
we denote by $\Dmod_{m, M}$ the set of probability distributions supported on $[m, M]$.
Then, the model corresponding to distributions with a bounded but unknown range is
the union of all such $\Dmod_{m, M}$:
\[
\Dmod_{-,+} = \displaystyle{\bigcup_{\substack{m,M \in \R : m < M}}} \Dmod_{m, M}\,.
\]

\subsection{Adaptation to the Range: Concept of Scale-Free Regret Bounds}
\label{sec:def}

Regret scales with the range length $M-m$, thus regret bounds involve a multiplicative factor $M-m$. We therefore consider
such bounds divided by the scale factor $M-m$ and call them scale-free regret bounds.
We denote by $\N$ the set of natural integers; rates on regret bounds will be given
by functions $\Phi : \N \to [0,+\infty)$. We define adaptation to the unknown range
in Definitions~\ref{def:free} and~\ref{def:dep} below.

\begin{definition}[Scale-free distribution-free regret bounds]
\label{def:free}
A strategy for stochastic bandits is
adaptive to the unknown range of payoffs with a scale-free distribution-free
regret bound $\Phi_{\free} : \N \to [0,+\infty)$
if for all real numbers $m < M$,
the strategy ensures, without the knowledge of $m$ and $M$:
\[
\forall \nup \ \mbox{\rm in } \Dmod_{m,M}, \ \ \forall T \geq 1,
\qquad R_T(\nup) \leq (M-m) \, \Phi_{\free}(T)\,.
\]
\end{definition}

We show in Section~\ref{sec:mainadahedge} that adaptation to the unknown range
may indeed be performed in the sense of Definition~\ref{def:free},
with a scale-free distribution-free regret bound of order \smash{$\sqrt{KT \ln K}$}.
The latter is optimal up to maybe a factor of $\smash{\sqrt{\ln K}}$ as
\citet{AuCBFrSc02} provided a lower bound \smash{$(1/20)\min\bigl\{\sqrt{KT},T\bigr\}$}
on the regret of any strategy against individual sequences in $[0,1]^K$, thus for bandit problems in $\cD_{0,1}$,
thus for scale-free distribution-free regret bounds.

\begin{definition}[Distribution-dependent rates for adaptation]
\label{def:dep}
A strategy for stochastic bandits is
adaptive to the unknown range of payoffs with a distribution-dependent
rate $\Phi_{\dep} : \N \to [0,+\infty)$
if for all real numbers $m < M$,
the strategy ensures, without the knowledge of $m$ and $M$:
\[
\forall \nup \ \mbox{\rm in } \Dmod_{m,M},
\qquad \limsup_{T \to +\infty} \frac{R_T(\nup)}{\Phi_{\dep}(T)} < +\infty\,.
\]
Put differently, the strategy ensures that $\limsup R_T(\nup)/\Phi_{\dep}(T) < +\infty$
for all $\nup \in \Dmod_{-,+}$.
\end{definition}

Definition~\ref{def:dep} does not add much to the classical notion of distribution-dependent rates
on regret bounds, as the scale factor $M-m$ does not appear in the definition;
it merely ensures that the strategy is not informed of the range.
Also, we are only interested in rates of convergence here, not
in the value of the finite limit of $R_T(\nup)/\Phi_{\dep}(T)$.
This limit however heavily depends on $\nup$, which justifies the terminology of
distribution-dependent \emph{rates} for adaptation $\Phi_{\dep}$.

In contrast, the bounds targeted in the distribution-free case have a finite-time, closed-form
expression, which is why we did not speak of rates in that case and rather referred to
scale-free distribution-free regret \emph{bounds} $\Phi_{\free}$.

\subsection{Scale-Free Distribution-Dependent Regret Bounds Considered in Isolation}
\label{sec:distrdep-isolation}

We first explain the impact of ignoring the range on distribution-dependent regret bounds.
What follows is discussed in greater detail in Appendix~\ref{sec:LB-nolnT}
as these results were already known or, at least, much expected.

When $m$ and $M$ are known, there exist several strategies ensuring
\[
\forall \nup \ \mbox{\rm in } \Dmod_{m,M},
\qquad \limsup_{T \to +\infty} \frac{R_T(\nup)}{\ln T} < +\infty\,,
\]
even with an optimal value of the limit; see the end of Section~\ref{sec:LRBK}.

Given Definition~\ref{def:dep}, one may therefore wonder whether $\Phi_{\dep} = \ln$
is achievable as a distribution-dependent rate for adaptation to the range.
Theorem~\ref{th:adaptmorethanlog} in Section~\ref{sec:A1} and the comment before its statement
provide a negative answer to this question.

However, a UCB-strategy with an increased exploration rate
given by a non-decreasing function $\varphi \gg \ln$ was suggested by \citet[Remark~8]{lattimore2017scale}
in the context of Gaussian bandits. It also works well in the setting of bounded bandits:
Theorem~\ref{th:UCB-incr} in Section~\ref{sec:A2}
states that it is adaptive to the unknown range of payoffs with a distribution-dependent
rate $\Phi_{\dep} = \varphi$. That is, any rate that is larger than a logarithm may be
achieved, including, for instance, $\varphi(t) = (\ln t) \ln \ln t$.

\subsection{Simultaneous Scale-Free Regret Bounds}
\label{sec:main-simult}

When the range $[m,M]$ of the payoffs is known, it is possible to simultaneously achieve
optimal distribution-free bounds, of order $\sqrt{KT}$,
and optimal distribution-dependent bounds, of order $\ln T$ with the optimal
constant recalled in Reminder~\ref{rem:LRBK} of Appendix~\ref{sec:LRBK}; see the KL-UCB-switch strategy by~\citet{garivier2018kl}.
Put differently, when the range of payoffs is known, one can achieve optimal asymptotic distribution-dependent
regret bounds while not sacrificing finite-time guarantees.
Simpler strategies like UCB strategies (see \citealp{UCB-orig})
also simultaneously achieve regret bounds of similar $\sqrt{T \ln T}$  and $\ln T$ orders of magnitude
but with suboptimal constants.
\citet{zimmert2018optimal} also provide a strategy, Tsallis-INF with $\alpha = 1/2$, that provides simultaneously
distribution-dependent regret guarantees of order $\ln T$, with suboptimal constants though, and
adversarial guarantees of order $\sqrt{KT}$, which are stronger than just distribution-free
guarantees.

\paragraph{First main result: existence of a trade-off.}
Our first main result states that getting simultaneously these $\ln T$ and $\sqrt{T}$ rates
is not possible anymore when the range of payoffs is unknown.

\begin{theorem}
\label{thm:adaptive_lower_bound}
Any strategy with a scale-free distribution-free regret bound
satisfying $\Phi_{\free}(T) = o(T)$ may only achieve distribution-dependent rates $\Phi_{\dep}$ for adaptation
satisfying $\Phi_{\dep}(T) \geq T/\Phi_{\free}(T)$.

More precisely,
the regret of such a strategy is lower bounded as follows: for all $\nup$ in $\Dmod_{-,+}$,
\begin{equation}
\label{eq:main-sim-LB}
\liminf_{T \to \infty } \, \frac{R_T( \nup)}{T / \Phi_{\free(T)}} \geq \frac{1}{16} \sum_{a=1}^K \Delta_a\,.
\end{equation}
\end{theorem}

The orders of magnitude of the scale-free distribution-free regret bounds $\Phi_{\free}(T)$ range between
the optimal $\sqrt{T}$ and the trivial $T$ rates. The distribution-dependent rates $\Phi_{\dep}$
for adaptation to the range are therefore at best $\sqrt{T}$ for strategies enjoying
scale-free distribution-free regret bounds; $\ln T$ rates are excluded.
More generally, Theorem~\ref{thm:adaptive_lower_bound} shows that there is a trade-off:
to force faster distribution-dependent rates for adaptation, one must suffer worsened
scale-free distribution-free regret bounds.

The proof of Theorem~\ref{thm:adaptive_lower_bound} is provided in Section~\ref{sec:proofthmainLB}.
It actually provides a finite-time, but messy, lower bound on
$R_T( \nup) \big/ \bigl( T / \Phi_{\free(T)} \bigr)$.

\paragraph{Second main result: achieving the trade-off.}
Our second main result consists of showing that the trade-off imposed by
Theorem~\ref{thm:adaptive_lower_bound} may indeed be achieved. Section~\ref{sec:AHB}
will introduce a strategy, relying on a parameter $\alpha \in [1/2,1)$ and called AHB---which stands for
AdaHedge for $K$--armed Bandits with extra-exploration; see Algorithm~\ref{algo:adahedge}.
Theorems~\ref{th:adv_talpha_bound} and~\ref{thm:distrib-dependent} show in particular
that AHB adapts to the unknown range, satisfies a scale-free distribution-free regret bound
\[
\Phi^{\AH}_{\free}(T) = \biggl( 3 + \frac{5}{\sqrt{1-\alpha}} \biggr) \sqrt{K \ln K} \,\, T^{\alpha} + 10 K\ln K = \cO(T^\alpha)\,,
\]
and achieves a distribution-dependent rate for adaptation $\Phi^{\AH}_{\dep}(T) = T/\Phi^{\AH}_{\free}(T) = \cO(T^{1-\alpha})$.
Like \citet{zimmert2018optimal}, we are actually able to prove
an adversarial regret bound, not only the mentioned distribution-free regret bound.

Even better, Theorem~\ref{thm:distrib-dependent} states that for all $\nup$ in $\Dmod_{-,+}$,
\begin{equation}
\label{eq:main-sim-UB}
\limsup_{T \to \infty} \frac{R_T(\nup)}{T / \Phi^{\AH}_{\free(T)}}
\leq  \frac{12 \ln K}{1-\alpha} \sum_{a = 1}^K \Delta_a \,.
\end{equation}

\paragraph{Discussion.}
The distribution-dependent constants in the right-hand sides of
\eqref{eq:main-sim-LB} and~\eqref{eq:main-sim-UB} are proportional to
the sums of the gaps,
\[
G(\nup) = \sum_{a = 1}^K \Delta_a\,,
\]
and differ from this sum only by
\changes{distribution-free} factors of $1/16$ and $(12 \ln K)/(1-\alpha)$.
The quantity $G(\nup)$ appears as a new measure of the underlying geometry of information.
We have no deep interpretation thereof, but may despite all
underline a fundamental difference in our setting
compared to the setting of a known range.

When \changes{the payoff range} is unknown,
the optimal distribution-dependent number of pulls of a suboptimal arm may be bounded independently of $\nup$.
The proof of Theorem~\ref{thm:adaptive_lower_bound} in Section~\ref{sec:proofthmainLB}
indeed shows that for all suboptimal arms $a \in [K]$,
\[
\liminf_{T \to +\infty} \, \frac{\E_{\nup}\big[N_a(T)\big]}{T/\Phi_{\free}(T)} \geq \frac{1}{16}\,.
\]
This is in contrast with the case of a known range,
for which the bound of Reminder~\ref{rem:LRBK} of Appendix~\ref{sec:LRBK} is optimal and
strongly depends on $\nu_a$ and $\mu^{\star}$.

The reason for this is that when ignoring the range, the player needs to be a lot more conservative in the exploitation and explore more often.
Indeed, to maintain the distribution-free regret bound, the player must avoid the catastrophic case in which an apparently suboptimal arm turns out to be good because of large rewards occurring with small probability, i.e., because of heavy-tail-like issues.
For this reason, the player must pull suboptimal arms more frequently than in the case of a known range.
This intuition is supported by the construction in the lower bound presented in Section~\ref{sec:proofthmainLB}:
the alternative bounded problem $\nup'$ against a problem $\nup$ has an arm satisfying $\P_{Y \sim \nu_a'}[Y \geq \mu_a + 2 \Delta_a / \eps] = \eps$.
This behavior is indeed reminiscent of issues arising with heavy-tailed distributions.

\section{Proof of Theorem~\ref{thm:adaptive_lower_bound}: Existence of a Trade-Off}
\label{sec:proofthmainLB}

We follow a proof technique introduced by
\citet{lai1985asymptotically} and \citet{BuKa96} and recently revisited by~\citet{garivier2018Explore}.
We fix some bandit problem $\nup$ in $\Dmod_{-,+}$ and construct an alternative
bandit problem $\nup'$ in $\Dmod_{-,+}$ by modifying the distribution of a single suboptimal arm~$a$
to make it optimal. This is always possible, as there is no bound on the upper end on the ranges of
the payoffs in the model. We apply a fundamental inequality that links the expectations of the numbers
of times $N_a(T)$ that $a$ is pulled under $\nup$ and $\nup'$. We then substitute inequalities
stemming from the definition of distribution-free scale-free regret bounds $\Phi_{\free}$,
and the result follows by rearranging all inequalities.

\paragraph{Step 1: Alternative bandit problem.} The lower bound is trivial---it equals~$0$---when all
arms of $\nup$ are optimal. We therefore assume that at least one arm is suboptimal and
fix such an arm~$a$.
For some $\epsilon \in [0,1]$ to be defined later by the analysis, we introduce
the alternative problem $\nup' = (\nu'_k)_{k \in [K]}$ with $\nu'_k = \nu_k$ for $j \neq a$ and
$\nu'_a = (1 - \epsilon ) \nu_a + \epsilon \delta_{\mu_a + 2\Delta_a/\epsilon}$.
This distribution $\nu'_a$ has a bounded range, so that $\nup'$ lies indeed in $\Dmod_{-,+}$.
The expectation of $\nu'_a$ equals $\mu'_a = \mu_a + 2 \Delta_a = \mu^\star + \Delta_a > \mu^\star$.
Thus, $a$ is the only optimal arm in $ \nup'$.
Finally, \changes{for $\epsilon < 2\Delta_a/(M-\mu_a)$, the point $\mu_a + 2 \Delta_a/\epsilon$ is larger than $M$
and thus} lies outside of the bounded support of $\nu_a$. In that case, the density of
$\nu_a$ with respect to $\nu'_a$ is given by $1/(1-\epsilon)$ on the support of $\nu_a$
and $0$ elsewhere, so that $\KL(\nu_a, \nu'_a) = \ln\bigl(1/(1-\epsilon)\bigr)$.

\paragraph{Step 2: Application of a fundamental inequality.}
We denote by $\kl(p,q)$ the Kullback-Leibler divergence between Bernoulli distributions with parameters $p$ and $q$.
We also index expectations in the rest of this proof only by the bandit problem they are relative to: for instance,
$\E_{ \nup}$ denotes the expectation of a random variable when the ambient randomness is given by the bandit problem $\nup$.
The fundamental inequality for lower bounds on the regret of stochastic bandits (\citealp{garivier2018Explore}, Section~2, Equation~6),
which is based on the chain rule for Kullback-Leibler divergence and on a data-processing inequality for expectations of $[0,1]$--valued
random variables, reads:
\[
\kl\! \left(\frac{\E_{ \nup}\big[N_a(T)\big]}{T}, \, \frac{\E_{ \nup'}\big[N_a(T)\big]}{T} \right)
\leq \E_{\nup}\big[N_a(T)\big] \, \KL(\nu_a, \nu'_a)
= \E_{\nup}\big[N_a(T)\big] \ln\bigl(1/(1-\epsilon)\bigr) \,.
\]
Now, since $u \in (-\infty,1) \mapsto -u^{-1}\ln (1-u)$ is increasing,
we have $\ln\bigl(1/(1-\epsilon)\bigr) \leq (2 \ln 2) \epsilon $
for $\epsilon \leq  1/ 2$. For all $(p,q) \in [0,1]^2$ and with the usual measure-theoretic conventions,
\[
\kl(p, q) = \underbrace{p \ln p + \changes{(1-p) \ln (1-p)}}_{\geq - \ln 2} + \underbrace{p \ln \frac{1}{q}}_{\geq 0} + (1-p) \ln \frac{1}{1-q}
\geq (1-p) \ln \frac{1}{1-q} - \ln 2\,,
\]
so that, putting all inequalities together, we have proved
\begin{equation}\label{eq:lower_bound_calc}
\left(1- \frac{\E_{ \nup}\big[N_a(T)\big]}{T} \right) \, \ln  \! \left(\frac{1}{1 - \E_{ \nup'}\big[N_a(T)\big] / T} \right)  - \ln 2
\leq (2 \ln 2) \, \epsilon \, \E_{ \nup}\big[N_a(T)\big]\,.
\end{equation}

\changes{In this step, we only imposed the constraint $\epsilon \in [0,1/2]$. We recall that in the
previous step, we imposed $\epsilon < 2\Delta_a/(M-\mu_a)$. Both conditions are implied
by $\epsilon \leq \Delta_a/\bigl(2(M-\mu_a)\bigr)$, which we will assume in the sequel.}

\paragraph{Step 3: Inequalities stemming from the definition of scale-free distribution-free regret bounds.}
We denote by $[m,M]$ a range containing the supports of all distributions of $\nup$.
By definition of $\Phi_{\free}$, given that $a$ is a suboptimal arm (i.e., $\Delta_a > 0$):
\[
\Delta_a \, \E_{\nup}[N_a(T)] \leq R_T( \nup) \leq (M-m) \, \Phi_{\free}(T)\,.
\]
\changes{We now prove a similar inequality for $\nup'$, for which
we recall that $a$ is the unique optimal arm.
We denote by $\Delta'_k = \mu'_a - \mu_k$ the gap of arm $k$ in $\nup'$.}
By the definition of $\nu'_a$,
the distributions of $\nup'$ have supports within the range $[m,M_\epsilon]$,
where we denoted \changes{$M_\epsilon = \max\{M,\,\mu_a+2\Delta_a/\epsilon\} = \mu_a+2\Delta_a/\epsilon$,
given the condition imposed on $\epsilon$.}
Therefore, by definition of
$\Phi_{\free}$, and given that all gaps $\Delta'_k$ are larger than the gap
$\Delta'_a = \mu'_a - \mu^\star = \Delta_a$ between the unique optimal arm $a$ of $\nup'$ and the second best
arm(s) of $\nup'$ (which were the optimal arms of $\nup$), we have
\[
\Delta_a \bigl( T - \E_{\nup'} [N_a(T)] \bigr) =
\Delta'_a \bigl( T - \E_{\nup'} [N_a(T)] \bigr) \leq
\sum_{j \ne a} \Delta'_j \, \E_{\nup'} [N_j(T)] = R_T(\nup') \leq (M_\epsilon-m) \, \Phi_{\free}(T)\,.
\]
By rearranging the two inequalities above, we get
\[
1- \frac{\E_{ \nup}\big[N_a(T)\big]}{T} \geq 1 - \frac{(M-m) \, \Phi_{\free}(T)}{T \Delta_a}
\qquad \mbox{and} \qquad
1- \frac{\E_{ \nup'}\big[N_a(T)\big]}{T} \leq \frac{(M_\epsilon-m) \, \Phi_{\free}(T)}{T \Delta_a}\,,
\]
thus, after substitution into~\eqref{eq:lower_bound_calc},
\begin{equation}
\label{eq:almostfinaltradeoff}
\left( 1 - \frac{(M-m) \, \Phi_{\free}(T)}{T \Delta_a} \right) \, \ln  \! \left( \frac{T \Delta_a}{(M_\epsilon-m) \, \Phi_{\free}(T)} \right)
- \ln 2 \leq (2 \ln 2) \, \epsilon \, \E_{ \nup}\big[N_a(T)\big]\,.
\end{equation}

\paragraph{Step 4: Final calculations.}
We take $\epsilon = \epsilon_T = \alpha^{-1} \, \Phi_{\free}(T) / T$ for some constant $\alpha > 0$; we will pick $\alpha = 1/8$.
By the assumption $\Phi_{\free}(T) = o(T)$, \changes{we have $\epsilon_T \leq \Delta_a/\bigl(2(M-\mu_a)\bigr)$, as needed, for $T$ large enough}, as well as
$M_{\epsilon_T} = \mu_a+2\Delta_a/\epsilon_T = \mu_a+2\alpha\Delta_a T/ \Phi_{\free}(T)$.
Substituting these values into~\eqref{eq:almostfinaltradeoff},
a finite-time lower bound on the quantity of interest is finally given by
\[
\frac{\E_{\nup}\big[N_a(T)\big]}{T/\Phi_{\free}(T)} \geq
\frac{\alpha}{2 \ln 2} \Biggl( - \ln 2 + \biggl( 1 - \underbrace{\frac{(M-m) \, \Phi_{\free}(T)}{T \Delta_a}}_{\to 0}
\biggr) \, \ln \biggl( \underbrace{\frac{T \Delta_a}{2\alpha\Delta_a T + (\mu_a-m) \Phi_{\free}(T)}}_{\to 1/(2\alpha)}
\biggr) \Biggr).
\]
It entails the asymptotic lower bound
\[
\liminf_{T \to +\infty} \, \frac{\E_{\nup}\big[N_a(T)\big]}{T/\Phi_{\free}(T)} \geq
\frac{\alpha}{2 \ln 2} \bigl( \ln(1/\alpha) - 2 \ln 2 \bigr) = \frac{1}{16}
\]
for the choice $\alpha = 1/8$.
The claimed result follows by adding these lower bounds for each suboptimal arm~$a$,
with a factor~$\Delta_a$, following the rewriting~\eqref{eq:defregrstoch}
of the regret.

\begin{remark}
\label{rk:th1OKDm+}
The proof above only exploits the fact that the upper end $M$
of the range is unknown: the alternative problems lie in
$\cD_{m,M'}$ for some $M'$ that can be arbitrarily large.
Yet, by definition of adaptation to the range, the strategy
needs to guarantee $(M'-m) \, \Phi_{\free}(T)$ distribution-free regret bounds
in that case.

We may note that therefore, Theorem~\ref{thm:adaptive_lower_bound}
also holds for the model of bounded distributions with a known lower end $m \in \R$
for the range:
\begin{equation}
\label{eq:Dmodm+}
\Dmod_{m,+} = \bigcup_{\substack{M \in \R: \\ M > m}} \Dmod_{m, M}\,.
\end{equation}
Definitions~\ref{def:free} and~\ref{def:dep} handle
the case of $\Dmod_{-,+}$ but can be adapted in an obvious way
to $\Dmod_{m,+}$ by fixing $m$, by having the strategy know $m$,
and requiring the bounds to hold for all $M \in [m,+\infty)$ and all
bandit problems in $\Dmod_{m,M}$, thus leading to the concept of
adaptation to the upper end of the range.

This observation is in line with the folklore knowledge that there is a difference in nature between dealing
with nonnegative payoffs, i.e., gains, or dealing with nonpositive payoffs, i.e., losses,
for regret minimization under bandit monitoring;
see~\citet[Remark~6.5, page 164]{cesa2006prediction} for an early reference
and~\citet{kwon2016gains} for a more complete literature review.
Actually, $0$ plays no special role, the issue is rather
whether one end of the payoff range is known.
\end{remark}

\section{Adaptation to Range Based on AdaHedge: The AHB Strategy}
\label{sec:AHB}
\label{sec:mainadahedge}

When the range of payoffs is known, \citet{AuCBFrSc02}
achieve a distribution-free regret bound of order $\sqrt{KT \ln K}$ with
exponential weights---the Hedge strategy---on estimated payoffs and with extra-exploration, i.e.,
by mixing exponential weights with the uniform distribution over arms.
Actually, it is folklore knowledge that the extra-exploration used in this case is unnecessary
(see, among others, \citealp{Sto05}).
To deal with the case of an unknown payoff range, we consider a self-tuned version of Hedge called AdaHedge
(\citealp{de2014follow}, see also an earlier work by \citealp{cesa2007improved}) and do add extra-exploration.
Just as \citet{AuCBFrSc02}, we will actually obtain regret guarantees for oblivious adversarial bandits,
not only distribution-free regret bounds for stochastic bandits. We therefore introduce now
the setting of oblivious adversarial bandits and define adaptation to the range in that case.

\subsection{Oblivious Adversarial Bandits}
\label{sec:setting-oblivious}
\label{sec:setting-equivalence}

In the setting of fully oblivious adversarial bandits (see \citealp{cesa2006prediction,audibert2009Minimax}),
a range $[m,M]$ is set by the environment, where $m,M$ are real numbers, not necessarily nonnegative.
The player is unaware of $[m,M]$ and will remain so.
The environment also picks beforehand a sequence $y_1,y_2,\ldots$ of reward vectors in $[m,M]^K$. We denote by
$y_t = (y_{t,a})_{a \in [K]}$ the components of these vectors. The player will observe a component of each of
these reward vectors in a sequential fashion, as follows.
\changes{Auxiliary randomizations $U_0,U_1,\ldots$ i.i.d.\ according to a uniform distribution over $[0,1]$
are available.} At each round $t \geq 1$, the player picks an arm $A_t \in [K]$, possibly at random
\changes{(thanks to $U_{t-1}$)} according to a probability distribution $p_t = (p_{t,a})_{a \in [K]}$,
and then receives and observes $y_{t,A_t}$.

More formally, a strategy of the player is a sequence of mappings from the observations to the action set,
$(U_0, \, y_{1, A_1}, U_1, \, \ldots, \, y_{t-1, A_{t-1}}, U_{t-1}) \mapsto A_t$. The strategy does not rely on $m$ nor $M$.

At each given time $T \geq 1$, denoting by $y_{1:T} = (y_1,\ldots,y_T)$ the reward vectors,
we measure the performance of a strategy through its expected regret:
\begin{equation}
\label{eq:regret-oblivious}
R_T(y_{1:T}) =  \max_{a \in [K]} \sum_{t = 1}^T y_{t,a} - \E \! \left[ \sum_{t = 1}^T y_{t,A_t} \right],
\end{equation}
where, as rewards are fixed beforehand, all randomness lies in the choice of the arms $A_t$ only,
i.e., where  the expectation is only over the choice of the arms $A_t$.

The counterpart of Definition~\ref{def:free} in this setting
is stated next.

\begin{definition}[Scale-free adversarial regret bounds]
\label{def:adv}
A strategy for oblivious adversarial bandits is
adaptive to the unknown range of payoffs with a scale-free adversarial
regret bound $\Phi_{\adv} : \N \to [0,+\infty)$
if for all real numbers $m < M$,
the strategy ensures, without the knowledge of $m$ and $M$:
\[
\forall y_1,y_2,\ldots \ \mbox{\rm in } [m,M]^K, \ \ \forall T \geq 1,
\qquad R_T(y_{1:T}) \leq (M-m) \, \Phi_{\adv}(T)\,.
\]
\end{definition}

\paragraph{Conversion of upper/lower bounds from one setting to the other.}
We recall that when applying Doob's optional skipping in Section~\ref{sec:setting-stochastic},
for each arm $a$, we denoted by $(Y_{t,a})_{t \geq 1}$ an i.i.d.\ sequence of rewards drawn beforehand, independently at random,
according to the distribution $\nu_a$ associated with that arm.
By the tower rule for the right-most equality below, we note that for all $m < M$ and
for all $\nup$ in $\Dmod_{m, M}$,
\begin{multline*}
R_T(\nup) = \max_{a \in [K]} \E \! \left[\sum_{t = 1}^{T} Y_{t,a} \right] - \E \! \left[\sum_{t = 1}^{T} Y_{t,A_t} \right]
\leq \E \! \left[ \max_{a \in [K]} \sum_{t = 1}^{T} Y_{t,a} - \sum_{t = 1}^{T} Y_{t,A_t} \right]
= \E \bigl[ R_T(Y_{1:T}) \bigr] \\
\leq \sup_{y_{1:T} \mbox{\tiny ~in } [m,M]^K} R_T(y_{1:T}) \,.
\end{multline*}
In particular, lower bounds on the regret for stochastic bandits are also lower bounds on the
regret for oblivious adversarial bandits, and strategies designed for oblivious adversarial bandits
obtain the same distribution-free regret bounds for stochastic bandits when the individual payoffs $y_{t,A_t}$
in their definition are replaced with the stochastic payoffs $Y_{t,A_t}$.

\subsection{The AHB Strategy}
\label{sec:algoAHB}

We state our main strategy, AHB---which stands for AdaHedge for $K$--armed Bandits, with extra-exploration---,
in the setting of oblivious adversarial bandits, see Algorithm~\ref{alg:exp3_adahedge}.
In a setting of stochastic bandits, it suffices to replace therein
$y_{t,A_t}$ with $Y_{t,A_t}$.
The AHB strategy relies on a payoff estimation scheme, which we discuss now.

\begin{figure}[t]
\renewcommand\footnoterule{}
\begin{algorithm}[H]
\stretchasneed
\begin{algorithmic}[1]
\caption{\label{algo:adahedge} AHB: AdaHedge for $K$--armed Bandits, with extra-exploration}
\label{alg:exp3_adahedge}
\STATE \textbf{Input:} a sequence $(\gamma_t)_{t \geq 1}$ in $[0,1]$ of extra-exploration rates;
a payoff estimation scheme, e.g.,~\eqref{eq:estscheme}\hspace{-1cm} \
\FOR{rounds $t = 1, \dots, K$}
\STATE Draw arm $A_t = t$
\STATE Get and observe the payoff $y_{t, t}$
\ENDFOR
\STATE \textbf{AdaHedge initialization:} $\eta_{K+1} = +\infty$ and $q_{K+1} = (1/K,\ldots,1/K) \eqdef \mathbf{1}/K$
\FOR{rounds $t = K+1,  \dots$}
\STATE Define $p_t$ by mixing $q_t$ with the uniform distribution according to $p_t = (1 - \gamma_t) q_t + \gamma_t \mathbf{1}/K$
\STATE Draw an arm $A_t \sim p_t$, i.e., independently at random according to the distribution $p_t$
\STATE Get and observe the payoff $y_{t, A_t}$
\STATE Compute estimates $\wh{y}_{t,a}$ of all payoffs with the payoff estimation scheme considered, e.g.,~\eqref{eq:estscheme}
\STATE Compute the mixability gap $\delta_t \geq 0$ based on the distribution $q_t$ and on these estimates:
\[
\underbrace{\delta_t =
- \sum_{a=1}^K q_{t,a}\,\wh{y}_{t,a} + \frac{1}{\eta_t} \log \Biggl( \sum_{a =1}^K q_{t, a} \e^{\eta_t \wh{y}_{t,a}} \Biggr)}_{\mbox{when $\eta_t \leq +\infty$
or $\eta_t = +\infty$}},
\qquad \mbox{i.e.,} \qquad
\underbrace{\delta_t = -\sum_{a=1}^K q_{t,a}\,\wh{y}_{t,a} + \max_{a \in [K]} \wh{y}_{t,a}}_{\mbox{when $\eta_t = +\infty$}}
\vspace{-.4cm}
\]
\STATE Compute the learning rate $\displaystyle{\eta_{t+1} = \Biggl( \sum_{s = K +1}^t \delta_s \Biggr)^{-1}} \ln K$
\STATE Define $q_{t+1}$ component-wise as \vspace{-.3cm}
\[
~ \hspace{2.5cm} q_{t+1,a} = \exp \! \left( \eta_{t+1} \sum_{s = K +1}^{t} \wh{y}_{a,s} \right)
\Bigg/ \sum_{k=1}^K \exp \! \left( \eta_{t+1} \sum_{s = K +1}^{t} \wh{y}_{k,s} \right) \vspace{-.35cm}
\]
\ENDFOR
\end{algorithmic}
\end{algorithm}
\end{figure}

In Algorithm~\ref{alg:exp3_adahedge},
some initial exploration lasting $K$ rounds is used to get a rough idea of the location of the payoffs
and to center the estimates used at an appropriate location.
Following~\citet{AuCBFrSc02}, we consider, for all rounds $t \geq K+1$
and arms $a \in [K]$,
\begin{equation}
\label{eq:estscheme}
\wh{y}_{t,a} = \frac{y_{t,A_t} -C }{p_{t,a}} \1{A_t = a} + C
\qquad \mbox{where} \qquad C \eqdef \frac{1}{K} \sum_{s=1}^K y_{s,s}\,.
\end{equation}
Note that all $p_{t,a} > 0$ for Algorithm~\ref{alg:exp3_adahedge} due to the use
of exponential weights.
As proved by~\citet{AuCBFrSc02}, the estimates $\wh{y}_{t,a}$ are conditionally unbiased.
Indeed, the distributions $q_{t}$ and $p_{t}$, as well as the constant $C$, are measurable functions of
the information $H_{t-1} = (U_0, \, y_{1, A_1}, U_1, \, \ldots, U_{t-2}, \, y_{t-1, A_{t-1}})$
available at the beginning of round~$t \geq K+1$,
and the arm $A_t$ is drawn independently at random according to $p_t$ based
on an auxiliary randomization denoted by $U_{t-1}$. Therefore, given that the payoffs are oblivious,
the conditional expectation of $\wh{y}_{t,a}$ with respect to $H_{t-1}$
amounts to integrating over the randomness given by the random draw $A_t \sim p_t$:
for $t \geq K+1$,
\begin{equation}
\label{eq:condexpwhy}
\E \bigl[ \wh{y}_{t,a} \,\big|\, H_{t-1} \bigr] =
\frac{y_{t,a} - C}{p_{t,a}} \,\, \mathbb{P} \bigl( A_t = a \,\big|\, H_{t-1} \bigr) + C
= \frac{y_{t,a}- C}{p_{t,a}} \, p_{t,a}  + C = y_{t,a}\,.
\end{equation}
These estimators are bounded:
assuming that all $y_{t,a}$, thus also $C$,
belong to the range $[m,M]$, and
given that the distributions $p_t$ were obtained by
a mixing with the uniform distribution, with weight $\gamma_t$, we have
$p_{t,a} \geq \gamma_t/K$, and therefore,
\begin{equation}
\label{eq:AH-bound}
\forall t \geq K+1, \ \ \forall a \in [K], \qquad\quad
\bigl| \wh{y}_{t,a} - C \bigr| \leq
\frac{|y_{t,a}- C|}{p_{t,a}} \leq
\frac{M-m}{\gamma_t/K}\,.
\end{equation}

\begin{remark}
Algorithm~\ref{algo:adahedge}
is invariant by affine changes, i.e., translations by real numbers and/or multiplications by positive factors, of the payoffs,
given that AdaHedge (see~\citealp[Theorem~16]{de2014follow}) and the payoff estimation scheme~\eqref{eq:estscheme}
are so. This is key for adaptation to the range.

This invariance is achieved, when ignoring the range $[m,M]$, thanks to a value $C \in [m,M]$.
\changes{Here}, we chose to have $K$ rounds of exploration in Algorithm~\ref{alg:exp3_adahedge}
and let $C$ equal the average of the payoffs achieved. However, it would of course have
been sufficient to pick one arm at random, observe a single reward $y_{1,A_1}$ and let $C = y_{1,A_1}$.
\end{remark}

\subsection{Regret Analysis, Part 1: Scale-Free Adversarial Regret Bound}
\label{sec:adahedge-distrfree}

\begin{theorem}
\label{th:adv_talpha_bound}
AdaHedge for $K$--armed bandits (Algorithm~\ref{alg:exp3_adahedge})
with a non-increasing extra-exploration sequence $(\gamma_t)_{t \geq 1}$ smaller than $1/2$
and the estimation scheme given by~\eqref{eq:estscheme}
ensures that for all bounded ranges $[m,M]$,
for all oblivious individual sequences $y_1,y_2,\ldots$ in $[m,M]^K$,
for all $T \geq 1$,
\[
R_T(y_{1:T}) \leq 3 (M-m) \, \sqrt{KT \ln K}
+ 5 (M-m)  \frac{K \ln K}{\gamma_T} + (M-m) \sum_{t = K+1}^T \gamma_t\,.
\]
In particular, given a parameter $\alpha \in (0, 1)$,
the extra-exploration $\smash{\gamma_t = \min \Bigl\{1/2, \, \sqrt{5(1-\alpha)K \log K} \big/ t^{\alpha} \Bigr\}}$
leads to the scale-free adversarial regret bound
\begin{equation}
\label{eq:defphiadvth2}
\Phi_{\adv}(T) =
\biggl( 3 + \frac{5}{\sqrt{1-\alpha}} \biggr) (M-m)\sqrt{K \ln K} \,\, T^{\max\{\alpha,1-\alpha\}} + 10 (M-m)K\ln K\,.
\end{equation}
For $\alpha = 1/2$, the bound reads $\Phi_{\adv}(T) = 7(M-m) \sqrt{TK \ln K} + 10(M-m) K \log K$.
\end{theorem}

This value $\alpha = 1/2$ is the best one to consider
if one is only interested in a distribution-free bound---i.e., if one is not interested in the distribution-dependent rates for the regret.
The proof of Theorem~\ref{th:adv_talpha_bound} is detailed
in Appendix~\ref{sec:proof:thUBadv} but we sketch its proof below.

\begin{remark}
\label{rk:lnN}
We strongly suspect that the $\sqrt{\log K}$ factor in the bound of Theorem~\ref{th:adv_talpha_bound} is superfluous.
In the case of a known range, the MOSS algorithm is known to be minimax optimal with a regret bound of order $\sqrt{KT}$.
One idea could thus be to use a MOSS-type index, together with a Bernstein-type upper confidence bound to account for the unknown variance and range. A final ingredient would be to add initial extra-exploration, pulling every arm $\sqrt{T / K}$ times before running the standard phase of the algorithm;
on a technical level, this automatically makes the sub-Poissonian term in Bernstein's inequality tractable.
We have not managed yet to fill in the technical details in order to prove this, although we believe a
variant of these ideas would get rid of the logarithmic factor.
In contrast, the algorithm discussed here, based on AdaHedge, enjoys a simple distribution-free analysis---as sketched below---,
as well as a distribution-dependent analysis (see Section~\ref{sec:adahedge-distrdep}), unlike an algorithm based on MOSS-type indices.

Another promising approach would be to use the Tsallis-INF algorithm introduced by \citet{audibert2009Minimax} and further studied by \citet{zimmert2018optimal}, which achieves a $(M-m)\sqrt{KT}$ adversarial regret bound when $M$ and $m$ are known. Unfortunately, current analyses of the algorithm rely crucially on the non-positivity of the reward estimates, or, equivalently on the knowledge of an upper bound on the rewards. \citet{zimmert2019connections} relax this requirement, but not enough for the relaxed version to be applied to our case.
However, when $M$ is known and $m$ is unknown, i.e., only adaptation to $m$ is needed, the reward estimates can be made non-positive by taking $C = M$ in the estimation scheme~\eqref{eq:estscheme}, and our techniques may be extended to
show that Tsallis-INF indeed enjoys an adversarial regret bound of order $(M-m)\sqrt{KT}$ in this case.
Details may be found in \secrefsuppl{Theorem~23 in~Appendix~F}{Theorem~\ref{thm:tsallis} in Appendix~\ref{sec:TsallisKTpossible}}.
\end{remark}

\bpc{sketch}
A direct application of the AdaHedge regret bound (Lemma~3 and Theorem~6 of~\citealp{de2014follow}),
bounding the variance terms of the form $\E\bigl[(X - \E[X])^2\bigr]$
by $\E\bigl[(X - C)^2\bigr]$, ensures that
\[
\max_{k \in [K]} \sum_{t=K+1}^T \wh{y}_{t,k} -
\sum_{\substack{t \geq K+1 \\ a \in [K]}}^T q_{t,a} \, \wh{y}_{t,a} \leq
2 \sqrt{\sum_{{\substack{t \geq K+1 \\ a \in [K]}}} q_{t,a} \bigl( \wh{y}_{t,a} - C \bigr)^2 \ln K} +
\frac{M-m}{\gamma_T/K} \left( 2 + \frac{4}{3}\ln K \right)\,.
\]
We take expectations, use the definition of the $p_t$ in terms of the $q_t$
in the left-hand side, and apply Jensen's inequality in the right-hand side
to get
\begin{multline*}
\E \Biggl[
\max_{k \in [K]} \sum_{t=K +1}^T \wh{y}_{t,k} -
\sum_{t=K+1}^T \overbrace{\sum_{a=1}^K p_{t,a} \, \wh{y}_{t,a}}^{= y_{t,A_t}}
+ \sum_{t=K+1}^T \gamma_t \overbrace{\sum_{a=1}^K (1/K - q_{t,a}) \, \wh{y}_{t,a}}^{\E[...] \in [m-M,M-m]} \Biggr] \\
\leq
2 \sqrt{\sum_{t=K+1}^T \sum_{a=1}^K \E \Bigl[ q_{t,a} \bigl( \wh{y}_{t,a} - C \bigr)^2 \Bigr] \ln K} +
\frac{M-m}{\gamma_T/K} \left( 2 + \frac{4}{3}\ln K \right).
\end{multline*}
Since $p_{t,a} \geq (1-\gamma_t) q_{t,a}$ with $\gamma_t \leq 1/2$ by assumption on the extra-exploration rate,
we have the bound $q_{t,a} \leq 2 p_{t,a}$.
Together with standard calculations similar to~\eqref{eq:condexpwhy}, we have
\[
\E \Bigl[ q_{t,a} \bigl( \wh{y}_{t,a} - C \bigr)^2 \Bigr]
\leq 2 \, \E \Bigl[ p_{t,a} (\wh{y}_{t,a} - \C)^2 \,\Big|\, H_{t-1} \Bigr]
= 2 \, \E \!\left[ \frac{(y_{t,A_t} - \C)^2}{p_{t,a}} \1{A_t = a} \right]
= 2 \underbrace{(y_{t,a} - \C)^2}_{\leq (M-m)^2}\,.
\]
The proof of the first regret bound of the theorem
is concluded by collecting all bounds and by taking care of the first $K$ rounds.
The second regret bound then follows from straightforward calculations.
\end{proof}

\subsection{Regret Analysis, Part 2: Distribution-Dependent Rates for Adaptation}
\label{sec:adahedge-distrdep}

Given the conversion explained in Section~\ref{sec:setting-equivalence},
Algorithm~\ref{alg:exp3_adahedge} tuned as in Corollary~\ref{th:adv_talpha_bound}
for $\alpha \in [1/2, 1)$ also enjoys the scale-free distribution-free regret bound
$\Phi^{\AH}_{\free}(T) = \Phi^{\AH}_{\adv}(T)$ of order $T^{\alpha}$.
The theorem below entails that AHB is adaptive to the unknown range with a
distribution-dependent regret rate $T / \Phi^{\AH}_{\free(T)}$ of order $T^{1-\alpha}$
that is optimal given the lower bound stated by Theorem~\ref{thm:adaptive_lower_bound}.

\begin{theorem}\label{thm:distrib-dependent}
Consider AHB (Algorithm~\ref{alg:exp3_adahedge}) tuned
with some $\alpha \in [1/2, 1)$ as in the
second part of Theorem~\ref{th:adv_talpha_bound}.
For all distributions $\nu_1,\ldots,\nu_K$ in $\cD_{-,+}$,
\begin{equation}
\limsup_{T \to \infty} \frac{R_T(\nup)}{T / \Phi^{\AH}_{\free(T)}}
\leq  \frac{12 \ln K}{1-\alpha} \sum_{a = 1}^K \Delta_a \,.
\end{equation}
\end{theorem}

The proof is provided in Appendix~\ref{sec:bernsrem}.
It follows quite closely that of Theorem~3 in \citet{seldin2017improved}, where the
authors study a variant of the Exp3 algorithm of~\citet{AuCBFrSc02} for stochastic rewards. It consists,
in our setting, in showing that the number of times the algorithm chooses suboptimal arms is almost only determined
by the extra-exploration. Our proof is simpler as we aim for cruder bounds. The main technical difference
and issue to solve lies in controlling the learning rates $\eta_t$, which heavily depend on data in our case.

\section{Numerical Illustrations}
\label{sec:numexpe}

We provide some numerical experiments on synthetic data to illustrate the qualitative behavior of
some popular algorithms like UCB strategies when they are incorrectly tuned, as opposed to
strategies that are less sensitive to ignoring the range or to the AHB strategy which adapts to it.
These experiments are only of an illustrative nature.

\begin{figure}[t]
	\includegraphics[width=1\textwidth]{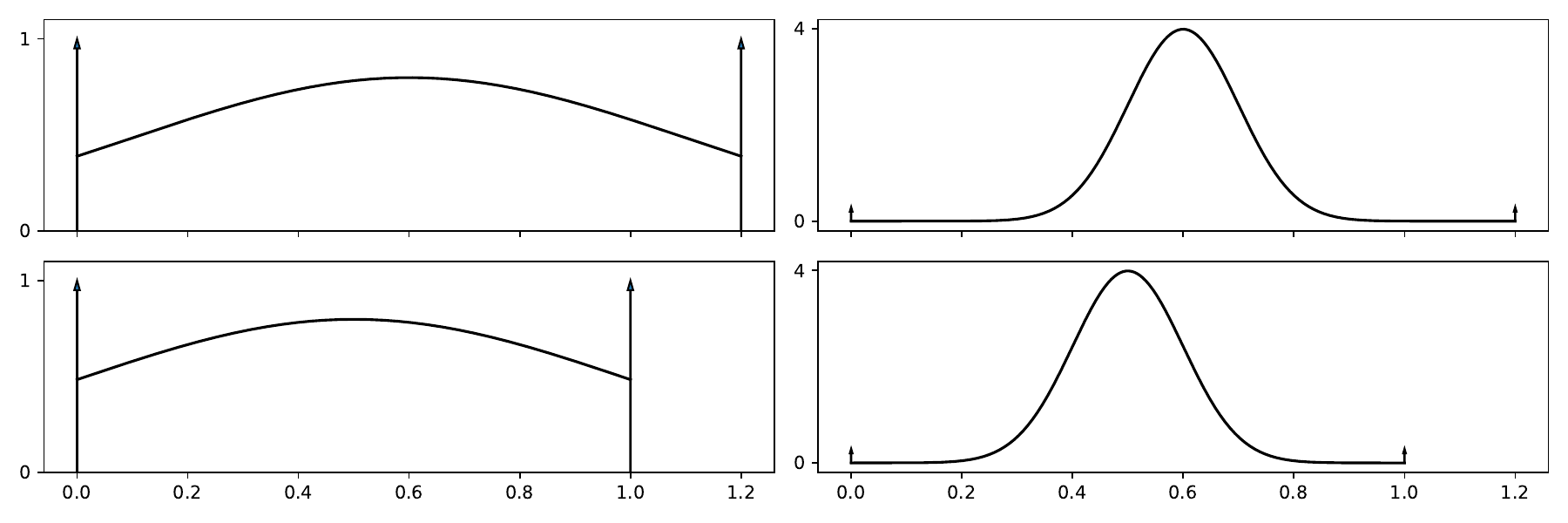}
	\caption{\label{fig:densities}
		Probability density functions of the reward distributions with respect to the sum of the Lebesgue measure and Dirac masses at $0$, $1$, and $1.2$.
        Left pictures: high-variance case; right pictures: low-variance case.
        Top pictures: first arm (optimal arm); bottom pictures: other arms.
		Arrows represent atoms and their lengths are only illustrative.
	}
\end{figure}

\paragraph{Bandit problems considered and UCB strategies.}
We consider stochastic bandit problems $\nup^{(\alpha)} = (\nu_a^{(\alpha)})_{a \in [K]}$
indexed by a scale parameter $\alpha \in \{ 0.01 , \, 1,    \, 100  \}$.
We take $K = 10$ arms, each arm~$a$ being
associated with a rectified Gaussian distribution. Precisely, the distribution $\nu_a^{(\alpha)}$ is the distribution of the variable
\begin{equation*}
	X_a^\alpha = \left\{
	\begin{split}
		\alpha \, \max\big\{ 0, \min \{Y, \,\,  1.2\} \big\} \quad \text{with} \quad  Y \sim \mathcal N(0.6\, , \, V)  \quad \text{if } a = 1, \\
        \alpha \, \max\big\{ 0, \min \{Y, \,\,  1 \} \big\} \quad \text{with} \quad  Y \sim \mathcal N(0.5\, , \, V) \quad \text{if } a \neq 1,
	\end{split} \right.
\end{equation*}
so that all distributions are commonly supported on $[m,M] = [0, \,\, 1.2 \, \alpha]$,
with arm~1 being the unique optimal arm. We will consider two values for $V$, namely $V = 0.01$ (low-variance case)
and $V = 0.25$ (high-variance case). See Figure~\ref{fig:densities} for a plot of the corresponding probability
density functions.

We denote by
$\mu^{(\alpha)}_1 = 0.6 \, \alpha$ and $\mu^{(\alpha)}_a = 0.5 \, \alpha$ if $a \neq 1$
the means associated with the distributions $\nu_1^{(\alpha)}$ and $\nu_a^{(\alpha)}$, respectively.
The gaps therefore equal $\Delta^{(\alpha)}_a = 0.1 \, \alpha$ for $a \geq 2$.

The main algorithm of interest is, of course, the AHB strategy with extra-exploration (Algorithm~\ref{algo:adahedge}),
which we tune as stated in Theorem~\ref{th:adv_talpha_bound} with parameter~$1/2$.
We now present the competitors.

\paragraph{UCB strategies at different scales.}
We consider instances of UCB (\citealp{UCB-orig}) using indices of the form
\[
\what \mu_a(t) +  \sqrt{\frac{8 \sigma^2  \log T}{N_a(t)}}\,,
\]
where $N_a(t)$ is the number of times arm $a$ was pulled up to round $t$,
and where $\what \mu_a(t)$ denotes the empirical average of payoffs obtained for arm $a$.
We hesitated between setting $\sigma^2$ based on the range $M-m = 1.2\alpha$,
namely, $\sigma^2 = (M-m)^2/4 = (1.2 \alpha)^2$, or based on a sub-Gaussian parameter,
which would be smaller. As distributions $\nu_a^{(\alpha)}$ are rectified Gaussians, it is not
immediately clear whether they are sub-Gaussian, but we considered despite all the choice $\sigma^2 = V$.
It turns out that this second choice outperformed the first one, which is why, in the rest of the study,
we consider the following three instances of UCB:
\[
\what \mu_a(t) +  s \sqrt{\frac{8 V \log T}{N_a(t)}}\,, \qquad \mbox{where} \qquad s \in \{0.01, \, 1, \, 100\}\,.
\]
When the scale parameter $\alpha$ is known, we would take $s = \alpha$.

\paragraph{Range-estimating UCB.}
We also study a version of UCB estimating the range, namely, using indices
\[
\what \mu_a(t) + \hat{r}_t \sqrt{\frac{2\log T}{N_a(t)}}\,,
\qquad \mbox{where} \qquad
\hat{r}_t = \max_{s \leq t} Y_{A_s,s} - \min_{s \leq t} Y_{A_s,s}
\]
estimates the range $M-m$.
We were unable to provide theoretical guarantees that match our lower bounds, and this algorithm does not perform particularly well in practice
as we will discuss below.

\paragraph{$\eps$--greedy.}
Finally, we also consider the $\epsilon$--greedy strategy, which, at round $t \geq K+1$, picks with probability $1 - \eps_t$ the arm with the best empirical mean, and otherwise, selects an arm uniformly at random. Following \cite{UCB-orig}, we used the tuning
\begin{equation*}
	\eps_t = \min \left\{ 1, \,\, \frac{ 5 K}{d^2 t} \right\}
	\quad \text{with} \quad
	 d = 1 / 12 \, .
\end{equation*}
Indeed, \cite{UCB-orig} exhibit theoretical guarantees for distributions over $[0,1]$ in the case where
$d$ is smaller than or equal to the smallest gap.
When rescaled on $[0,1]$, the smallest gap equals $0.1 \alpha (M-m)/ (1.2 \alpha) = 1/12$ in our setting;
this explains our choice $d = 1/12$, but note that the $\epsilon$--greedy strategy defined above relies on some extra knowledge encompassed in the choice $d = 1/12$, compared to the completely agnostic
AHB strategy.
Interestingly, for any fixed-in-advance sequence of $\eps_t$,
the $\eps$--greedy strategy is scale-free. Of course, its strong downside is that a proper
tuning of the $\eps_t$ requires knowledge of a scaled lower bound on the gaps.

\begin{figure}[t!h!]
    \center
	\includegraphics[width=0.95\textwidth]{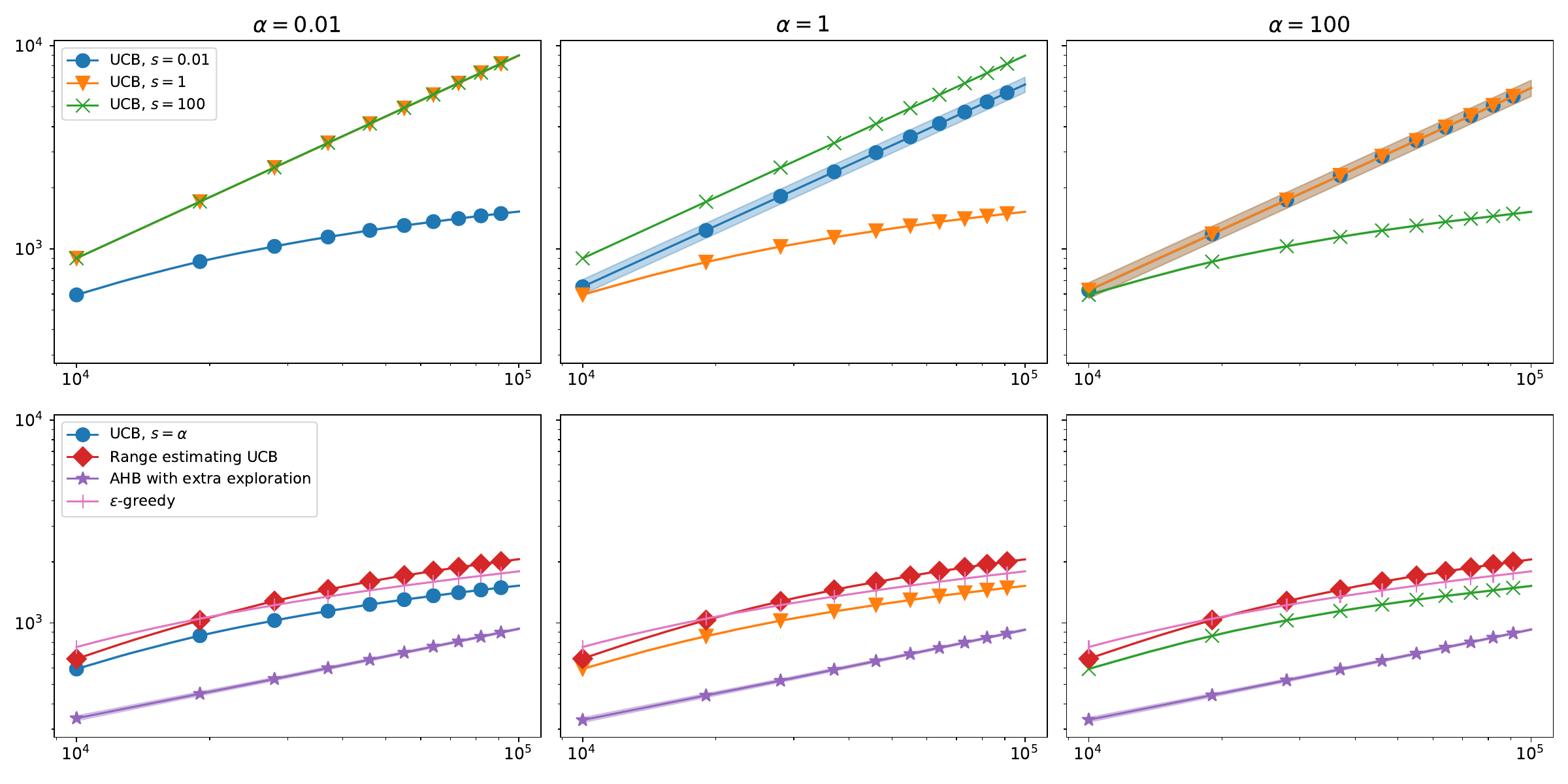}
	\caption{\label{fig:exp}
		Comparison of the (estimated) regrets of various strategies
		over bandit problems $\nup^{(\alpha)}$ in the high variance case, where $\alpha$ ranges in
		$\{ 0.01 , \, 1 , \, 100  \}$ and $V = 0.25$.
		Each algorithm was run $N = 100 $ times on every problem for $T = 100,\!000$ time steps.
		Solid lines report the values of the estimated regrets, while
		shaded areas correspond to $\pm 2$ standard errors of the estimates.
	}
\end{figure}

\begin{figure}[t!h!]
    \center
	\includegraphics[width=0.95\textwidth]{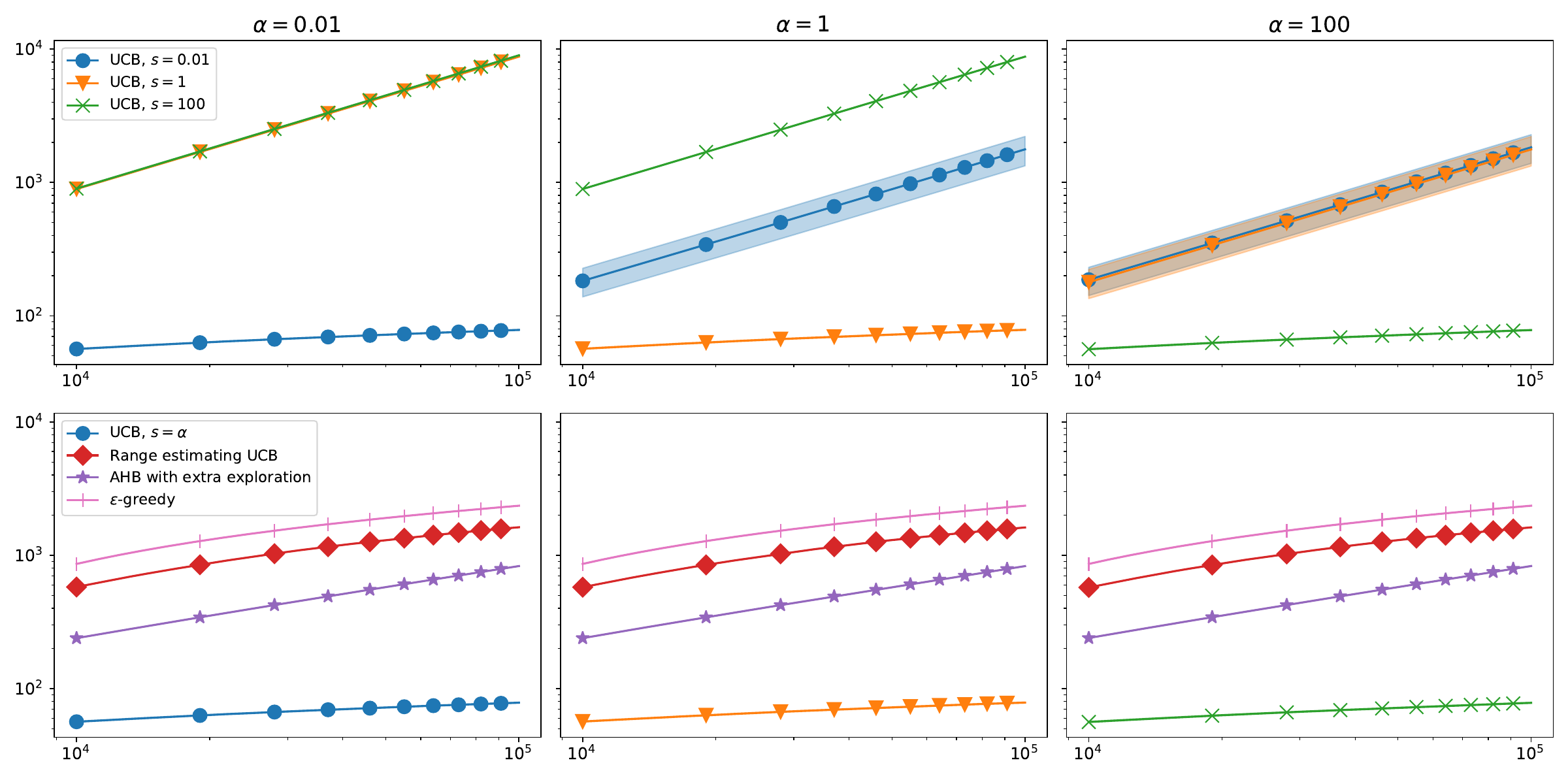}
	\caption{\label{fig:exp2}
		Same legend as for Figure~\ref{fig:exp}, but in the low variance case.
	}
\end{figure}

\paragraph{Experimental setting.}
Each algorithm is run $N = 100$ times, on a time horizon $T = 100,\!000$.
We plot estimates of the rescaled regret $R_T(\nup^{(\alpha)}) / \alpha$ to have a meaningful comparison between the bandit
problems.
These estimates are constructed as follows.
We index the arms picked in the $n$--th run by an additional subscript $n$,
so that $A_{T,n}$ refers to the arm picked by some strategy
at time $t$ in the $n$--th run.
The expected regret of a given strategy can be rewritten as
\[
R_T(\nup^\alpha) = T \max_{a \in [K]} \mu^{(\alpha)}_a - \E \! \left[\sum_{t = 1}^{T} \mu^{(\alpha)}_{A_t} \right]
= T \times ( 0.6 \, \alpha) - \E \! \left[\sum_{t = 1}^{T} \mu^{(\alpha)}_{A_t} \right]
\]
and is estimated by
\[
\what R_T(\alpha) = \frac{1}{N} \sum_{n=1}^N \what R_T(\alpha, n)
\qquad \mbox{where} \qquad
\what R_T(\alpha, n)= T \times ( 0.6 \, \alpha) - \sum_{t = 1}^T \mu^{(\alpha)}_{A_{t, n}}\,.
\]
On Figures~\ref{fig:exp} and~\ref{fig:exp2} we plot the estimates ${\what R_T(\alpha)}/\alpha$
of the rescaled regret as solid lines.
The shaded areas correspond to $\pm 2$ standard errors of the
sequences $\big(\what R_T(\alpha, n) / \alpha \big)_{n \in [N]}$.

\paragraph{Discussion of the results.}
An initial observation is that, as expected, the performance of AHB,
the range-estimating UCB, and $\eps$-greedy
is unaffected by the scale of the problems
(see the second lines of Figures~\ref{fig:exp} and~\ref{fig:exp2}).
It turns out that out of these three algorithm, AHB performs best.

A second observation is that the performance of UCB depends dramatically on the value of the parameter $s$. UCB performs like follow-the-leader when $s$ is too small, and like random play when $s$ is too large; both of these strategies suffer linear regret and
UCB incorrectly scaled also does so (see the first lines of Figures~\ref{fig:exp} and~\ref{fig:exp2}).

It remains to compare AHB to UCB tuned with the correct scale:
the ranking between the two depends on the value of $V$,
with AHB outperforming UCB tuned with the correct scale in the high-variance case
and vice versa in the low variance case.

Our last observation is that in the low-variance case,
the range-estimating version of UCB is far off from UCB tuned with the correct scale.
This is because of the large difference between the sub-Gaussian parameter and its upper bound
given by the squared half-range, which the range-estimating version of UCB is targeting.

\appendix
\section{More on Scale-Free Distribution-Dependent Regret Bounds \\ ~~~~~~~~~~~~~~~~~~Considered in Isolation}
\label{sec:LB-nolnT}

This section details the claims of Section~\ref{sec:distrdep-isolation}:
no strategy may be adaptive to the range and achieve $\Phi_{\dep} = \ln T$
(Section~\ref{sec:A1}) but we may construct a strategy adaptive to the range and achieving
$\Phi_{\dep} \gg \ln T$ (Section~\ref{sec:A2}). Before we do so,
we provide a reminder on a general, and optimal, distribution-dependent regret lower bound for $K$--armed stochastic bandits
(Section~\ref{sec:LRBK}).

\subsection{Reminder of a General Regret Lower Bound for $K$--Armed Bandits}
\label{sec:LRBK}

This section considers some general model $\cD$.
It also rules out poor strategies by restricting its attention to
so-called consistent strategies---according to the terminology
introduced by \citet{lai1985asymptotically},
while \citealp{BuKa96} rather speak of uniformly fast convergent strategies.

\begin{definition}
\label{def:UCF}
A strategy is consistent on a model $\Dmod$ if for all bandit problems $\nup$ in $\Dmod$,
it achieves a subpolynomial regret bound, that is, $R_T(\nup)/T^\alpha \to 0$ for \changes{all $\alpha \in (0,1]$}.
\end{definition}

A lower bound on the distribution-dependent rates that such a strategy may achieve is provided by a general, and optimal,
result of \citet{lai1985asymptotically} and \citet{BuKa96}; see also its rederivation by~\citet{garivier2018Explore}.
It involves a quantity defined as an infimum of Kullback-Leibler divergences:
we recall that for two probability distributions $\nu,\nu'$ defined on the same probability space $(\Omega,\mathcal{F})$,
\[
\KL(\nu,\nu') = \left\{\begin{array}{ll}
	\displaystyle \bigintsss_{\Omega} \ln\!\left(\frac{\d \nu}{\d \nu'}\right)\! \d \nu & \textrm{if $\nu \ll \nu'$}, \smallskip \\
	+ \infty & \textrm{otherwise},
\end{array} \right.
\]
where $\nu \ll \nu'$ means that $\nu$ is absolutely continuous with respect to $\nu'$ and
$\d\nu/\d\nu'$ then denotes the Radon-Nikodym derivative.
Now, for any probability distribution $\nu$,
any real number $x$, and any model $\Dmod$, we define
\[
\Kinf(\nu,x,\Dmod) = \inf\bigl\{ \KL(\nu,\nu') : \nu' \in \Dmod \mbox{ and } \Ed(\nu') > x \bigr\}\,,
\]
where by convention, the infimum of an empty set equals $+\infty$ and where we denoted by $\Ed(\nu')$
the expectation of $\nu'$.
The quantity $\Kinf(\nu,x,\Dmod)$ can be null. With the usual measure-theoretic conventions,
in particular, $0/0 = 0$, we then have the following lower bound.

\begin{reminder}
\label{rem:LRBK}
For all models $\Dmod$, for all consistent strategies on $\Dmod$,
for all bandit problems $\nup$ in $\Dmod$,
\[
\liminf_{T \to +\infty} \frac{R_T(\nup)}{\ln T} \geq \sum_{a \in [K]} \frac{\Delta_a}{\Kinf(\nu_a,\mu^\star,\Dmod)}\,.
\]
\end{reminder}

\paragraph{The case of a known payoff range $[m,M]$.}
When the payoff range $[m,M]$ is known, i.e., when the model is $\Dmod_{m,M}$,
there exist strategies achieving the lower bound of Reminder~\ref{rem:LRBK},
like the DMED strategy of~\citet{hondatakemura-first,honda15anon-asymptotic}
or the KL--UCB strategy of~\citet{cappe_kullbackleibler_2013} and~\citet{garivier2018kl}.

\paragraph{The case of a known payoff upper bound $M$.}
The DMED strategy of~\citet{honda15anon-asymptotic} actually achieves
the lower bound of Reminder~\ref{rem:LRBK} even for the model
\[
\Dmod_{-,M} = \bigcup_{\substack{m \in \R: \\ m < M}} \Dmod_{m, M}
\]
and for the model $\Dmod_{-\infty,M}$ of all distributions upper bounded by~$M$
but not necessarily lower bounded. This suggests that adaptation to~$M$
is much more difficult than adaptation to~$m$ as far as distribution-dependent regret bounds
are considered, and is in line with Remark~\ref{rk:th1OKDm+}.

That $M$ is more important than $m$ for distribution-dependent
bounds is also reflected in the lower bound of Reminder~\ref{rem:LRBK}:
this lower bound does not depend on whether $\Dmod$ equals some $\Dmod_{m, M}$,
or $\Dmod_{-,M}$, or even $\Dmod_{-\infty,M}$.
We may indeed easily show (\changes{see \secrefsuppl{Appendix~E}{see Appendix~\ref{sec:app-a-e}}}) that given $M \in \R$,
for all $m \leq M$, for all $\nu \in \Dmod_{m, M}$ and all $\mu > \Ed(\nu)$,
\[
\Kinf\big( \nu, \mu, \Dmod_{m, M} \big)
= \Kinf\big( \nu, \mu, \Dmod_{-\infty, M} \big)\,.
\]

\subsection{Adaptation to the Range Impossible at Logarithmic Distribution-Dependent Rate}
\label{sec:A1}

A strategy that would be adaptive to the range with a distribution-dependent
rate $\Phi_{\dep} = \ln$ would, by definition and in particular, be consistent on $\Dmod_{-,+}$.
The following theorem therefore shows, by contradiction,
that no strategy may be adaptive to the range with a distribution-dependent
rate $\Phi_{\dep} = \ln$. A similar phenomenon was discussed by \citet{lattimore2017scale}
in the case of stochastic bandits with Gaussian distributions.

\begin{theorem}
\label{th:adaptmorethanlog}
For all distributions $\nu_a \in \Dmod_{-,+}$ with expectation $\mu_a$,
and all $\mu^\star > \mu_a$, we have
\[
\Kinf(\nu_a,\mu^\star,\Dmod_{-,+}) = 0\,.
\]
As a consequence, all consistent strategies on $\Dmod_{-,+}$ are such that,
for all bandit problems $\nup$ in $\Dmod_{-,+}$ with at least one suboptimal arm~$a$,
\[
\liminf_{T \to +\infty} \frac{R_T(\nup)}{\ln T} = +\infty\,.
\]
\end{theorem}

Interestingly, \citet{cowan2015asymptotically} observe that for the model of uniform distributions over bounded intervals,
the $\Kinf$ is positive, and thus the lower bound of Reminder~\ref{rem:LRBK} does not prevent logarithmic regret
bounds. In fact, they also provide an algorithm enjoying optimal distribution-dependent bounds---thus being, in a sense,
adaptive to the range in that very restricted model.
\medskip

\bp
We denote by $[m,M]$ an interval containing the support of $\nu_a$.
We remind the reader of the model $\Dmod_{m,+}$ defined in~\eqref{eq:Dmodm+},
composed of all bounded distributions with unknown upper end on the range
but known lower end $m$ on the range.
As $\Dmod_{m,+} \subset \Dmod_{-,+}$ and by definition of $\Kinf$,
\[
\Kinf(\nu_a,\mu^\star,\Dmod_{-,+}) \leq \Kinf(\nu_a,\mu^\star,\Dmod_{m,+})\,,
\]
so that it suffices to show that $\Kinf(\nu_a,\mu^\star,\Dmod_{m,+}) = 0$.

We have in particular $\mu_a \geq m$.
We use the same construction as in the proof of Theorem~\ref{thm:adaptive_lower_bound}.
Let $\nu'_\epsilon = (1-\epsilon)\nu_a + \epsilon \delta_{\mu_a + 2 \Delta_a/\epsilon}$
for $\epsilon \in (0,1)$: it is
a bounded probability distribution, with lower end of support larger than $m$, that is,
$\nu'_\epsilon \in \Dmod_{m,+}$. For $\epsilon$ small enough, $\mu_a + 2 \Delta_a/\epsilon$
lies outside of the bounded support of $\nu_a$. In that case, the density of
$\nu_a$ with respect to $\nu'_\epsilon$ is given by $1/(1-\epsilon)$ on the support of $\nu_a$
and $0$ elsewhere, so that
\[
\KL\bigl(\nu_a,\nu'_\epsilon\bigr) = \ln \biggl( \frac{1}{1-\epsilon} \biggr).
\]
Moreover, $\Ed\bigl(\nu'_\epsilon\bigr) = (1-\epsilon)\mu_a + \epsilon\bigl(\mu_a + 2 \Delta_a/\epsilon\bigr)
= \mu_a + 2\Delta_a = \mu^\star + \Delta_a > \mu^\star$. Therefore,
by definition of $\Kinf$ as an infimum,
\[
\Kinf(\nu_a,\mu^\star,\Dmod_{m,+}) \leq \KL\bigl(\nu_a,\nu'_\epsilon\bigr) = \ln \biggl( \frac{1}{1-\epsilon} \biggr).
\]
This upper bound holds for all $\epsilon > 0$ small enough and thus shows that
$\Kinf(\nu_a,\mu^\star,\Dmod_{m,+}) = 0$.

The second part of the theorem follows from
Reminder~\ref{rem:LRBK}, from the existence of an arm $a$ with
$\Delta_a = \mu^\star - \mu_a > 0$, and from the fact that
$\Kinf(\nu_a,\mu^\star,\Dmod_{-,+}) = 0$, as we established above.
\end{proof}

\begin{remark}
\label{rk:m-noln}
Recall that Remark~\ref{rk:th1OKDm+} defined a notion of adaptation to the upper end~$M$ of the payoff range.
The proof above reveals that Theorem~\ref{th:adaptmorethanlog}
holds with all occurrences of $\Dmod_{-,+}$ replaced by
$\Dmod_{m,+}$, for some $m \in \R$.
We may therefore similarly exclude a $\ln T$ distribution-dependent rate for adaptation
to the upper end~$M$ of the payoff range.

This observation is yet another example that the knowledge of the lower end~$m$ of the payoff
range does not critically change the picture, and the difficulty in
ignoring a payoff range lies in ignoring the upper end thereof.
\end{remark}

\subsection{UCB with an Increased Exploration Rate Adapts to the Range}
\label{sec:A2}

The impossibility result implied by Theorem~\ref{th:adaptmorethanlog}
does not prevent distribution-dependent rates for adaptation that are larger than a logarithm.
Let $\varphi$ be a non-decreasing function such that $\varphi(t) \gg \ln t$,
like $\varphi(t) = (\ln t)^2$ or even $\varphi(t) = (\ln t) \ln \ln t$.
\citet[Remark~8]{lattimore2017scale} introduced and studied,
in the case of Gaussian bandits with unknown variances, the following variant of UCB,
which we refer to in this section as UCB with an increased exploration rate $\varphi$:
\[
\what\mu_a(t) + \sqrt{\frac{\varphi(t)}{N_a(t)}}
\qquad \mbox{where} \quad \frac{\varphi(t)}{\ln t} \to +\infty
\quad \mbox{and} \quad \frac{\varphi(t)}{t} \to 0 \,,
\]
and where $\what\mu_a(t)$ denotes the empirical average of payoffs obtained till round $t$
when playing arm~$a$.
The (asymptotic only) analysis of \citet[Remark~8]{lattimore2017scale} relies on the fact that
$\varphi(t) \geq 2(M-m)\ln t$ for $t$ larger than some unknown threshold $T_0$,
and that after $T_0$, the indexes are thus larger than the ones of the original version
of UCB based on the knowledge of $m$ and $M$. This argument
readily extends to the case of sub-Gaussian distributions, where
we recall that a distribution $\nu$ with expectation $\mu$
is $v$--sub-Gaussian, with $v > 0$, if
\[
\forall t \in \R, \qquad \int \e^{t (x-\mu)} \d\nu(x) \leq \e^{v t^2/2}\,.
\]
Hoeffding's lemma proves that distributions over a bounded range $[m,M]$ are
$(M-m)^2/4$--sub-Gaussian.
Based on a slightly different proof than the one of \citet[Remark~8]{lattimore2017scale}, one can prove the following finite-time result---where
we did not aim for tight numerical constants.

\begin{theorem}
\label{th:UCB-incr}
UCB with an increased exploration rate given by a non-decreasing function $\varphi$ ensures that for all $v > 0$,
for all distributions $\nu_1,\ldots,\nu_K$ that are $v$--sub-Gaussian, for all $T \geq K+1$,
\[
R_T(\nup) \leq \underbrace{\sum_{a \in [K]: \Delta_a > 0} \frac{4}{\Delta_a} \varphi(T)}_{\text{\rm main term}} +
\underbrace{\sum_{a \in [K]} 2 \Delta_a \max \! \left\{ \frac{32 v}{\Delta_a^2}, \, 1 \right\}
\left( 1 + \sum_{t=K}^{T-1} \e^{- \varphi(t) / (2v)} \right)}_{\text{\rm smaller-order term: typically, a } \mathcal{O}(1)}
\]
Whenever $\varphi \gg \ln$,
this strategy is therefore adaptive to the unknown range of payoffs with a distribution-dependent
rate~$\Phi_{\dep} = \varphi$.
\end{theorem}

The second part of the statement follows from the claimed bound given that
$\varphi \gg \ln$ entails $\varphi(t) \geq 4 v \log t$ for $t$ large enough,
and therefore, $\e^{-\varphi(t) / (2v)} \leq 1 / t^2$. As a consequence, the sum tagged as
smaller-order term in the bound is finite.
Possible such choices are $\varphi : t \mapsto (\log t)^2$, or even $\varphi : t \mapsto (\log t) (\log \log t)$.

However, as already mentioned in \citet{lattimore2017scale}, as the distribution-dependent rate approaches $\log t$, the smaller-order
term blows up.
For example, if $\varphi(t) = (\log t)^2 $, the summands $\e^{- (\ln t)^2 / (2v)}$ in the smaller-order term are larger than $\e^{-1}$
for all $t \leq \e^{\sqrt{2v}}$: the smaller-order term is at least of the order of $\e^{\sqrt{2v}}$, and the regret thus
carries an exponential dependence on $\sqrt{v}$. In the case of a bounded range, this means an exponential dependence on the range $M-m$.
This is probably not an artifact of the proof: in the case of a bounded range,
as long as $\varphi(t) \ll (M-m) \log t$, the lack of exploration bonus entails that
the strategy behaves similarly to a follow-the-leader strategy, which is known to suffer catastrophic, i.e., linear, regret. \\

\changes{
\begin{proof}
As indicated above,
we did not aim for tight numerical constants here and we somehow simplified the standard analysis of UCB
by not considering thresholds of the form $\mu^\star - \varepsilon$ but rather $\mu^\star - \Delta_a/4$.
Hence the non-standard (much increased) numerical factor in front of $\sum_a \ln T/\Delta_a$ when we
specify $\varphi: t \mapsto 2 (M-m)^2 \log t$ into the bound.

In this proof, we repeatedly use that i.i.d.\ random variables $X_1,\ldots,X_n$
with a $v$--sub-Gaussian distribution with expectation $\mu$ satisfy, by the Cram{\'e}r-Chernoff inequality:
for all $\epsilon > 0$,
\[
\P \! \left[ \frac{1}{n} \sum_{i=1}^n X_i \geq \mu + \varepsilon \right]
\leq \inf_{\lambda > 0} \,\, \e^{- n \lambda \varepsilon} \,\, \E \! \left[ \exp \! \left(
- \lambda \sum_{i=1}^n (X_i - \mu) \right) \right]
\leq \inf_{\lambda > 0} \e^{- n \lambda \varepsilon} \, \bigl( \e^{v \lambda^2/2} \bigr)^n
= \e^{-n \varepsilon^2/(2v)}\,;
\]
and we obtain a similar inequality for deviations of the form ``$\leq \mu - \epsilon$''.

Let $a^\star$ and $a$ be an optimal and a suboptimal arm, respectively.
Each arm is pulled once in the first $K$ round.
We bound $\E\big[ N_a(T)]$ by using that for $t \geq K$,
an arm $A_{t+1}$ is pulled only if its index has the highest value,
and then introduce the threshold $\mu^\star - \Delta_a/4$ to separate the $U_{a}(t)$ and the $U_{a^\star}(t)$:
\begin{align*}
	\lefteqn{\E\big[ N_a(T)]} \\
	& \leqslant 1 + \sum_{t=K}^{T-1} \P\bigl[ U_a(t) \geq U_{a^\star}(t) \text{ and } A_{t+1} = a \bigr]  \\
	&\leqslant 1 + \sum_{t=K}^{T-1} \mathbb P\biggl[ U_a(t) \geqslant \mu^\star - \frac{\Delta_a}{4} \text{ and } A_{t+1} = a \biggr] +
	\sum_{t=K}^{T-1} \mathbb P\biggl[ U_{a^\star}(t) \leqslant \mu^\star - \frac{\Delta_a}{4} \biggr]  \\
	&\leq 1+\sum_{t=K}^{T-1} \P\biggl[ \hat \mu_a(t) + \sqrt{\frac{\varphi (t)}{N_a(t)} }\geq \mu^\star - \frac{\Delta_a}{4} \text{ and } A_{t+1} = a \biggr] +
	\sum_{t=K}^{T-1} \mathbb P\biggl[ \hat \mu_{a^\star}(t) \leqslant \mu^\star - \frac{\Delta_a}{4} - \sqrt{\frac{\varphi(t)}{N_{a^\star}(t)}} \biggr]  \\
	&\leq 1 +
    \sum_{t=K}^{T-1} \P\biggl[ \hat \mu_a(t) + \sqrt{\frac{\varphi (T)}{N_a(t)} }\geq \mu^\star - \frac{\Delta_a}{4} \text{ and } A_{t+1} = a \biggr] +
	\sum_{t=K}^{T-1} \mathbb P\biggl[ \hat \mu_{a^\star}(t) \leqslant \mu^\star - \frac{\Delta_a}{4} - \sqrt{\frac{\varphi(t)}{N_{a^\star}(t)}} \biggr]  \\
    &\leq 1 +
    \underbrace{\sum_{n=1}^{T-K+1} \P\biggl[ \hat \mu_{a,n} \geq \mu^\star - \frac{\Delta_a}{4} - \sqrt{\frac{\varphi (T)}{n} } \biggr]}_{\text{Sum($a$)}} +
	\underbrace{\sum_{t=K}^{T-1} \sum_{n=1}^{t-K+1} \mathbb P\biggl[ \hat \mu_{a^\star,n} \leqslant \mu^\star - \frac{\Delta_a}{4} - \sqrt{\frac{\varphi(t)}{n}}\biggr]}_{\text{Sum($a^\star$)}}.
\end{align*}
Note that we used the fact that $\varphi$ is non-decreasing to get to the last but one inequality,
and we used optional skipping for the last one; we denote by $\hat \mu_{a,n}$ and $\hat \mu_{a^\star,n}$ the
average of $n$ i.i.d.\ rewards distributed according to $\nu_a$ and $\nu_{a^\star}$, respectively.

We first deal with Sum($a$).
Let $N_0 = \lceil 4\,  \varphi(T) / \Delta_a^2 \rceil$. For $n \geq N_0$,
\[
\P\biggl[ \hat \mu_{a,n} \geq \mu^\star - \frac{\Delta_a}{4} - \sqrt{\frac{\varphi (T)}{n} } \biggr]
\leq \P\biggl[ \hat \mu_{a,n} \geq \mu_a + \frac{\Delta_a}{4} \biggr] \leq
\e^{-  n\Delta_a^2 / (32v)}\,.
\]
Therefore,
\[
\text{Sum($a$)} \leq N_0 - 1 + \sum_{n \geq N_0} \e^{- n \Delta_a^2 / (32 v)} \leq \frac{4 \varphi(T)}{\Delta_a^2} +
\frac{1}{1 - \e^{- \Delta_a^2 / (32 v)}} \leq \frac{4 \varphi(T)}{\Delta_a^2} + 2 \max \! \left\{ \frac{32 v}{\Delta_a^2}, \, 1 \right\},
\]
where we used\footnote{For a bounded distribution, the case $x > 1$ does not occur as
$x = \Delta_a^2 / (32 v) = \Delta_a^2 / \bigl( 8 (M-m)^2 \bigr) \leq 1/8$; but it may occur
for other sub-Gaussian distributions.}
in the last step $1 / (1 - \e^{-x}) \leq 2 / x$ for $x \in (0,1]$ and $1 / (1 - \e^{-x}) \leq 2$ for $x \geq 1$.

For Sum($a^\star$),
we apply the Cram{\'e}r-Chernoff inequality, then use $(x+y)^2 \geq x^2+y^2$ for $x,y \geq 0$,
and finally apply the same inequalities on $1 / (1 - \e^{-x})$ as for the other sum:
\begin{multline*}
\sum_{t=K}^{T-1} \sum_{n=1}^{t-K+1} \mathbb P\biggl[ \hat \mu_{a^\star,n} \leqslant \mu^\star - \frac{\Delta_a}{4} - \sqrt{\frac{\varphi(t)}{n}}\biggr]
\leq \sum_{t=K}^{T-1} \sum_{n=1}^{t-K+1} \e^{- n \big( \Delta_a / 4 + \sqrt{\varphi(t) / n} \, \big)^2 / (2v) }\\
\leq \sum_{t=K}^{T-1} \sum_{n=1}^{t-K+1} \e^{- n \Delta_a^2 / (32v)} \,\, \e^{- \varphi(t) / (2v)}
\leq \sum_{t=K}^{T-1} \e^{- \varphi(t) / (2v)} \, \frac{1}{1 - \e^{- \Delta_a^2 / (32v)}}
\leq 2 \max \! \left\{ \frac{32 v}{\Delta_a^2}, \, 1 \right\} \sum_{t=K}^{T-1} \e^{- \varphi(t) / (2v)} \,.
\end{multline*}

The proof is concluded by substituting the bounds in $R_T(\nup) = \displaystyle{\sum_{a \in [K]} \Delta_a \, \E\big[ N_a(T)]}$.
\vspace{-.85cm}
\end{proof}
}

\section{Proof of Theorem~\ref{th:adv_talpha_bound}}
\label{sec:proof:thUBadv}

\paragraph{How the second regret bound follows from the first one.}
We substitute the stated values of the $\gamma_t$.
We have, first,
\begin{equation}
\label{eq:calcgt}
\sum_{t=K+1}^T \gamma_t \leq \sqrt{5(1-\alpha)K \log K} \sum_{t=K+1}^T t^{-\alpha}
\leq \sqrt{5(1-\alpha )K \log K} \int_0^T \frac{1}{t^\alpha} \,\d t
= \sqrt{\frac{5 K \ln K }{1-\alpha}} T^{1-\alpha}\,,
\end{equation}
second, using the definition of $\gamma_T$ as a minimum,
\[
\frac{K \ln K}{\gamma_T}
\leq \frac{K \ln K}{1/2} + \frac{T^{\alpha} K \ln K}{\sqrt{5(1-\alpha)K \log K}}
= 2 K \ln K + \sqrt{\frac{K \ln K}{5(1-\alpha)}} \, T^{\alpha}\,,
\]
and third, $\sqrt{T} \leq T^{\max\{\alpha,1-\alpha\}}$,
so that the first regret bound of Theorem~\ref{th:adv_talpha_bound} is further bounded by
\[
(M-m) \sqrt{K \ln K} \Biggl( 3 + 2 \sqrt{\frac{5}{1-\alpha}} \Biggr) T^{\max\{\alpha,1-\alpha\}}
+ 10 (M-m) K \ln K \,.
\]
The claimed expression for $\Phi_{\adv}(T)$ is obtained by bounding $2\sqrt{5}$ by~$5$.

\paragraph{First regret bound.}
In Algorithm~\ref{algo:adahedge},
for time steps $t \geq K+1$, the weights $q_t$ are obtained by using the AdaHedge algorithm of~\citet{de2014follow}
on the payoff estimates $\wh{y}_{t,a}$. AdaHedge is designed for the
case of a full monitoring---not of a bandit monitoring---, but the use of these estimates
emulates a full monitoring. Section~2.2 of~\citet{de2014follow}---see also
an earlier analysis by~\citet{cesa2007improved}---ensures the bound
stated next in Reminder~\ref{rem:adahedge}.

We call pre-regret the quantity at hand in Reminder~\ref{rem:adahedge}:
it corresponds to some regret defined in terms of the payoff estimates.

\begin{reminder}[Application of Lemma~3 and Theorem~6 of~\citealp{de2014follow}]
\label{rem:adahedge}
For all sequences of payoff estimates $\wh{y}_{t,a}$ lying in some bounded real-valued interval,
denoted by $[b,B]$, for all $T \geq K+1$, the pre-regret of AdaHedge satisfies
\begin{multline*}
\max_{k \in [K]} \sum_{t=K+1}^T \wh{y}_{t,k} -
\sum_{t=K+1}^T \sum_{a=1}^K q_{t,a} \, \wh{y}_{t,a} \leq 2 \sum_{t=K+1}^T \delta_t \\
\mbox{where} \qquad \sum_{t=K+1}^T \delta_t
\leq \underbrace{
	\sqrt{\sum_{t=K+1}^T \sum_{a=1}^K q_{t,a} \! \left( \wh{y}_{t,a} - \sum_{k \in [K]} q_{t,k} \, \wh{y}_{t,k} \right)^{\!\! 2} \ln K}
			}_{\leq \sqrt{\sum\limits_{t=K+1}^T \sum\limits_{a=1}^K q_{t,a} (\wh{y}_{t,a} - c)^2 \ln K} \quad \text{for any  } c \, \in \,  \R}
+ (B-b) \left( 1 + \frac{2}{3}\ln K \right)
\end{multline*}
and AdaHedge does not require the knowledge of $[b,B]$ to achieve this bound.
\end{reminder}

The bound of Reminder~\ref{rem:adahedge}
will prove itself particularly handy for three reasons: first, it is valid for
real-valued payoffs; second, it is adaptive to the range of
payoffs; third, the right-hand side looks at first sight not intrinsic enough a bound,
as it also depends on the weights $q_t$, but we will see later that this dependency
is particularly useful in our specific case. To the best of our knowledge,
this is the first direct application of the AdaHedge bound depending on the weights $q_t$
(previous applications were rather solving inequations on the regret, e.g.,
to get improvements for small losses; see \citealp{cesa2007improved}
and \citealp{de2014follow}).

We recall that we
start the summation in Reminder~\ref{rem:adahedge}
at $t = K+1$ because the AdaHedge algorithm is only started at this time, after the initial exploration.
The bound holding ``for any $c \in \R$'' is obtained by a classical bound on the variance. \medskip

\bpc{of the first bound of Theorem~\ref{th:adv_talpha_bound}}
We deal with the contribution of the initial exploration by using the inequality $\max (u + v)\leq \max u + \max v$,
together with the fact that $y_{t, a} - y_{t, A_T} \leq M-m$ for any $a \in [K]$:
\begin{equation}
\label{eq:exp-regret-max}
R_T(y_{1:T})
\leq \underbrace{\max_{a \in [K]} \sum_{t = 1}^K y_{t,a} - \E \! \left[ \sum_{t = 1}^K y_{t,A_t} \right]}_{\leq K(M-m)}
+ \max_{a \in [K]} \sum_{t = K+1}^T y_{t,a} - \E \! \left[ \sum_{t = K+1}^T y_{t,A_t} \right].
\end{equation}
We now transform the pre-regret bound of Reminder~\ref{rem:adahedge}, which is stated with the distributions~$q_t$,
into a pre-regret bound with the distributions~$p_t$; we do so
while substituting the bounds $B = C+ KM/\gamma_T$ and $b = C + Km /\gamma_T$ implied by~\eqref{eq:AH-bound}
and the fact that $(\gamma_t)$ is non-increasing, and by using the definition $q_{t,a} = p_{t,a} - \gamma_t (1/K - q_{t,a})$
for all~$a \in [K]$:
\begin{equation}
\begin{split}
\label{eq:preregret-p}
\max_{k \in [K]} \sum_{t=K +1}^T \wh{y}_{t,k} -
\sum_{t=K +1}^T \sum_{a=1}^K p_{t,a} \, \wh{y}_{t,a}
+ \sum_{t=K+1}^T \gamma_t \sum_{a=1}^K (1/K - q_{t,a}) \, \wh{y}_{t,a} \leq 2 \sum_{t=K+1}^T \delta_t \\
\mbox{where} \qquad \sum_{t=K+1}^T \delta_t
\leq \sqrt{\sum_{t=K+1}^T \sum_{a=1}^K q_{t,a} (\wh{y}_{t,a}  - \C)^2 \ln K}
+ \frac{(M-m)K}{\gamma_T} \left( 1 + \frac{2}{3}\ln K \right).
\end{split}
\end{equation}
As noted by~\citet{AuCBFrSc02},
by the very definition~\eqref{eq:estscheme} of the estimates,
\[
\sum_{a=1}^K p_{t,a} \, \wh{y}_{t,a} = y_{t,A_t}\,.
\]
By~\eqref{eq:condexpwhy}, the tower rule and the fact that $q_t$ is
$H_{t-1}$--measurable, on the one hand,
and the fact that the expectation of a maximum is larger
than the maximum of expectations, on the other hand,
the left-hand side of the first inequality in~\eqref{eq:preregret-p} thus satisfies
\begin{align*}
\lefteqn{\E \! \left[
\max_{k \in [K]} \sum_{t=K +1}^T \wh{y}_{t,k} -
\sum_{t=K+1}^T \sum_{a=1}^K p_{t,a} \, \wh{y}_{t,a}
+ \sum_{t=K+1}^T \gamma_t \sum_{a=1}^K (1/K - q_{t,a}) \, \wh{y}_{t,a} \right]} \\
&\geq
\max_{k \in [K]} \sum_{t=K+1}^T y_{t,k}
- \E \! \left[ \sum_{t=K+1}^T y_{t,A_t} \right]
+ \sum_{t=K+1}^T \gamma_t  \Biggl( \underbrace{\sum_{a = 1}^K y_{t, a} / K}_{\in [m, M]} - \underbrace{\sum_{a = 1}^K \E\bigl[q_{t,a}\bigr] y_{t, a} }_{\in [m, M ]} \Biggr)\\
&\geq \max_{k \in [K]} \sum_{t=K+1}^T y_{t,k}
- \E \! \left[ \sum_{t=K+1}^T y_{t,A_t} \right] - (M-m) \sum_{t=1}^T \gamma_t\,.
\end{align*}
As for the right-hand side of the second inequality in~\eqref{eq:preregret-p}, we first note
that by definition (see line~4 in Algorithm~\ref{alg:exp3_adahedge}),
$p_{t,a} \geq (1-\gamma_t) q_{t,a}$ with $\gamma_t \leq 1/2$ by assumption on the extra-exploration rate,
so that $q_{t,a} \leq 2 p_{t,a}$; therefore, by substituting first this inequality
and then by using Jensen's inequality,
\begin{equation}
\label{eq:Jensen-mainproofUB}
\begin{split}
\E \! \left[ \sqrt{\sum_{t=K+1}^T \sum_{a=1}^K q_{t,a} (\wh{y}_{t,a} - \C)^2 \ln K} \right]
\leq \sqrt{2} \,\, \E \! \left[ \sqrt{\sum_{t=K+1}^T \sum_{a=1}^K p_{t,a} (\wh{y}_{t,a}-\C)^2 \ln K} \right]\\
\leq \sqrt{2} \sqrt{\sum_{t=K+1}^T \sum_{a=1}^K \E \Bigl[ p_{t,a} (\wh{y}_{t,a} - \C)^2 \Bigr] \ln K}\,.
\end{split}
\end{equation}
Standard calculations (see~\citealp{AuCBFrSc02} again) show, similarly to~\eqref{eq:condexpwhy}, that
for all $a \in [K]$,
\[
\E \Bigl[ p_{t,a} (\wh{y}_{t,a} - \C)^2 \,\Big|\, H_{t-1} \Bigr]
= \E \!\left[ \frac{(y_{t,A_t} - \C)^2}{p_{t,a}} \1{A_t = a} \right]
= (y_{t,a} - \C)^2 \leq (M-m)^2\,,
\]
where the last inequality comes from~\eqref{eq:AH-bound}.
By the tower rule, the same upper bound holds for the (unconditional) expectation.
Therefore, taking the expectation of both sides of~\eqref{eq:preregret-p}
and collecting all bounds together, we proved so far
\[
	R_T(y_{1:T}) \leq \underbrace{2\sqrt{2}}_{\leq 3} (M-m) \, \sqrt{KT \ln K}
	+ (M-m)  \frac{K \ln K}{\gamma_T} \underbrace{\left( \frac{2+\gamma_T}{\ln K}+ \frac{4}{3}\right)}_{\leq 5} + (M-m) \sum_{t = K+1}^T \gamma_t\,,
\]
where we used $\gamma_T \leq 1/2$ and $\ln K \geq \ln 2$ as $K \geq 2$.
\end{proof}

\section{Proof of Theorem~\ref{thm:distrib-dependent}}
\label{sec:bernsrem}

Given the decomposition~\eqref{eq:defregrstoch} of the regret, it is necessary and sufficient
to upper bound the expected number of times $\E[N_a(t)]$ any suboptimal arm $a$ is drawn,
where by definition of Algorithm~\ref{alg:exp3_adahedge},
\[
\E[N_a(t)] = 1 + \E \! \left[\sum_{t = K+1}^T \biggl( (1 - \gamma_t)q_{t, a} + \frac{\gamma_t}{K} \biggr)\right]
\leq  1 + \sum_{t =K+1}^T \E [q_{t, a}] + \frac{1}{K} \sum_{t= K+1}^T \gamma_t\,.
\]
We show below (and this is the main part of the proof) that
\begin{equation}
\label{eq:sumlnT}
\sum_{t = K+1}^T \E [q_{t, a}] = \cO(\ln T)\,.
\end{equation}
The straightforward calculations~\eqref{eq:calcgt} already showed that
\[
\frac{1}{K} \sum_{t=K+1}^T \gamma_t \leq \sqrt{\frac{5 \ln K }{(1-\alpha)K}} \, T^{1-\alpha}\,.
\]
Substituting the value~\eqref{eq:defphiadvth2}
of $\Phi^{\AH}_{\free}(T) = \Phi_{\adv}(T)$
and using the decomposition~\eqref{eq:defregrstoch} of $R_T(\nup)$
into $\sum \Delta_a \, \E[N_a(t)]$ then yield
\[
\frac{R_T(\nup)}{T / \Phi^{\AH}_{\free(T)}}
\leq \sum_{a \in [K]} \Delta_a \sqrt{\frac{5 \ln K }{(1-\alpha)K}} \biggl( 3 + \frac{5}{\sqrt{1-\alpha}} \biggr) \sqrt{K \ln K} \bigl( 1 + o(1) \bigr)
+ \cO\!\left( \frac{\ln T}{T^{1-\alpha}} \right),
\]
from which the stated bound follows, via the crude inequality
$3\sqrt{5}\sqrt{1-\alpha}+5 \leq 12$.

\paragraph{Structure of the proof of~\eqref{eq:sumlnT}.}
Let $a^\star$ denote an optimal arm.
By definition of $q_{t, a}$ and by lower bounding a sum of exponential terms
by any of the summands, we get
\[
q_{t, a} = \frac{\exp \! \left( \eta_t \displaystyle{\sum_{s = K+1}^{t-1} \what y_{t, a}}  \right) }{
\displaystyle{\sum_{k = 1}^K \exp \! \left( \eta_t \displaystyle{\sum_{s = K+1}^{t-1} \what y_{t, k}} \right)}}
\leq \exp \!\left( \eta_t \displaystyle{\sum_{t=K+1}^{t-1} (\what y_{t, a} - \what y_{t, a^\star})} \right).
\]
Then, by separating cases, depending on whether $\sum_{t=K+1}^{t-1} (\what y_{t, a} - \what y_{t, a^\star})$
is smaller or larger than the threshold $ - (t-1-K) \Delta_a/2$, and by remembering that the probability $q_{t, a}$ is always smaller than $1$,
we get
\begin{align}
\label{eq:weights_bound}
\sum_{t = K+1}^T \E[ q_{t, a} ] \leq
& \ \ \ \sum_{t = K+1}^T \E \! \left[ \exp \! \bigg( - \eta_t \frac{(t -1 -K) \Delta_a}{2} \bigg) \right] \\
\nonumber
& + \sum_{t = K+1}^T  \P \!  \left[ \sum_{s=K+1}^{t-1} (\what y_{s, a} - \what y_{s, a^\star}) \geq - \frac{(t-1-K) \Delta_a}{2} \right].
\end{align}
We show that the sums in the right-hand side of~\eqref{eq:weights_bound} are respectively $\cO(1)$ and $\cO(\ln T)$.

\paragraph{First sum in the right-hand side of~\eqref{eq:weights_bound}.}
Given the definition
of the learning rates (see the statement of Algorithm~\ref{alg:exp3_adahedge}), namely,
\begin{equation}
\label{eq:reminderdefetat}
\eta_t = \ln K \Bigg/ \sum_{s = K+1}^{t-1} \delta_s\,,
\end{equation}
we are interested in upper bounds on the sum of the $\delta_s$.
Such upper bounds were already derived in the proof of Theorem~\ref{th:adv_talpha_bound};
the second inequality in~\eqref{eq:preregret-p} together with the bound $q_{t,a} \leq 2 p_{t,a}$
stated in the middle of the proof immediately yield
\begin{align*}
\sum_{s = K+1}^{t-1} \delta_s & \leq \sqrt{ \sum_{s = K+1}^t \sum_{a = 1}^K q_{s, a} \big(  \what y_{s, a} - C \big)^2 \ln K}
+ \frac{(M-m)K}{\gamma_t} \left( 1 + \frac{2}{3} \ln K\right) \\
& \leq \sqrt{2} \sqrt{\sum_{s = K+1}^t \sum_{a = 1}^K p_{s, a} \big( \what y_{s, a} - C \big)^2 \ln K}
+ \frac{(M-m)K}{\gamma_t} \left( 1 + \frac{2}{3} \ln K\right).
\end{align*}
Unlike what we did to complete the proof of Theorem~\ref{th:adv_talpha_bound},
we do not take expectations and rather proceed with deterministic bounds.
By the definition~\eqref{eq:estscheme}
of the estimated payoffs for the equality below,
by~\eqref{eq:AH-bound} for the first inequality below,
and by the fact that the exploration rates are non-increasing
for the second inequality below,
we have, for all $s \geq K+1$,
\begin{equation}
\label{eq:manipyhat}
\sum_{a = 1}^K p_{s, a} \big( \what y_{s, a} - C \big)^2
= \frac{\big(  y_{s, A_s} - C \big)^2}{p_{s, A_s}}
\leq \frac{(M-m)^2}{\gamma_s/K}
\leq \frac{(M-m)^2}{\gamma_t/K}\,.
\end{equation}
Therefore,
\[
\sum_{s = K+1}^{t-1} \delta_s \leq \sqrt{2} (M-m) \sqrt{\frac{t \, K \ln K}{\gamma_t}}
+ \frac{(M-m)K}{\gamma_t} \left( 1 + \frac{2}{3} \ln K\right)
\eqdef D_t = \Theta\Bigl(\sqrt{t/\gamma_t} + 1/\gamma_t \Bigr)\,.
\]
For the sake of concision, we denoted by $D_t$ the obtained bound.
Via the definition~\eqref{eq:reminderdefetat} of $\eta_t$,
the sum of interest is in turn bounded by
\[
\sum_{t = K+1}^T
\exp \! \left( - \eta_t \big(t-1-K\big) \frac{\Delta_a}{2}\right)
\leq \sum_{t = K+1}^T
\exp \! \left( - \frac{\Delta_a \ln K}{2} \,\, \frac{t-1-K}{D_t} \right)
= \cO(1)\,,
\]
where the equality to $\cO(1)$, i.e., the fact that the considered series is bounded,
follows from the fact that
\[
- (t-1-K)/D_t = \Theta \Bigl( \sqrt{t \gamma_t} + t \gamma_t \Bigr) = \Theta\bigl( t^{(1-\alpha)/2} + t^{1-\alpha} \bigr)\,.
\]

\paragraph{Second sum in the right-hand side of~\eqref{eq:weights_bound}.}
We will use Bernstein's inequality for martingales, and more specifically,
the formulation of the inequality by \citet[Thm. 1.6]{freedman1975tail}---see
also \citet[Section~2.2]{Mas03}---, as stated next.
\begin{reminder}
\label{rm:bernstein_mart}
Let $(X_n)_{n \geq 1}$ be a martingale difference sequence with respect to a filtration $(\mathcal{F}_n)_{n \geq 0}$,
and let $N \geq 1$ be a summation horizon.
Assume that there exist real numbers $b$ and $v_N$ such that, almost surely,
\begin{equation*}
\forall n \leq N, \quad X_n \leq b
\qquad \mbox{and} \qquad
\sum_{n = 1}^N \E\big[ X_n^2 \, \big| \,  \mathcal F_{n-1}\big] \leq v_N\,.
\end{equation*}
Then for all $\delta \in (0,1)$,
\[
\P \! \left[  \sum_{n = 1}^N X_n \geq \sqrt{2 v_N \log \frac{1}{\delta}}
+ \frac{b}{3} \log \frac{1}{\delta} \right] \leq \delta\,.
\]
\end{reminder}
For $s \geq K+1$, we consider the increments $X_s = \Delta_a -  \what y_{s, a^\star} + \what y_{s,a}$,
which are adapted to the filtration $\cF_s = \sigma(A_1,Z_1,\ldots,A_s,Z_s)$,
where we recall that $Z_1,\ldots,Z_s$ denote the payoffs obtained in rounds $1,\ldots,s$.
Also, as $p_s$ is measurable with respect to past information $\cF_{s-1}$ and
since payoffs are drawn independently from everything else (see Section~\ref{sec:setting-stochastic}), we have,
by the definition~\eqref{eq:estscheme}
of the estimated payoffs (where we rather denote by $Y_{s,a}$ the payoffs drawn at random according
to $\nu_a$, to be in line with the notation of Section~\ref{sec:setting-stochastic} for stochastic bandits):
for all $a \in [K]$,
\[
\E\big[ \,\what y_{s, a} \,\big|\, \mathcal F_{s-1} \big]
= \frac{\E[ Y_{s, a} \,|\, \mathcal F_{s-1}] - C}{p_{s,a}} \1{A_s = a} + C
= \frac{\mu_a - C}{p_{s,a}} \1{A_s = a} + C = \mu_a\,.
\]
As a consequence, $\E[X_s\,|\,\cF_{s-1}] = \E\bigl[\Delta_a -  \what y_{s, a^\star} + \what y_{s,a}\,|\,\cF_{s-1}\bigr] = 0$.
Put differently, $(X_s)_{s \geq K+1}$ is indeed a martingale difference sequence with respect to the
filtration $(\mathcal{F}_s)_{s \geq K}$.

We now check that the additional assumptions of Reminder~\ref{rm:bernstein_mart}
are satisfied. Manipulations and arguments similar to the ones used in~\eqref{eq:AH-bound} and~\eqref{eq:manipyhat}
show that for all $s \geq K+1$,
\begin{align*}
\Delta_a -  \what y_{s, a^\star} + \what y_{s,a}
& \leq \Delta_a - \frac{Y_{s, a^\star} - C}{p_{s,a}} \1{A_s = a^\star} + \frac{Y_{s, a} - C}{p_{s,a}} \1{A_s = a} \\
& \leq (M-m) (1 + K/\gamma_s) \leq b \eqdef (M-m) (1 + K/\gamma_t) \,.
\end{align*}
For the variance bound, we first note that for all $s \leq t-1$,
we have $(\what y_{s, a} - C)(\what y_{s, a^\star} - C) = 0$ because of the indicator functions,
and therefore,
\begin{align*}
\E  \Bigl[ \big(\Delta_a -  \what y_{s, a^\star} + \what y_{s,a }\big)^2 \,\Big|\, \cF_{s-1} \Bigr]
& \leq \E  \Bigl[ \big(\what y_{s, a^\star} + \what y_{s,a }\big)^2 \,\Big|\, \cF_{s-1} \Bigr] \\
& \leq \E \Bigl[ \big(\what y_{s, a^\star} - C \big)^2 \,\Big|\, \cF_{s-1} \Bigr]
+ \E \Bigl[ \big(\what y_{s, a} - C \big)^2 \,\Big|\, \cF_{s-1} \Bigr]\,;
\end{align*}
in addition, for all $a \in [K]$, including $a^\star$,
\[
\E \Bigl[ \big(\what y_{s, a} - C \big)^2 \,\Big|\, \cF_{s-1} \Bigr]
= \E\biggl[ \frac{(Y_{s, A_s}-C)^2}{p_{s,a}^2} \1{A_s = a} \,\bigg|\, \mathcal F_{s-1} \biggr]
\leq  \frac{(M-m)^2}{p_{s, a}} \leq \frac{(M-m)^2 K}{\gamma_t}\,.
\]
Therefore
\begin{equation*}
\sum_{s = K+1}^{t-1} \E \Bigl[ \big(\Delta_a -  \what y_{s, a^\star} + \what y_{s,a }\big)^2 \,\Big|\, \cF_{s-1} \Bigr]
\leq \frac{2 K(M-m)^2 (t - 1- K)}{\gamma_t} \leq v_t \eqdef \frac{2 (M-m)^2 t K}{\gamma_t}\,.
\end{equation*}

Bernstein's inequality (Reminder~\ref{rm:bernstein_mart}) may thus be applied; the choice $\delta = 1/t$
therein leads to
\[
\P\!\Biggl[ \sum_{s = K+1}^{t-1} \big(\Delta_a - (\what y_{s, a^\star}  - \what y_{s, a} )\big)
\geq \underbrace{2(M-m)\sqrt{\frac{t K}{\gamma_t} \log t}+\frac{M-m}{3} \left( 1 + \frac{K}{\gamma_t} \right) \log t}_{\eqdef \normalfont D'_t} \Biggr]
\leq \frac{1}{t}\,.
\]

As $\sqrt{t / \gamma_t} = \mathcal O(t^{(1+\alpha)/2})$ and $1/\gamma_t = \mathcal O(t^\alpha)$ as $t \to \infty$,
where $\alpha < 1$, and as $\Delta_a > 0$ given that we are considering a suboptimal arm~$a$,
there exists $t_0 \in \N$ such that for all $t \geq t_0$,
\[
D'_t \leq \frac{(t-1-K) \Delta_a}{2}
\]
thus
\begin{align*}
\P\!\Biggl[ \sum_{s = K+1}^{t-1} (\what y_{s, a}  - \what y_{s, a^\star})
\geq - \frac{(t-1-K) \Delta_a}{2} \Biggr]
& =
\P\!\Biggl[ \sum_{s = K+1}^{t-1} \big(\Delta_a - (\what y_{s, a^\star}  - \what y_{s, a} )\big)
\geq \frac{(t-1-K) \Delta_a}{2} \Biggr] \\
& \leq
\P\!\Biggl[ \sum_{s = K+1}^{t-1} \big(\Delta_a - (\what y_{s, a^\star}  - \what y_{s, a} )\big)
\geq D'_t \Biggr] \leq \frac{1}{t}\,.
\end{align*}
Therefore, as $T \to \infty$
\[
\sum_{t = 1}^T \P \! \left[ \sum_{t=K+1}^{t-1} (\what y_{t, a} - \what y_{t, a^\star}) \geq - \frac{(t-1-K) \Delta_a}{2} \right]
= \mathcal O (\ln T)\,,
\]
as claimed. This concludes the proof.

\bibliography{biblio_range_adaptation}

\clearpage
\begin{center}
{\bf \LARGE Supplementary material for \\[7pt]
``Adaptation to the Range in $K$--Armed Bandits''} \\[15pt]
{\Large by H{\'e}di Hadiji and Gilles Stoltz}
\end{center}

\vspace{2cm} \

\section{Adaptation to the Range for Linear Bandits}
\label{sec:linear-app}

To illustrate the generality of the techniques discussed in this paper, we quickly describe how these
can be used to obtain range adaptive algorithms for linear bandits. This section is meant for illustration
and not for completeness. In particular, we focus on the case of (oblivious) adversarial linear bandits:
we refer the reader to \citet[Chapter~27]{lattimore_szepesvari_2020}, which we follow closely, for a more
thorough description of the setting; we do not describe the application of our techniques to stochastic linear bandits.

\paragraph{Learning protocol.}
A finite action set $\cA \subset \R^d$, of cardinality $K$, is given. (The setting of vanilla $K$--armed bandits considered in
the rest of the article corresponds to $\cA$ formed by the vertices of the probability simplex of $\R^K$.)
The environment selects beforehand a sequence $(y_t)_{t \geq 1}$ of vectors in $\R^d$ satisfying a boundedness assumption:
there exists an interval $[m,M]$ such that
\begin{equation}
\label{eq:bounded-linearbandits}
\forall t \geq 1, \ \
\forall x \in \cA, \qquad
\dotptr{x}{y_t} \in [m,M]\,.
\end{equation}
We assume that the player does not know in advance $m$ nor $M$.
To simplify the exposition, we also assume that $m \leq 0 \leq M $.

At every time step, the player chooses an action $X_t \in \mathcal A$ and receives and only observes the payoff $\dotptr{X_t}{y_t}$.
It does not observe $y_t$ nor the payoffs $\dotptr{x}{y_t}$ associated with choices $x \ne X_t$. The action
$X_t$ is chosen independently at random according to a distribution over $\cA$ denoted by $p_t = \bigl( p_t(a) \bigr)_{a \in \cA}$.

The expected regret is defined as
\[
R_T(y_{1:T}) = \max_{x \in \mathcal A} \sum_{t = 1}^T \dotptr{x}{y_t} - \E \! \left[ \sum_{t = 1}^T \dotptr{X_t}{y_t} \right].
\]

\paragraph{Estimating the unobserved payoffs.}
As in the case of vanilla $K$--armed bandits, the key is to estimate unobserved payoffs.
We may actually build an estimate $\what y_t$ of the vectors $y_t$, from which we form
the estimates $\dotptr{x}{\what{y}_t}$. This estimate takes advantage of the linear structure
of the problem.

Fix a distribution $\pi$ such that the non-negative symmetric matrix
\[
M(\pi) \eqdef \sum_{x \in \mathcal A} \pi(x) \, x x^{\!\top}
\]
is invertible: such a distribution exists whenever $\cA$ spans $\R^d$,
which we may assume with no loss of generality; see Lemma~\ref{lem:lin_extra_expl}
below. This distribution $\pi$ will be used to explore the arms; it is in general not uniform
over the arms.
For all distributions $q$ over $\cA$ and
all $\gamma \in (0,1]$, the distribution $p = (1 - \gamma) q + \gamma \pi$
is such that the non-negative symmetric matrix $M(p)$ is invertible as well
(as it is larger than $\gamma\,M(\pi)$, in the sense of the partial inequality $\succcurlyeq$
over non-negative symmetric matrices). We only use distributions $p_t$ of this form.
We may then define
\begin{equation}
\label{eq:lin_est_formula}
\what y_t = M(p_t)^{-1} \, X_t \dotptr{X_t}{y_t}
\end{equation}
and note that
\begin{equation}
\label{eq:ylin-unbiased}
\E\bigl[ \what y_t \,\big|\, p_t \bigr]
=  M(p_t)^{-1} \Biggl( \underbrace{\sum_{x \in \mathcal A} p_t(x) \, x x^{\!\top}}_{= M(p_t)} y_t \Biggr) = y_t\,;
\end{equation}
indeed, conditioning on $p_t$ amounts to integrating over the random choice of $X_t$ according to $p_t$.

\paragraph{An algorithm adaptive to the unknown range.}
When the range is given, a well-known strategy is to use plain exponential weights over actions in $\mathcal A$
with the estimates $\dotptr{x}{\what{y}_t}$ to obtain distributions $q_t$ that are then mixed with $\pi$
to form the final distributions $p_t$.
When the range is unknown, we suggest to simply replace plain exponential weights
with AdaHedge (the difference lies in the tuning of the rates $\eta_t$), which leads to Algorithm~\ref{algo:adahedge_lin}.
In this algorithm, we refer to rates $\gamma_t$ as exploration rates (and not as extra-exploration rates as in
Algorithm~\ref{algo:adahedge}) and similarly, to $\pi$ as the exploration distribution. This is because
for adversarial linear bandits, exploration was always required even to get expected results (unlike
for $K$--armed bandits, see the introduction of Section~\ref{sec:mainadahedge}).

\begin{figure}[t]
	\renewcommand\footnoterule{}
	\begin{algorithm}[H]
		\begin{algorithmic}[1]
			\caption{\label{algo:adahedge_lin} AdaHedge for adversarial linear bandits}
			\STATE \textbf{Input:} an exploration distribution $\pi$ over $\mathcal A$ and exploration rates $(\gamma_t)_{t \geq 1}$ in $[0,1]$
			\STATE \textbf{Initialization:} $\eta_1 = +\infty$ and $q_{1}$ is the uniform distribution over $\cA$
			\FOR{rounds $t = 1,  \dots$}
			\STATE Define $p_t$ by mixing $q_t$ with $\pi$ according to \vspace{-.25cm} \\
			\[
			p_t = (1 - \gamma_t) q_t + \gamma_t \pi
			\vspace{-.3cm}
			\]
			\STATE Draw an arm $X_t \sim p_t$ (independently at random according to the distribution $p_t$)
			\STATE Get and observe the payoff $\dotptr{X_t}{y_t}$
			\STATE Compute estimates $\dotptr{x}{\what{y}_t}$ of all payoffs according to~\eqref{eq:lin_est_formula}
			\STATE Compute the mixability gap $\delta_t$ based on the distribution $q_t$ and on these estimates:
			\[
			\delta_t =\left\{
			\begin{split}
			& - \sum_{x \in \mathcal A} q_t(x)\,\dotptr{x}{\what{y}_t} +
                \frac{1}{\eta_t} \log \Biggl( \sum_{x \in \mathcal A} q_t(x) \e^{\eta_t \dotptr{x}{\what{y}_t}} \Biggr)
			\quad & \text{if } \eta_t < + \infty \\
			& - \sum_{x \in \mathcal A} q_t(x)\, \dotptr{x}{\what{y}_t} + \max_{x \in \mathcal A} \dotptr{x}{\what{y}_t} \quad & \text{if } \eta_t = + \infty
			\end{split}\right.
			\]
			\STATE Compute the learning rate $\displaystyle{\eta_{t+1} = \Biggl( \sum_{s = 1}^t \delta_s \Biggr)^{-1}} \ln K$
			\STATE Define $q_{t+1}$ component-wise as \vspace{-.3cm}
			\[
			~ \hspace{2.5cm} q_{t+1}(a) = \exp \! \left( \eta_{t+1} \sum_{s =1}^{t} \dotptr{a}{\what{y}_s} \right)
			\Bigg/ \sum_{x \in \mathcal A} \exp \! \left( \eta_{t+1} \sum_{s = 1}^{t} \dotptr{x}{\what{y}_s} \right) \vspace{-.3cm}
			\]
			\ENDFOR
		\end{algorithmic}
	\end{algorithm}
\end{figure}

The analysis of this algorithm relies on the same ingredients as the ones already encountered
in Section~\ref{sec:adahedge-distrfree}, with the addition of the following lemma, that quantifies
the quality of the exploration.
This lemma requires that $\mathcal A$ spans $\R^d$, which we may assume with no loss
of generality (otherwise, we just replace $\R^d$ by the vector space generated by $\cA$).

\begin{lemma}[{\citealp[Theorem 21.1]{lattimore_szepesvari_2020}}]
\label{lem:lin_extra_expl}
There exists a distribution $\pi$ over $\mathcal A$ such that
\[
M(\pi) = \sum_{x \in \mathcal A} \pi(x) \, x x^{\!\top} \;  \mbox{is invertible} \qquad  \mbox{and} \qquad
\max_{x \in \mathcal A}  x^\top M(\pi)^{-1} x  = d\,.
\]
\end{lemma}

We are now ready to state the main result of this section. It is the counterpart of Theorem~\ref{th:adv_talpha_bound};
for the sake of simplicity, we only state it for the value $\alpha = 1/2$.

\begin{theorem}
AdaHedge for adversarial linear bandits (Algorithm~\ref{algo:adahedge_lin})
with the extra-exploration
\[
\gamma_t = \min \Bigl\{1/2, \, \sqrt{2.5 \,  d (\log K) t^{-1/2}} \Bigr\}
\]
ensures that for all bounded ranges $[m,M]$ containing~$0$, for all oblivious individual
sequences $y_1,\,y_2,\ldots$ satisfying the boundedness condition~\eqref{eq:bounded-linearbandits},
\[
R_T(y_{1:T}) \leq 12(M-m) \sqrt{d T \ln K} + 18(M-m) d \log K\,.
\]
\end{theorem}

The proof starts by following closely the one of Theorem~\ref{th:adv_talpha_bound}
(provided in Appendix~\ref{sec:proof:thUBadv});
the differences are underlined and dealt with in the second part of the proof.
\medskip 

\bp
By Reminder~\ref{rem:adahedge}, since the player plays the AdaHedge strategy over the payoff estimates $\dotptr{x}{\what{y}_t}$,
the pre-regret satisfies
\[
\max_{x \in \mathcal A} \sum_{t = 1}^T  \dotptr{x}{\what{y}_t} - \sum_{t = 1}^T \sum_{a \in \mathcal A} q_t(a) \, \dotptr{a}{\what{y}_t}
\leq 2 \sqrt{V_T \log K} +  M_T \left(2 + \frac{4}{3} \log K \right)
\]
with $\displaystyle{V_T = \sum_{t = 1}^T \sum_{x \in \mathcal A} q_t(x) \, \bigl( \dotptr{x}{\what{y}_t} \bigr)^2}$
and
\[
M_T = \max \bigl\{ \dotptr{x}{\what{y}_t} : t \leq T \ \mbox{and} \ x \in \cA \bigr\}
- \min \bigl\{ \dotptr{x}{\what{y}_t} : t \leq T \ \mbox{and} \ x \in \cA \bigr\}\,.
\]
Since $\gamma_t \leq 1/2$, we have
$q_t(x) \leq 2 \, p_t(x)$ for all $x \in \cA$. We therefore define
\[
V'_T = \sum_{t = 1}^T \sum_{x \in \mathcal A} p_t(x) \, \bigl( \dotptr{x}{\what{y}_t} \bigr)^2
\]
and have $V_t \leq 2 V'_T$.
By the tower rule, based on the equality~\eqref{eq:ylin-unbiased},
and given that the expectation of a maximum is larger than the maximum of the
expectations (for the first inequality),
and by the definition of the $p_t$ (for the second inequality), we have
proved so far that
\begin{align*}
R_T(y_{1:T})
& \leq \E\! \left[ \max_{x \in \mathcal A} \sum_{t = 1}^T  \dotptr{x}{\what{y}_t}
- \sum_{t = 1}^T \sum_{a \in \mathcal A} p_t(a) \, \dotptr{a}{\what{y}_t} \right] \\
& \leq \E\! \left[ \max_{x \in \mathcal A} \sum_{t = 1}^T  \dotptr{x}{\what{y}_t}
- \sum_{t = 1}^T \sum_{a \in \mathcal A} q_t(a) \, \dotptr{a}{\what{y}_t} \right]
+ \E \! \left[ \sum_{t = 1}^T \gamma_t \sum_{a \in \mathcal A} \bigl( \pi(a) - q_t(a) \bigr) \, \dotptr{a}{\what{y}_t} \right] \\
& \leq \E\Bigg[2\sqrt{2 V'_T \log K} + M_T \left(2 + \frac{4}{3} \log K \right) \Bigg]
+ \sum_{t = 1}^T \gamma_t \underbrace{\sum_{a \in \mathcal A} \bigl( \pi(a) - q_t(a) \bigr) \, \dotptr{a}{y_t}}_{ \leq (M-m)} \,.
\end{align*}
Hence by Jensen's inequality
and by the bounds $\E[V'_T] \leq (M-m)^2 d T$ and $M_T \leq 2(M-m)d/\gamma_T$
proved below, we finally get
\begin{align*}
R_T(y_{1:t})
& \leq 2 \sqrt{2 \, \E[V'_T] \log K} + \E[M_T] \left(2 + \frac{4}{3} \log K \right)
+ (M-m) \sum_{t = 1}^T \gamma_t \\
&\leq 2\sqrt{2} (M-m) \sqrt{d T \ln K} + \left(2 + \frac{4}{3} \ln K \right) \frac{2(M-m)d}{\gamma_T} + (M-m) \sum_{t = 1}^T \gamma_t \\
& \leq 3(M-m) \sqrt{d T \ln K} + 9(M-m) \frac{d\ln K}{\gamma_T} + (M-m) \sum_{t = 1}^T \gamma_t\,.
\end{align*}
Replacing the $\gamma_t$ by their values and using the same bounds as at
the beginning of Appendix~\ref{sec:proof:thUBadv} yields the claimed result;
the factor $12$ in the bound comes from
\[
3 + \sqrt{10} + 9 \sqrt{\frac{2}{5}} \leq 12\,.
\]
We only need to prove the two claimed bounds to complete the proof;
they can be extracted from the proof of Theorem 27.1 by \citet{lattimore_szepesvari_2020}
but we provide derivations for the sake of completeness. \smallskip

\noindent \emph{Proof of $M_T \leq 2(M-m)d/\gamma_T$.}
We fix $x \in \cA$ and $t \leq T$. We recall that
$M(p_t)$ and thus $M(p_t)^{-1}$ are positive definite symmetric matrices.
By the Cauchy-Schwarz inequality applied with the norm induced by
the positive $M(p_t)^{-1}$,
\[
\bigl| x^{\!\top} M(p_t)^{-1} X_t \bigr|
\leq \sqrt{x^{\!\top} M(p_t)^{-1} x} \,\, \sqrt{X_t^{\!\top} M(p_t)^{-1} X_t}
\leq \max_{x \in \mathcal A} \Big\{ x^{\!\top} M(p_t)^{-1} x \Big\}\,.
\]
As indicated right before~\eqref{eq:lin_est_formula},
we have $M(p_t) \succcurlyeq \gamma_t \, M(\pi)$ and therefore
$M(p_t)^{-1} \curlyeqprec M(\pi)^{-1}/\gamma_t$. This entails
\[
\bigl| x^{\!\top} M(p_t)^{-1} X_t \bigr|
\leq \frac{1}{\gamma_t} \max_{x \in \mathcal A} \Big\{ x^{\!\top} M(\pi)^{-1} x \Big\}
= \frac{d}{\gamma_t}
\leq \frac{d}{\gamma_T}\,,
\]
where the equality follows from Lemma~\ref{lem:lin_extra_expl}
and where we used $\gamma_T \leq \gamma_t$ for the second inequality.
Finally, keeping in mind that we assumed $m \leq 0 \leq M$,
\[
x^{\!\top} \what{y}_t =
\underbrace{x^{\!\top} M(p_t)^{-1} \, X_t}_{\in [-d/\gamma_t,d/\gamma_t]} \underbrace{\dotptr{X_t}{y_t}}_{\in [m,M]}
\in \left[ - \frac{d \max\{-m,M\}}{\gamma_T}, \, \frac{d \max\{-m,M\}}{\gamma_T} \right],
\]
from which the bound
\[
M_t = 2 \, \frac{d \max\{-m,M\}}{\gamma_T} \leq \frac{2 d (M-m)}{\gamma_T}
\]
follows, as desired. \smallskip

\noindent
\emph{Proof of $\E[V'_T] \leq (M-m)^2 d T$.}
Since $\bigl| \dotptr{X_t}{y_t} \bigr| \leq \max\{-m,M\} \leq M-m$,
the definition~\eqref{eq:lin_est_formula} leads to
\begin{align*}
\bigl( \dotptr{x}{\what{y}_t} \bigr)^2
= \Bigl( x^{\!\top} M(p_t)^{-1} \, X_t \dotptr{X_t}{y_t} \Bigr)^2
& \leq (M-m)^2 \, \Bigl( x^{\!\top} M(p_t)^{-1} \, X_t \Bigr)^2 \\
& = (M-m)^2 \, \dotptr{X_t}{M(p_t)^{-1} x} \dotptr{x}{M(p_t)^{-1} X_t}\,.
\end{align*}
Therefore, summing over $x \in \cA$ and using the very definition of $M(p_t)$, we get
\begin{align*}
\sum_{x \in \cA} p_t(x) \, \bigl( \dotptr{x}{\what{y}_t} \bigr)^2
& \leq (M-m)^2 \, X_t^{\!\top} M(p_t)^{-1} \left( \sum_{x \in \cA} p_t(x) \, x x^{\!\top} \right) M(p_t)^{-1} X_t \\
& = (M-m)^2 \, X_t^{\!\top} M(p_t)^{-1} X_t
= (M-m)^2  \Tr \Bigl( M(p_t)^{-1} X_t X_t^\top \Bigr)\,.
\end{align*}
Now, by the linearity of the trace,
\[
\E \biggl[ \Tr \Bigl( M(p_t)^{-1} X_t X_t^\top \Bigr) \biggr]
= \E\! \left[ \sum_{x \in \mathcal A} p_t(x) \Tr \Bigl( M(p_t)^{-1} x x^\top \Bigr) \right]
= \E \bigl[ \Tr(I_d) \bigr] = d\,,
\]
where $I_d$ is the $d$--dimensional identity matrix.
Collecting all bounds together and summing over $t$
yields the claimed inequality $\E[V'_T] \leq (M-m)^2 d T$.
\end{proof}


\clearpage
\section{Proof of $\Kinf\big( \nu, \mu, \Dmod_{m, M} \big) = \Kinf\big( \nu, \mu, \Dmod_{-\infty, M} \big)$ in Appendix~\ref{sec:LRBK}}
\label{sec:app-a-e}

For the sake of readability, we formally restate the equality to be proved. 

\begin{proposition}
\label{prop:Dmod-m}
Fix $M \in \R$.
For all $m \leq M$, for all $\nu \in \Dmod_{m, M}$ and all $\mu > \Ed(\nu)$,
\[
\Kinf\big( \nu, \mu, \Dmod_{m, M} \big)
= \Kinf\big( \nu, \mu, \Dmod_{-\infty, M} \big)\,.
\]
\end{proposition}

\bp
The inequality $\geq$ is immediate, as the right-hand side of the equality is an infimum over the larger set $\Dmod_{-\infty, M}$.
For the inequality $\leq$, we may assume with no loss of generality that $\mu < M$, as otherwise,
there is no distribution $\nu'$ neither in $\Dmod_{m, M}$ nor in $\Dmod_{-\infty, M}$ with
$\Ed(\nu') > \mu \geq M$, so that both $\Kinf$ quantities equal $+\infty$.

We fix $M$, $m$, $\nu$ and $\mu$ as in the statement of the proposition.
It suffices to show that in the case $\mu < M$, for all $\nu' \in \mathcal \Dmod_{- \infty, M}$ with $\Ed(\nu')> \mu$
and $\nu \ll \nu'$, there exists $\nu'' \in \mathcal \Dmod_{m, M}$ with $\Ed(\nu'')> \mu$ and $\KL(\nu, \nu'') \leq \KL(\nu, \nu')$.
(If $\nu$ is not absolutely continuous with respect to $\nu'$, then $\KL(\nu, \nu') = +\infty$ and taking $\nu''$ as
the Dirac mass $\delta_M$ at $M$ is a suitable choice.)
To do so, given such a distribution $\nu'$, we first note that $\nu \ll \nu'$ and $\nu \in \Dmod_{m, M}$, i.e., $\nu([m,M]) = 1$,
entail that $\nu'([m,M]) > 0$, so that we may define the restriction $\nu'' = \nu'_{[m,M]}$ of $\nu'$ to $[m, M]$; its density
with respect to $\nu'$ is given by
\[
\frac{\dint \nu''}{\dint \nu'}(x)  = \nu'\bigl( [m,M] \bigr)^{-1} \,  \1{x \in [m,M]} \qquad \nu'\text{--a.s. for all $x \in \R$.}
\]
We have the absolute-continuity chain $\nu \ll \nu'' \ll \nu'$, and the Radon-Nykodym derivatives thus defined satisfy
\begin{equation}
\label{eq:chainruledensity}
\frac{\dint \nu}{\dint \nu'}(x) = \frac{\dint \nu}{\dint \nu''}(x) \frac{\dint \nu''}{\dint \nu'}(x) =
\nu'\bigl( [m,M] \bigr)^{-1} \, \frac{\dint \nu}{\dint \nu''}(x) \,  \1{x \in [m,M]} \qquad \nu'\text{--a.s. for all $x \in \R$.}
\end{equation}
Moreover $\Ed(\nu'') \geq \Ed(\nu')$, and thus $\Ed(\nu'') > \mu$, as
\begin{align*}
\Ed(\nu') & = \int_{(-\infty,m)} x \dint \nu'(x) + \int_{[m,M]} x \dint \nu'(x) \\
& \leq \Big( 1 - \nu'\big( [m,M] \bigr) \Big) m + \nu'\big( [m,M] \bigr) \, \Ed(\nu'')
\leq \Ed(\nu'')\,.
\end{align*}
Finally, by~\eqref{eq:chainruledensity}, which also holds $\nu$--almost surely,
and the definition of Kullback-Leibler divergences,
\begin{align*}
\KL(\nu, \nu')
= \int_{(-\infty,M]} \ln \! \left( \frac{\dint \nu}{\dint \nu'} \right) \! \dint \nu
& = - \ln \nu'\bigl( [m,M] \bigr) + \int_{[m,M]} \ln \! \left( \frac{\dint \nu}{\dint \nu''} \right) \! \dint \nu \\
& = - \ln \nu'\bigl( [m,M] \bigr) + \KL(\nu, \nu'') \geq \KL(\nu, \nu'')\,.
\end{align*}
This concludes the proof.
\end{proof}

\clearpage
\section{Known~$M$ but Unknown~$m$: Adaptation to the Range with a $\sqrt{KT}$ \\ 
~~~~~~~~~~~~~~~~~~Scale-Free Distribution-Free Regret Bound}
\label{sec:TsallisKTpossible}

This appendix details a claim made in the second part of Remark~\ref{rk:lnN}:
that when the upper end $M$ of the range is known, and adaptation is only 
with respect to the lower end $m$ of the range, then a $\sqrt{KT}$ scale-free distribution-free regret upper bounds may be achieved, 
which exactly matches the distribution-free lower bound.
This is the main result of this appendix, to be stated as Theorem~\ref{thm:tsallis} in Appendix~\ref{app:proofs_tsallis}.
The full outline of this section of the appendices is detailed below.

\paragraph{Disclaimer.}
In the case where $M$ is known and adaptation is only to $m$, 
we could not exhibit a strategy that would simultaneously achieve both optimal distribution-dependent
and distribution-free regret bounds, unlike what is known in the case of a known payoff range (the KL-UCB-switch
strategy by~\citealp{garivier2018kl}) and unlike
what we achieved in the main body of the article when adapting to the unknown range $[m,M]$, or unknown
upper end $M$ on the range but known lower end $m$. 

We however conjecture that this should be possible and that, at least, no trade-off exists
between the two bounds (i.e., we conjecture that Theorem~\ref{thm:adaptive_lower_bound} should not hold).

\paragraph{Outline of Appendix~\ref{sec:TsallisKTpossible}.}
All results of this appendix rely on the AdaFTRL methodology of~\citet{orabona2018scale},
which we recall first in Appendix~\ref{sec:AdaFTRL-full}.
AdaFTRL stands for adaptive follow-the-regularized-leader
and it was partially built on and inspired by the analysis for AdaHedge, which is a special case of AdaFTRL with entropic
regularizer (see \citealp{de2014follow} for AdaHedge, as well as the earlier analysis by \citealp{cesa2007improved}).
\citet{koolen_2016} proposes an alternative analysis of AdaFTRL, closer to the AdaHedge formulation,
namely, using directly some mixability gaps instead of upper bounds thereon; this is the analysis we actually recall
in Section~\ref{sec:AdaFTRL-full}.

Appendix~\ref{app:ada_hedge_proof} specializes the general results of 
Appendix~\ref{sec:AdaFTRL-full} to an entropic regularizer, leading to AdaHedge.
It provides a proof of the AdaHedge bound, i.e., Reminder~\ref{rem:adahedge},
in order to make this article self-complete. 

An interesting observation, described in Appendix~\ref{app:proofs_exp3},
is that (as in the case of a fully known payoff range)
AdaHedge does not require any extra-exploration (i.e., any mixing with the uniform
distribution) to achieve a scale-free distribution-free regret bound of order
$\smash{(M-m)\sqrt{KT \ln K}}$. 

We then turn to the main result of this appendix, stated in Appendix~\ref{app:proofs_tsallis},
which is a $\sqrt{KT}$ scale-free distribution-free regret bound for AdaFTRL
with $1/2$--Tsallis entropy, in the case where $M$ is known and adaptation is only to~$m$.
The choice of the $1/2$--Tsallis entropy as a regularizer is motivated by
the INF strategy of \citet{audibert2009Minimax}, which can be seen as an instance of FTRL
with $1/2$--Tsallis entropy, as essentially noted by~\citet{audibert2014regret}.
Now, the INF strategy provides a distribution-free regret bound of order $\sqrt{KT}$
in case of a known payoff range $[m,M]$. Up to some technical issues, which we could solve,
it may be extended to provide a similar scale-free distribution regret bound, which is optimal
as it does not contain any superfluous $\sqrt{\ln K}$ factor.
The exact statement (Theorem~\ref{thm:adaptive_lower_bound}) proved in Appendix~\ref{app:proofs_tsallis}
is the following: AdaFTRL with $1/2$--Tsallis entropy relying on an upper bound $M$ on the payoffs
ensures that for all $m \in \R$ with $m \leq M$,
for all oblivious individual sequences $y_1,y_2,\ldots$ in $[m,M]^K$,
for all $T \geq 1$,
\[
R_T(y_{1:T}) \leq 4 (M-m)\sqrt{KT} + 2(M-m)\,.
\]
We now give a high-level idea of the technical issues that were solved to obtain the bound above.
We consider estimates $\what{y}_{t, a}$ obtained from~\eqref{eq:estscheme} by replacing the constant $C$ therein by the
known upper end~$M$. We however could not simply derive the regret bound from
some generic full-information regret guarantee for AdaFTRL with $1/2$--Tsallis entropy, as
to the best of our knowledge, there are no meaningful full-information regret bounds for Tsallis entropy
in the first place, and as these would anyway scale with the effective range of the estimates.
We instead provide a more careful analysis exploiting special properties of the estimates, namely, that 
$\what{y}_{t, a} = M$ for all $a \ne A_t$ and $\what{y}_{t, A_t} \leq M$.

\paragraph{Disclaimer, continued.}
We were unable so far to provide a non-trivial distribution-dependent regret
bound for our strategy AdaFTRL with $1/2$--Tsallis entropy. Note that there exist
$\mathcal{O}(\ln T)$ bounds for FTRL with $1/2$--Tsallis entropy, i.e.,
with a different tuning of the learning rates (namely, $\eta_t$ of order $1/\sqrt{t}$, but then, the scale-free distribution-free
guarantees are lost); see \citet{zimmert2018optimal}.
We would have liked to prove such a $\mathcal{O}(\ln T)$ scale-free distribution-dependent regret
bound for AdaFTRL with $1/2$--Tsallis entropy (or even achieve a more modest aim like a poly-logarithmic bound),
as this seems possible and would have shown with certainty that the trade-off imposed
by~Theorem~\ref{thm:adaptive_lower_bound} does not hold anymore when the upper end~$M$ on the payoff range
is known.
The techniques of \citet{seldin2017improved}, which consist in a precise tuning of the extra-exploration in their
variant of the Exp3 algorithm of~\citet{AuCBFrSc02} together with a gap estimation scheme, or
the ones of \citet{zimmert2018optimal} might be helpful to that end. We leave this problem
for future research.

\subsection{AdaFTRL for Full Information (Reminder of Known Results)}
\label{sec:AdaFTRL-full}

To avoid confusion with the notation used in the main body of the paper,
we first describe the considered setting of prediction of oblivious individual sequences with full information.

\paragraph{Full-information setting.}
The game between the player and the environment is actually the same as the one described in Section~\ref{sec:setting-oblivious},
except that the player observes at each step the entire payoff vector, not just the obtained payoff.
More formally (and with a different piece of notation $z$ instead of $y$, to better distinguish the two settings),
the environment first picks a sequence of payoff vectors $z_t \in \R^K$, for all $t \geq 1$.
Then, in a sequential manner, at every time step $t$, the player picks an action $A_t$,
distributed according to a probability $p_t$ over the action set $[K]$,
obtains the payoff $z_{t, A_t}$, and observes the entire vector $z_t$ (i.e., also
the payoffs $z_{t,a}$ corresponding to the actions $a \neq A_t$).

In the sequel, we denote by $\cS$ the simplex of probability distributions over $[K]$
and we use the short-hand notation, for $p \in \cS$ and $z \in \R^K$,
\[
\dotp{p}{z} = \sum_{a \in [K]} p_a z_a\,.
\]

\paragraph{FTRL (follow-the-regularized-leader).}
The FTRL method consists in choosing $p_t$ according to
\[
p_t \in \argmin_{p \in \mathcal S : F(p) < +\infty} \left\{\frac{F(p)}{\eta_t} - \sum_{s =1}^{t-1} \dotp{p}{z_s} \right\},
\]
where $F : \R^K \to \R \cup \{+ \infty\}$ is a convex function, called the regularizer,
and $\eta_t$ is a non-negative learning rate in $(0, + \infty]$, which may depend on past observations.
The condition $F(p) < +\infty$ will always be satisfied for some $p \in \cS$ by the considered regularizers (see below)
and is only meant to avoid the undefined $+\infty/+\infty$ in the case $\eta_t = +\infty$. For the sake of concision
we will however omit it in the sequel.

Let us give a succint account of the convex analysis results we use here,
following the exposition of \citet[Chapter 26]{lattimore_szepesvari_2020}. Using their terminology,
the domain $\dom L$ of a convex function $L : \R^K \to \R \cup \{+ \infty\}$
is the set $\{ x \in \R^K : L(x) < +\infty \}$ of those points where it takes finite values.
A convex function $L : \R^K \to \R \cup \{+ \infty\}$ is said to be Legendre if the interior of its domain $\interior( \dom L)$ is non-empty,
if $L$ is strictly convex and differentiable on $\interior( \dom L)$, and if its gradient $\nabla L$
blows up on the boundary of $\dom L$. The minimizers of
Legendre functions may be seen to satisfy the following properties.

\begin{proposition}[Special case of {\citealp[Proposition~26.14]{lattimore_szepesvari_2020}}]
\label{prop:opt_condition}
Let $L$ be a Legendre function and $A \subseteq \R^d$ be a convex set that intersects $\interior(\dom L)$.
Then $L$ possesses a unique minimizer $x^\star$ over $A$, which belongs to $\interior (\dom L )$,
therefore ensuring that $L$ is differentiable at $x^\star$. Furthermore,
\[
\forall x \in A \cap \dom L, \qquad \dotp{\nabla L(x^\star)}{x-x^\star} \geq 0\,.
\]
\end{proposition}
\smallskip

Finally, for $x, y \in \R^d$, if $F : \R^K \to \R \cup \{+ \infty\}$ is differentiable at $y$,
we define the Bregman divergence between $x$ and $y$ as
\begin{equation}
B_F(x, y) = F(x) - F(y) - \dotp{\nabla F(y)}{x- y}\,;
\end{equation}
when $F$ is convex, we have $B_F(x, y) \geq 0$ for all $x \in \R^d$.

We are now ready to state our first reminder, which is a classical
regret bound for FTRL (see, e.g., \citealp[Chapter~28, Exercise 28.12]{lattimore_szepesvari_2020}
for references, and \citealp{mcmahan2017survey} for more general versions). It involves
the diameter $D_F$ of the action set (the $K$--dimensional simplex $\cS$ in our case):
\[
D_F = \max_{p, q \in \mathcal S} \big\{ F(p) - F(q) \big\}\,.
\]

\begin{reminder}[Generic full-information FTRL bound over the simplex]
\label{rem:ftrl}
The FTRL method with a Legendre regularizer $F$ (of finite diameter $D_F$) and
with any rule for picking the learning rates so that they form a non-increasing sequence
satisfies the following guarantee:
for all sequences $z_1,\,z_2,\ldots$ of vector payoffs in $\R^K$,
the regret is bounded by
\begin{align}
\nonumber
\max_{a \in [K]} \sum_{t = 1}^T z_{t, a} - \sum_{t = 1}^T \dotp{p_t}{z_t}
\leq \frac{D_F}{\eta_T}
& + \sum_{t = 1}^{T-1} \left(\dotp{p_t - p_{t+1}}{ - z_t}  - \frac{B_{F}(p_{t+1}, p_t)}{\eta_t}\right) \\
\label{eq:ftrl_generic_bound}
& + \left(\dotp{p_T - p^\star}{ - z_T}  - \frac{B_{F}(p^\star, p_T)}{\eta_T}\right), \\
\nonumber
\mbox{where} \qquad p^\star \in \argmax_{p \in \mathcal S} \sum_{t=1}^T \dotp{p}{z_t} &
\end{align}
and where the regret bound is well defined, thanks to the following observations and conventions:
for rounds $t \geq 1$ where $\eta_t < +\infty$, the function $F$ is indeed differentiable at $p_t$
so that $B_{F}(p_{t+1}, p_t)$ is well defined;
for rounds $t \geq 1$ where $\eta_t = +\infty$, we set $B_{F}(p_{t+1}, p_t)/\eta_t = 0$
irrespectively of the fact whether $F$ is differentiable at $p_t$.
\end{reminder}

\bpc{of Reminder~\ref{rem:ftrl}}
Denote by $S_t$ the cumulative vector payoff up to time $t \geq 1$.
Fix $T \geq 1$.
For the sake of concision of the equations, we define $p_{T+1} = p^\star$,
which is a Dirac mass at some arm (that is, $p_{T+1}$ is not given by FTRL).
The regret can therefore be rewritten as
\begin{align*}
\max_{a \in [K]} \sum_{t = 1}^T z_{t, a} - \sum_{t = 1}^T \dotp{p_t}{z_t}
& = \max_{p \in \mathcal S} \sum_{t = 1}^T \dotp{p}{z_{t}}
- \sum_{t = 1}^T \dotp{p_t}{z_t} \\
& = \sum_{t = 1}^T \dotp{p_{T+1}}{z_{t}}
- \sum_{t = 1}^T \dotp{p_t}{z_t}
= \sum_{t = 1}^T \dotp{p_t - p_{T+1}}{-z_t}\,.
\end{align*}
By summation by parts,
\begin{align}
\nonumber
\MoveEqLeft \sum_{t = 1}^T \dotp{p_t - p_{T+1}}{-z_t} \\
\nonumber
& = \sum_{t = 1}^T \sum_{s = t}^T \dotp{p_s - p_{s+1}}{-z_t}
= \sum_{s = 1}^T \sum_{t= 1}^s \dotp{p_s - p_{s+1}}{-z_t}
= \sum_{s = 1}^T \dotp{p_s - p_{s+1}}{-S_s} \\
\label{eq:decregFTRL}
& = \sum_{t = 1}^T \dotp{p_t - p_{t+1}}{-z_t} + \sum_{t = 1}^T \dotp{p_t - p_{t+1}}{-S_{t-1}}\,.
\end{align}
If $\eta_t < + \infty$, then by the optimality condition from Proposition~\ref{prop:opt_condition} applied
to the Legendre function $L : x \mapsto \eta_t^{-1} F(x) - \dotp{S_{t-1}}{x} $, we know that $L$ thus $F$ are differentiable at $p_t$ and that
\begin{align*}
& \dotp{\eta_t^{-1}\nabla F(p_t) - S_{t-1}}{p_{t+1}-p_t} \geq 0\,, \\
\mbox{that is,} \qquad
& \dotp{p_t -p_{t+1}}{- S_{t-1}} \leq \dotp{\eta_t^{-1}\nabla F(p_t)}{p_{t+1} - p_t} \, .
\end{align*}
If $\eta_t = + \infty$, the previous inequality holds too,
as by definition of $p_t$, we have $\dotp{p_t - p_{t+1}}{-S_{t-1}} \leq 0$
and as we set by convention $\eta_t^{-1} \nabla F(p_t) = 0$ regardless of whether $F$ is differentiable at $p_t$ or not.
Substituting in~\eqref{eq:decregFTRL}, we proved so far
\begin{equation}
\sum_{t = 1}^T \dotp{p_t - p_{T+1}}{-z_t}
\leq \sum_{t = 1}^T \dotp{p_t - p_{t+1}}{-z_t}
+ \dotp{\eta_t^{-1}\nabla F(p_t)}{p_{t+1} - p_t}\,.
\end{equation}
This inequality can be rewritten in terms of Bregman divergences:
\[
\sum_{t = 1}^T \dotp{p_t - p^\star}{-z_t}
\leq \sum_{t = 1}^T \left( \dotp{p_t - p_{t+1}}{-z_t} - \frac{B_{F}(p_{t+1}, p_t)}{\eta_t} \right)
+ \sum_{t = 1}^T\frac{F(p_{t+1}) - F(p_{t})}{\eta_{t}}
\]
We now upper bound the second sum in the right-hand side:
again by summation by parts, with the convention $\eta_0 = +\infty$ and $1/\eta_0 = 0$:
\begin{align*}
& \sum_{t = 1}^T\frac{F(p_{t+1}) - F(p_{t})}{\eta_{t}}  = \sum_{t = 1}^T\big(F(p_{t+1}) - F(p_{t})\big) \sum_{s = 1}^{t} \left( \frac{1}{\eta_{s}}  - \frac{1}{\eta_{s-1}} \right) \\
= & \sum_{s=1}^{T} \sum_{t = s}^T \big(F(p_{t+1}) - F(p_{t})\big) \left( \frac{1}{\eta_{s}}  - \frac{1}{\eta_{s-1}} \right) = \sum_{s=1}^T
\big(\underbrace{F(p_{T+1}) - F(p_s)}_{\leq D_F} \big)
\biggl( \underbrace{\frac{1}{\eta_{s}}  - \frac{1}{\eta_{s-1}}}_{\geq 0} \biggr)
\leq \frac{D_F}{\eta_T}\,,
\end{align*}
where the final equality is obtained by a telescoping sum,
using that the sequence of learning rates is non-increasing.
\end{proof}

\paragraph{AdaFTRL, an adaptive version of FTRL.}
The AdaFTRL approach consists in tuning the learning rate in a way that scales with the observed data.
More precisely, it relies on a quantity called the (generalized) mixability gap, which naturally appears
as an upper bound on the summands in the FTRL bound of Reminder~\ref{rem:ftrl}:
\begin{equation}
\label{eq:def_gen_mix_gap}
\delta^{F}_t \eqdef \max_{p \in \mathcal S} \left\{ \dotp{p_t - p}{ - z_t } - \frac{B_F(p, p_t)}{\eta_t} \right\} \geq 0\,.
\end{equation}
That mixability gaps are always nonnegative can be seen by taking
$p = p_t$ in the definition.
We may further upper bound \eqref{eq:ftrl_generic_bound} when it holds by using this mixability gap:
\begin{equation}
\label{eq:ada_ftrl_intermediate}
\max_{a \in [K]} \sum_{t = 1}^T z_{t, a} - \sum_{t = 1}^T \dotp{p_t}{z_t} \leq \frac{D_F}{\eta_T} + \sum_{t = 1}^T \delta_T^F\,.
\end{equation}
The AdaFTRL learning rate balances the two terms in the above regret bound by taking
\begin{equation}
\label{eq:adaftrl-rates}
\eta_t = D_F \Bigg/  \, \sum_{s = 1}^{t-1} \delta_s^F  \quad \in (0,  + \infty]
\end{equation}
Note that this rule for picking learning rates indeed leads to non-increasing sequences thereof, as the mixability gaps are non-negative.
We summarize the discussion above in the theorem stated next,
from which subsequent (closed-from) regret bounds will be derived by
using the specific properties of the regularizer $F$ at hand to upper bound the mixability gaps.

\begin{theorem}[AdaFTRL tool box]\label{thm:adaftrl_toolbox}
Under the assumptions of Reminder~\ref{rem:ftrl} and with its conventions,
the regret of the FTRL method based on the learning rates~\eqref{eq:adaftrl-rates}
satisfies
\begin{equation}\label{eq:gen_ftrl_mixsum}
\max_{a \in [K]} \sum_{t = 1}^T z_{t, a} - \sum_{t = 1}^T \dotp{p_t}{z_t}
\leq 2 \sum_{t = 1}^T \delta_t^F
\end{equation}
where, moreover,
\begin{equation}\label{eq:gen_sq_mix_gaps}
\left( \sum_{t = 1}^T \delta_t^F \right)^{\!\! 2}
= 2 D_F \sum_{t = 1}^T \frac{\delta_t^F}{\eta_t} +   \sum_{t = 1}^T \bigl( \delta_t^F \bigr)^2\,.
\end{equation}
\end{theorem}

\bp
Inequality~\eqref{eq:gen_ftrl_mixsum} follows
from~\eqref{eq:ada_ftrl_intermediate} and~\eqref{eq:adaftrl-rates}.
The equality~\eqref{eq:gen_sq_mix_gaps} is obtained by expanding the squared sum,
\[
\left( \sum_{t = 1}^T \delta_t^F \right)^{\!\! 2} = \sum_{t = 1}^T \bigl( \delta_t^F \bigr)^2
+ 2 \sum_{t = 1}^T \sum_{s = 1}^{t-1} \delta_t^F \delta_s^{F} = \sum_{t = 1}^T (\delta_t^F)^2 + 2 \sum_{t =1}^T \delta_t^F \frac{D_F}{\eta_t}
\]
where the final equality is obtained by substituting the definition~\eqref{eq:adaftrl-rates} of $\eta_t$.
\end{proof}

\subsection{AdaHedge for Full Information (Reminder of Known Results)}
\label{app:ada_hedge_proof}

The content of this section is extracted from various sources, out of which the most important is~\citet{koolen_2016}.
We claim no novelty. This section recalls how the bound for AdaHedge (Reminder~\ref{rem:adahedge},
for which a direct proof was provided by~\citealp{de2014follow}) can also be seen as a special case of the results of
Section~\ref{sec:AdaFTRL-full}.

It is well-known (see \citealp{FSSW97,KW99,Aud09}),
and can be found again by a simple optimization under a linear constraint,
that the Hedge weight update corresponds to FTRL with the negentropy as a regularizer:
\[
H_{\negent}(p) = \sum_{a = 1}^K p_a \log p_a\,,
\]
with value $+\infty$ whenever $p_a = 0$ for some $a \in [K]$. That is,
\begin{multline}
\label{eq:defadahedgefullinfo}
\argmin_{p \in \mathcal S} \left\{\frac{H_{\negent}(p)}{\eta_t} - \sum_{s =1}^{t-1} \dotp{p}{z_s} \right\} = \{ p_t \} \\
\mbox{with} \qquad
p_{t,a} = \exp \! \left( \eta_{t} \sum_{s = 1}^{t-1} z_{a,s} \right)
\Bigg/ \sum_{k=1}^K \exp \! \left( \eta_{t} \sum_{s = 1}^{t-1} z_{k,s} \right).
\end{multline}
Straightforward calculation show that the regularizer $H_{\negent}$ is indeed Legendre
(see~\citealp{lattimore_szepesvari_2020}, Example 26.11)
and the $H_{\negent}$--diameter of the simplex equals $D_{H_{\negent}} = \log K$. Reminder~\ref{rem:ftrl}
and Theorem~\ref{thm:adaftrl_toolbox} can therefore be applied.

AdaHedge is exactly AdaFTRL with $H_{\negent}$ as a regularizer. Indeed,
the mixability gap~\eqref{eq:def_gen_mix_gap} can be computed in closed form (as noted
by~\citealp[Lemma 5]{reid15generalized}) and reads in this case:
\begin{equation}
\label{eq:def_delta_t}
\delta^{\negent}_t =\left\{
\begin{split}
& - \dotp{p_t}{z_t} + \eta_t^{-1} \log \!\left( \sum_{  a =1}^K p_{t, a} \e^{\eta_t z_{t,a}} \right) \quad & \text{if } \eta_t < + \infty,  \\
& - \dotp{p_t}{z_t} +  \max_{a \in [K]} z_{t,a} \quad & \text{if } \eta_t = + \infty.
\end{split}\right.
\end{equation}

\bpc{of the rewriting~\eqref{eq:def_delta_t}}
When $\eta_t = +\infty$, the mixability gap equals, by definition,
\[
\delta^{F}_t = \max_{p \in \mathcal S} \bigl\{ \dotp{p_t - p}{ - z_t } \bigr\}
= - \dotp{p_t}{z_t} +  \max_{p \in \mathcal S} \dotp{p}{z_t }
= - \dotp{p_t}{z_t} +  \max_{a \in [K]} z_{t,a}\,.
\]
For the case $\eta_t < +\infty$,
the following formula, which is at the heart of the closed-form formula for the Hedge updates~\eqref{eq:defadahedgefullinfo},
will be useful: for any $S \in \R^d$,
\begin{equation}\label{eq:negent_dual}
\min_{p \in \mathcal S} \Big\{  H_{\negent}(p)  -\dotp{p}{S}  \Big\}
= \sum_{i = 1}^K \frac{\e^{S_i}}{\sum_{j = 1}^K \e^{S_j}} \left(\ln \Bigg( \frac{\e^{S_i}}{\sum_{j = 1}^K \e^{S_j}  } \Bigg) - S_i \right)
= - \ln \Bigg( \sum_{i = 1}^K \e^{S_i} \Bigg)\,.
\end{equation}
When $\eta_t < +\infty$, Equation~\eqref{eq:defadahedgefullinfo} shows that $p_t$
lies in the interior $\interior(\cS)$ of $\cS$.
The Bregman divergence at hand in the definition~\eqref{eq:def_gen_mix_gap} of the mixability gaps
may be simplified into
\[
B_F(p, p_t)
= H_{\negent}(p) - H_{\negent}(p_t) - \dotp{\nabla H_{\negent}(p_t)}{p - p_t}
= H_{\negent}(p) - \dotp{\nabla H_{\negent}(p_t)}{p} + 1\,,
\]
where the second inequality holds by
taking into account the fact that $H_{\negent}$ is twice differentiable at any $p \in \interior(\cS)$, with
\[
\nabla H_{\negent}(p) = \big(1 + \ln p_i\big)_{i \in [K]}
\qquad \mbox{so that} \qquad
\dotp{\nabla H_{\negent}(p)}{p} = 1+ \sum_{i = 1}^K p_i \log p_i = 1 + H_{\negent}(p)\,.
\]
The mixability gaps can therefore be rewritten
\begin{align*}
\delta^{F}_t & = \max_{p \in \mathcal S} \left\{ \dotp{p_t - p}{ - z_t } - \frac{B_F(p, p_t)}{\eta_t} \right\} \\
& = - \dotp{p_t}{z_t} - \frac{1}{\eta_t} + \frac{1}{\eta_t} \max_{p \in \mathcal S} \bigl\{ \eta_t \dotp{p}{z_t} -
H_{\negent}(p) + \dotp{\nabla H_{\negent}(p_t)}{p} \bigr\} \\
& = - \dotp{p_t}{z_t} - \frac{1}{\eta_t} - \frac{1}{\eta_t} \min_{p \in \mathcal S} \Bigl\{ H_{\negent}(p) -
\dotpb{p}{\eta_t z_t + \nabla H_{\negent}(p_t)} \Bigr\}
\end{align*}
Now by \eqref{eq:negent_dual}, specialized with $S = \eta_t z_t + \nabla H_{\negent}(p_t)$,
we can compute the value of the minimum:
\[
\min_{p \in \mathcal S} \Bigl\{ H_{\negent}(p) -
\dotpb{p}{\eta_t z_t + \nabla H_{\negent}(p_t)} \Bigr\}
= - \ln \left( \sum_{ i =1}^K \e^{\eta_t z_i + 1 + \ln p_{i}} \right)
= -1 - \ln \left( \sum_{ i =1}^K p_{i} \e^{\eta_t z_i} \right).
\]
Collecting all equalities together concludes the proof.
\end{proof}

Reminder~\ref{rem:adahedge} is thus a special case of the following bound.

\begin{theorem}[See Lemma~3 and Theorem~6 of~\citealp{de2014follow}]
\label{thm:adahedge}
For all sequences of payoffs $z_{t,a}$ lying in some bounded real-valued interval,
denoted by $[b,B]$, for all $T \geq 1$, the regret of the AdaHedge algorithm with full information,
as defined by~\eqref{eq:defadahedgefullinfo} and~\eqref{eq:def_delta_t}, satisfies
\begin{multline*}
\max_{k \in [K]} \sum_{t=1}^T z_{t,k} -
\sum_{t=1}^T \sum_{a=1}^K p_{t,a} \, z_{t,a} \leq 2 \sum_{t=1}^T \delta^{\negent}_t \\
\mbox{where} \qquad \sum_{t=1}^T \delta^{\negent}_t
\leq \sqrt{\sum_{t=1}^T \sum_{a=1}^K p_{t,a} \! \left( z_{t,a} - \sum_{k \in [K]} q_{t,k} \, z_{t,k} \right)^{\!\! 2} \ln K}
+ (B-b) \left( 1 + \frac{2}{3}\ln K \right),
\end{multline*}
and AdaHedge does not require the knowledge of $[b,B]$ to achieve this bound.
\end{theorem}

The quantities
\[
v_t \eqdef \sum_{a=1}^K p_{t,a} \! \left( z_{t,a} - \sum_{k \in [K]} q_{t,k} \, z_{t,k} \right)^{\!\! 2}
\]
in the bound correspond to the variance of the random variables taking values $z_{t,a}$ with probability $p_{t,a}$;
the variational formula for variances indicates that
\[
\sum_{a=1}^K p_{t,a} \! \left( z_{t,a} - \sum_{k \in [K]} q_{t,k} \, z_{t,k} \right)^{\!\! 2}
= \min_{c \in \R} \sum_{a=1}^K p_{t,a} \bigl( z_{t,a} - c \bigr)^2\,,
\]
which entails the final bound given as a note in the statement of Reminder~\ref{rem:adahedge}.

The following formulation of Bernstein's inequality will be useful in the proof of Theorem~\ref{thm:adahedge}.
\begin{lemma}[Bernstein's inequality tailored to our needs]
\label{lm:Bern-Hedi}
Let $X$ be a random variable in $[0, 1]$, with variance denoting by $\Var(X)$.
Then for all $\eta > 0$,
\[
\frac{\log \Bigl( \E \bigl[ \e^{\eta (X - \E[X]) } \bigr] \Bigr)}{\eta^2}
\leq \frac{1}{2} \Var(X) +
\frac{1}{3} \, \frac{\log \Bigl( \E \bigl[ \e^{\eta (X - \E[X]) } \bigr] \Bigr)}{\eta}\,.
\]
\end{lemma}

\bp
Denote by $ \psi_X(\eta) = \log \bigl( \E \bigl[ \e^{\eta (X - \E[X]) } \bigr]$ the log-moment generating function of $X$.
A version of Bernstein's inequality with an appropriate control of the moments
(as stated by \citealp[Section~2.2.3]{Mas03} and applied to $X$ with $c = 1/3$) indicates that for all $\eta \in (0, 3)$,
\[
\Big( 1 - \frac{\eta}{3} \Big) \psi_X(\eta) \leq \frac{\eta^2}{2} \Var(X)\,.
\]
Actually, this inequality also holds for $\eta \geq 3$ as its left-hand side is non-positive while its right-hand side is nonnegative.
The claimed result is derived by rearraging the terms
\[
\psi_X(\eta) \leq \frac{\eta^2}{2} \Var(X) + \frac{\eta}{3} \psi_X(\eta)
\]
and by dividing both sides by $\eta^2$.
\end{proof}

\bpc{of Theorem~\ref{thm:adahedge}}
We apply Theorem~\ref{thm:adaftrl_toolbox}.
To that end, we first bound the mixability gaps. The rewriting~\eqref{eq:def_delta_t}
(and Jensen's inequality) directly shows that $0 \leq \delta_t^{\negent} \leq B-b$.
We may also prove the bound
\begin{equation}
\label{eq:BernGap}
\frac{\delta_t^{\negent}}{\eta_t} \leq  \frac{v_t}{2} + \frac{1}{3}(B-b)\delta_t^{\negent}\,.
\end{equation}
It suffices to do so for $\eta_t < +\infty$.
Consider the random variable $X$ taking values
$(z_{t,a} - b)/(B-b)$ with probability $p_{t,a}$, for $a \in \{1,\ldots,K\}$.
The mixability gap can be rewritten as
\[
\delta_t^{\negent} = \frac{1}{\eta_t} \psi_X\bigl( \eta_t(B-b)\bigr)
\]
with the notation of the proof of Lemma~\ref{lm:Bern-Hedi}.
The variance of $X$ equals $v_t/(B-b)^2$.
Lemma~\ref{lm:Bern-Hedi} with $\eta = \eta_t(B-b)$ yields
\[
\frac{\delta_t^{\negent}}{\eta_t (B-b)^2} \leq \frac{v_t}{2(B-b)^2} + \frac{\delta_t^{\negent}}{3(B-b)}\,.
\]
from which we obtain~\eqref{eq:BernGap} by rearranging.

From~\eqref{eq:gen_sq_mix_gaps} and~\eqref{eq:BernGap}, we deduce, together with
the bound $(\delta_t^{\negent})^2 \leq (B-b)\delta_t^{\negent}$, that
\[
\left(\sum_{t = 1}^T \delta_t^{\negent} \right)^{\!\! 2}
\leq (\log K) \sum_{t = 1}^T v_t + (B-b)\left( \frac{2}{3} \log K + 1  \right) \sum_{t = 1}^T \delta_t^{\negent}\,.
\]
Therefore, using the fact that $x^2 \leq a + bx$ implies $x \leq \sqrt{a} + b$ for all $a, b, x \geq 0$,
\[
\sum_{t = 1}^T \delta_t^{\negent} \leq \sqrt{\ln K \sum_{t = 1}^T v_t} + (B-b)\left( \frac{2}{3} \log K + 1  \right),
\]
which, thanks to~\eqref{eq:gen_ftrl_mixsum}, concludes the proof of Theorem~\ref{thm:adahedge}.
\end{proof}

\subsection{AdaHedge with Known Payoff Upper Bound $M$ (Application of Section~\ref{app:ada_hedge_proof})}
\label{app:proofs_exp3}

We show how to obtain a scale-free distribution-free regret bound of order $(M-m)\sqrt{K T \ln K}$
with no extra-exploration (including no initial exploration)
when an upper bound $M$ on the payoffs is given to the player.
We consider Algorithm~\ref{algo:adahedgeM}, where no mixing takes place (unlike
in Algorithm~\ref{algo:adahedge}) and where the probability distributions $p_t$
are directly computed via an AdaHedge update (no need for intermediate
probabilities $q_t$). Note also that we use the estimates~\eqref{eq:estscheme}
with the choice $C_t = M$, that is,
\begin{equation}
\label{eq:loss_estimates_known_M}
\what{y}_{t, a} = \frac{y_{t, a}-M}{p_{t, a}} \1{A_t = a} + M\,.
\end{equation}
The following observation is key in the analysis below:
$\what{y}_{t, a} = M$ for all $a \ne A_t$ and $\what{y}_{t, A_t} \leq M$.
We will also use, as in the proof of Theorem~\ref{th:adv_talpha_bound},
\[
\sum_{a=1}^K p_{t,a}\,\wh{y}_{t,a} = y_{t, A_t}\,.
\]

\begin{figure}[h]
\renewcommand\footnoterule{}
\begin{algorithm}[H]
\begin{algorithmic}[1]
\caption{\label{algo:adahedgeM} AdaHedge for $K$--armed bandits, when an upper bound on the payoffs is given}
\label{alg:exp3_adahedgeM}
\STATE \textbf{Input:} an upper bound $M$ on the payoffs
\STATE \textbf{AdaHedge initialization:} $\eta_{1} = +\infty$ and $p_{1} = (1/K,\ldots,1/K)$
\FOR{rounds $t = 1,\,2,\, \dots$}
\STATE Draw an arm $A_t \sim p_t$ (independently at random according to the distribution $p_t$)
\STATE Get and observe the payoff $y_{t, A_t}$
\STATE Compute the estimates of all payoffs
\[
\what y_{t, a} = \frac{y_{t, a}-M}{p_{t, a}} \1{A_t = a} + M
\]
\STATE Compute the mixability gap $\delta_t$ based on the distribution $p_t$ and on these estimates:
\[
\delta_t =\left\{
\begin{split}
& - \sum_{a=1}^K p_{t,a}\,\wh{y}_{t,a} + \frac{1}{\eta_t} \log \Biggl( \sum_{a =1}^K p_{t, a} \e^{\eta_t \wh{y}_{t,a}} \Biggr)
\quad & \text{if } \eta_t < + \infty \\
& - \sum_{a=1}^K p_{t,a}\,\wh{y}_{t,a} + \max_{a \in [K]} \wh{y}_{t,a} \quad & \text{if } \eta_t = + \infty
\end{split}\right.
\]
\STATE Compute the learning rate $\displaystyle{\eta_{t+1} = \Biggl( \sum_{s = 1}^t \delta_s \Biggr)^{-1}} \ln K$
\STATE Define $p_{t+1}$ component-wise as \vspace{-.3cm}
\[
~ \hspace{2.5cm} p_{t+1,a} = \exp \! \left( \eta_{t+1} \sum_{s = 1}^{t} \wh{y}_{a,s} \right)
\Bigg/ \sum_{k=1}^K \exp \! \left( \eta_{t+1} \sum_{s = 1}^{t} \wh{y}_{k,s} \right) \vspace{-.3cm}
\]
\ENDFOR
\end{algorithmic}
\end{algorithm}
\end{figure}

The performance bound for this simpler algorithm is stated next.
\begin{theorem}
\label{thm:adv_any_exploration--simpler}
AdaHedge for $K$--armed bandits relying on an upper bound $M$
on the payoffs (Algorithm~\ref{algo:adahedgeM})
ensures that for all $m \in \R$ with $m \leq M$,
for all oblivious individual sequences $y_1,y_2,\ldots$ in $[m,M]^K$,
for all $T \geq 1$,
\[
R_T(y_{1:T}) \leq 2 (M-m)\sqrt{KT \log K} + 2(M-m)\,.
\]
\end{theorem}

The main technical difference with respect to the analysis of Algorithm~\ref{algo:adahedge}
is that the mixability gaps are directly bounded by the range $M-m$.
We no longer need to artificially control the size of the estimates (which we did via extra-exploration) to
get, in turn, a control of the mixability gaps.

\begin{lemma}[Improved mixability gap bound]
\label{lem:adahedge_up_bound_mix_gap}
The mixability gaps of AdaHedge for $K$--armed bandits relying on an upper bound $M$
on the payoffs (Algorithm~\ref{algo:adahedgeM})
are bounded, for all $m \in \R$ with $m \leq M$,
for all oblivious individual sequences $y_1,y_2,\ldots$ in $[m,M]^K$,
for all $t \geq 1$, by
\[
0 \leq \delta_t \leq  M-m  \qquad \mbox{and} \qquad \frac{\delta_t}{\eta_t} \leq \frac{1}{2} \,  p_{t, A_t}^{-1}(M - y_{t,A_t})^2  .
\]
\end{lemma}

\bp
The fact that $\delta_t  \geq 0$ holds by definition of the gaps and Jensen's inequality.
For $\delta_t \leq M-m$,
the observations after~\eqref{eq:loss_estimates_known_M} indicate that when $\eta_t = +\infty$,
\[
\delta_t = - \sum_{a=1}^K p_{t,a}\,\wh{y}_{t,a} + \max_{a \in [K]} \wh{y}_{t,a}
= M - \wh{y}_{t,A_t}\,,
\]
while for $\eta_t < +\infty$,
\begin{align*}
\delta_t & = - y_{t, A_t} + \frac{1}{\eta_t} \log \left( (1 - p_{t, A_t}) \e^{\eta_t M}
+ p_{t, A_t} \e^{\eta_t M} \e^{\eta_t (y_{t, A_t}-M) / p_{t, A_t}}\right) \\
& \leq M - y_{t, A_t} + \frac{1}{\eta_t} \log \Bigl( (1 - p_{t, A_t})
+ p_{t, A_t} \underbrace{\e^{\eta_t (y_{t, A_t}-M) / p_{t, A_t}}}_{\leq 1} \Bigr)\,,
\end{align*}
which entails $\delta_t \leq M - y_{t, A_t} \leq M-m$.

Furthermore, in the case $\eta_t < +\infty$,
using the inequality $\e^{-x} \leq  1 - x +  x^2/2$ valid for $x \geq 0$, followed by the inequality
$\log(1 + u) \leq u$, valid for all $u > -1$, we get
\[
\delta_t
\leq M - \what y_{t, A_t}  + \frac{1}{\eta_t} \log \! \biggl(
\underbrace{1 - p_{A_t, t,} + p_{A_t, t}}_{=1} \underbrace{- \eta_t (M -  y_{t, A_t} ) + \eta_t^2 \frac{ (M -  y_{t, A_t} )^2}{2 p_{A_t, t}}}_{= u} \biggr)
\leq \eta_t \frac{(M -  y_{t, A_t} )^2}{2 p_{t, A_t}}\,.
\]
The second inequality is trivial in case $\eta_t = +\infty$, as $\delta_t/\eta_t = 0$.
\end{proof}

We are now ready to prove Theorem~\ref{thm:adv_any_exploration--simpler}. \medskip 

\bpc{of Theorem~\ref{thm:adv_any_exploration--simpler}}
As indicated in Section~\ref{app:ada_hedge_proof}, AdaHedge is a special case of AdaFTRL
and the bound of Theorem~\ref{thm:adaftrl_toolbox} is applicable.

Equation~\eqref{eq:gen_sq_mix_gaps} and Lemma~\ref{lem:adahedge_up_bound_mix_gap},
which entails in particular that $\delta_t^2 \leq (M-m) \delta_t$,
yield
\[
\bigg(\sum_{t = 1}^T \delta_t \bigg)^{\!\! 2} = 2 (\ln K) \sum_{t = 1}^T \frac{\delta_t}{\eta_t}
+ \sum_{t = 1}^T (\delta_t)^2 \leq (\ln K)  \sum_{t = 1}^T p_{t, A_t}^{-1} (M- y_{t, A_t})^2+ (M-m) \sum_{t=1}^T \delta_t\,,
\]
which, through the fact that $x^2 \leq a + bx$ implies $x \leq \sqrt{a} + b$ for all $a, b, x \geq 0$,
leads in turn to
\[
\sum_{t = 1}^T \delta_t \leq \sqrt{\sum_{t=1}^T p_{t,A_t}^{-1} \! \left( M - \wh{y}_{t,A_t} \right)^{ 2} \ln K} + (M-m)\,.
\]
Therefore, Equation~\eqref{eq:gen_ftrl_mixsum} guarantees that
\begin{equation}\label{eq:almostfin}
\max_{k \in [K]} \sum_{t=1}^T \wh{y}_{t,k} -
\sum_{t=1}^T \underbrace{\sum_{a=1}^K p_{t,a} \, \wh{y}_{t,a}}_{= y_{t,A_t}}
\leq
2\sqrt{\sum_{t=1}^T p_{t,A_t}^{-1} \! \left( M - \wh{y}_{t,A_t} \right)^{ 2} \ln K}
+ 2(M-m)\,.
\end{equation}
We conclude the proof by integrating the inequality above
and using Jensen's inequality, exactly as in the proof of Theorem~\ref{th:adv_talpha_bound}.
Indeed, Equation~\eqref{eq:exp-regret-max} therein indicates that
\[
R_T(y_{1:T})
= \max_{k \in [K]} \sum_{t=1}^T y_{t,k} - \E \! \left[ \sum_{t=1}^T  y_{t,A_t}  \right]
\leq \E \! \left[ \max_{k \in [K]} \sum_{t=1}^T \wh{y}_{t,k} - \sum_{t=1}^T y_{t,A_t} \right]
\]
and, by the same manipulations as in~\eqref{eq:Jensen-mainproofUB} and in the equation
that follows it,
\begin{align*}
	\E \! \left[ \sqrt{\sum_{t=1}^T p_{t,A_t}^{-1} \! \left( M - \wh{y}_{t,A_t} \right)^{ 2} \ln K} \right]
	& \leq \sqrt{ \E \! \left[ \sum_{t=1}^T p_{t,A_t}^{-1} \! \left( M - y_{t,A_t} \right)^{ 2} \ln K \right]} \\
	& = \sqrt{ \E \! \left[ \sum_{t=1}^T \sum_{a=1}^K \! \left( M - y_{t,a} \right)^{ 2} \ln K \right]}
	\leq (M-m)\sqrt{KT \ln K}
\end{align*}
The claimed result is obtained by collecting all bounds together.
\end{proof}

\subsection{AdaFTRL with Tsallis Entropy in the Case of a Known Payoff Upper Bound $M$}
\label{app:proofs_tsallis}

In this section we describe how the AdaHedge learning rate scheme can be used in the
FTRL framework with a different regularizer, namely Tsallis entropy, to improve the scale-free distribution-free regret bound into a bound of optimal order $(M-m)\sqrt{KT}$, i.e.,
without any superfluous $\sqrt{\ln K}$ factor.

\paragraph{Tsallis entropy.}
We focus on the (rescaled) $1/2$--Tsallis entropy, which is defined by
\[
H_{1/2}(p) = -\sum_{a = 1}^K 2 \sqrt{p_a}\,.
\]
This regularizer is Legendre over the domain $[0, +\infty)^K$ (see \citealp[Example~26.10]{lattimore_szepesvari_2020}).
Its diameter equals
\begin{equation}
	D_{H_{1/2}} = \max_{p \in \mathcal S} H_{1/2}(p) -  \min_{q \in \mathcal S} H_{1/2}(q)  =  -2 - \bigl( -2 \sqrt K \bigr) = 2 \bigl( \sqrt K -1 \bigr)\,,
\end{equation}
as for all $p \in \mathcal S$, we have (by concavity of the square root for the right-most inequality)
\[
1 \leq \sum_{ a = 1}^K  p_a \leq \sum_{ a = 1}^K \sqrt{p_a} \leq \sqrt K\,,
\]
where $1$ is achieved with $p = (1, 0, \dots, 0)$ and $\sqrt{K}$ with the uniform distribution.

The function $H_{1/2}$ is differentiable at all $q \in (0, +\infty)^K$, with $\nabla H_{1/2}(q) = \big(- 1 / \sqrt{q_a} \big)_{a \in [K]}$ .
The Bregman divergence associated with $H_{1/2}$ equals, for $p, q \in \mathcal S$ such that $q_{a} > 0$ for all $a$:
\begin{multline*}
	B_{H_{1/2}}(p, q)
	= -2 \sum_{ a = 1}^K \sqrt{p_a} + 2 \sum_{ a = 1}^K \sqrt{q_a} + \sum_{ a = 1}^K \frac{1}{\sqrt{q_a}} (p_a - q_a)  \\
	= -2 \sum_{a = 1}^K \frac{ \sqrt{p_a} - \sqrt{q_a}}{2 \sqrt{q_a}}  \Big( 2 \sqrt q_a - (\sqrt{p_a} + \sqrt{q_a}) \Big)
	= \sum_{  a =1}^K \frac{(\sqrt{p_a} - \sqrt{q_a})^2}{\sqrt {q_a}}\,.
\end{multline*}

\paragraph{AdaFTRL with $1/2$--Tsallis entropy.}
We consider FTRL with the $1/2$--Tsallis entropy on the estimated losses~\eqref{eq:loss_estimates_known_M}:
\[
p_t \in
\argmin_{p \in \mathcal S} \left\{\frac{H_{1/2}(p)}{\eta_t} - \sum_{s =1}^{t-1} \dotp{p}{\what y_{s}} \right\}
= \argmin_{p \in \mathcal S}  \left\{ - \frac{1}{\eta_t} \sum_{a = 1}^K 2\sqrt{p_a}
- \sum_{a = 1}^K p_{a} \sum_{s =1}^{t-1} \what y_{s, a} \right\}.
\]

FTRL with the $1/2$--Tsallis entropy was essentially introduced by~\citet{audibert2009Minimax} to
get rid of a $\sqrt{\ln K }$ factor in the distribution-free regret bound of $K$--armed adversarial bandits
(with known payoff range). It was later noted by~\citet{audibert2014regret} that it actually is an instance of
mirror descent with Tsallis entropy as a regularizer.
More recently, \citet{zimmert2018optimal} showed that this regularizer can obtain quasi-optimal regret bounds
for both stochastic and adversarial rewards.

We more precisely consider AdaFTRL with the $1/2$--Tsallis,
that is, we compute the learning rates $\eta_t$ based on the mixability
gaps~\eqref{eq:def_gen_mix_gap}; see Algorithm~\ref{alg:tsallis_ftrl}.
We denote by $\delta_t^{\ts}$ the mixability gaps~\eqref{eq:def_gen_mix_gap}.

\paragraph{On the implementation.}
For Tsallis entropy, the optimization problems involved in the computation of the updates $p_t$ and
of the mixability gaps~$\delta_t^{\ts}$ admit a (semi-)explicit formula.
Indeed, $p_t$ can be computed thanks to the formula, for all $z \in \R^K$,
\begin{equation}\label{eq:opt_pb}
\argmin_{p \in \mathcal S} \bigl\{  H_{1/2}(p) - \dotp{p}{z} \bigr\}
= \argmax_{p \in \mathcal S} \left\{  \dotp{p}{z} + \sum_{a = 1}^K 2 \sqrt{p_a} \right\}
= \left(\frac{1}{\bigl( c(z) - z_a \bigr)^2}\right)_{a \in K}\,,
\end{equation}
where $c(z)$ is an implicit normalization constant, such that the vector lies in the simplex~$\cS$ and $c(z) > z_a$ for all $a \in [K]$.
This constant $c(z)$ is in fact the Lagrange multiplier associated with the constraint $ p_1 + \ldots + p_K = 1$.
See \citet{zimmert2018optimal} for more details on how to compute $c(z)$ efficiently, see also \citet{audibert2014regret}.
To compute the mixabity gap, rewrite
\begin{align}
\nonumber
\delta^{\ts}_t & =
\max_{p \in \mathcal S} \left\{ \dotp{p_t - p}{ - \what y_t } - \frac{H_{1/2}(p) - H_{1/2}(p_t) - \dotp{\nabla H_{1/2}(p_t)}{p - p_t}}{\eta_t} \right\} \\
\label{eq:computeBHgap}
& = \dotp{p_t}{- \what y_t} +  \frac{H_{1/2}(p_t) }{\eta_t}- \frac{\dotp{\nabla H_{1/2} (p_t)}{p_t}}{\eta_t}
+ \frac{1}{\eta_t} \max_{p \in \mathcal S} \Big\{ \dotp{p}{\nabla H_{1/2}(p_t) + \eta_t \what y_t} - H_{1/2}(p) \Big\}\,,
\end{align}
where the maximum in the left-most side of these equalities
can be computed efficiently, thanks to~\eqref{eq:opt_pb}.

\begin{figure}[t]
	\renewcommand\footnoterule{}
	\begin{algorithm}[H]
		\begin{algorithmic}[1]
			\caption{\label{algo:tsallis_ftrlM} AdaFTRL with Tsallis entropy for $K$--armed bandits with a known payoff upper bound}
			\label{alg:tsallis_ftrl}
			\STATE \textbf{Input:} an upper bound $M$ on the payoffs
			\STATE \textbf{Initialization:} $\eta_{1} = +\infty$ and $p_{1} = (1/K,\ldots,1/K)$
			\FOR{rounds $t = 1,\,2,\, \dots$}
			\STATE Draw an arm $A_t \sim p_t$ (independently at random according to the distribution $p_t$)
			\STATE Get and observe the payoff $y_{t, A_t}$
			\STATE Compute the estimates of all payoffs
			\[
			\what y_{t, a} = \frac{y_{t, a}-M}{p_{t, a}} \1{A_t = a} + M
			\]
			\STATE Compute the mixability gap $\delta^{\ts}_t$ based on the distribution $p_t$ and on these estimates, e.g., using the
efficient implementation stated around~\eqref{eq:computeBHgap}:
			\[
				\delta^{\ts}_t = \max_{p \in \mathcal S} \left\{ \dotp{p_t - p}{ - \what y_t } - \frac{B_{H_{1/2}}(p, p_t)}{\eta_t} \right\}
			\]
			\STATE Compute the learning rate $\displaystyle{\eta_{t+1} = 2 \Biggl( \sum_{s = 1}^t \delta^{\ts}_s \Biggr)^{-1}} \bigl( \sqrt K -1 \bigr)$
			\STATE Define $p_{t+1}$ as
			\[
			p_{t+1} \in \argmin_{p \in \mathcal S}  \left\{ - \sum_{a = 1}^K p_{a} \sum_{s =1}^{t} \what y_{s, a}
- \frac{1}{\eta_{t+1}} \sum_{a = 1}^K 2\sqrt{p_a} \, \right\},
			\]
where an efficient implementation is provided by, e.g., \eqref{eq:opt_pb}
			\ENDFOR
		\end{algorithmic}
	\end{algorithm}
\end{figure}

\paragraph{Analysis of the algorithm.} We provide the following performance bound.

\begin{theorem}
\label{thm:tsallis}
AdaFTRL with $1/2$--Tsallis entropy for $K$--armed bandits relying on an upper bound $M$
on the payoffs (Algorithm~\ref{alg:tsallis_ftrl})
ensures that for all $m \in \R$ with $m \leq M$,
for all oblivious individual sequences $y_1,y_2,\ldots$ in $[m,M]^K$,
for all $T \geq 1$,
\[
R_T(y_{1:T}) \leq 4 (M-m)\sqrt{KT} + 2(M-m)\,.
\]
\end{theorem}

As in Section~\ref{app:proofs_exp3}, the proof scheme is a combination of the AdaFTRL bound
of Theorem~\ref{thm:adaftrl_toolbox} (which is indeed applicable), together with an improved bound on
the mixability gap that exploits the specific shape of the estimates. This bound is stated in the
next lemma, which is much similar to Lemma~\ref{lem:adahedge_up_bound_mix_gap}.

\begin{lemma}\label{eq:tsall_mix_gap_bound}
The mixability gaps of AdaFTRL with Tsallis entropy for $K$--armed bandits relying on an upper bound $M$
on the payoffs (Algorithm~\ref{alg:tsallis_ftrl})
are bounded, for all $m \in \R$ with $m \leq M$,
for all oblivious individual sequences $y_1,y_2,\ldots$ in $[m,M]^K$,
for all $t \geq 1$, by
\[
0 \leq \delta^{\ts}_t \leq M - m \qquad \mbox{and} \qquad \frac{\delta^{\ts}_t}{\eta_t} \leq  p_{t, A_t}^{-1/2} (M- y_{t, A_t})^2\,.
\]
\end{lemma}

The proof of Lemma~\ref{eq:tsall_mix_gap_bound} is postponed to the end of this section
and we now proceed with the proof of Theorem~\ref{thm:tsallis}. \medskip

\bpc{of Theorem~\ref{thm:tsallis}}
The structure of the proof is much similar to the one of Theorem~\ref{thm:adv_any_exploration--simpler},
which is why we only sketch our arguments.
The bound of Theorem~\ref{thm:adaftrl_toolbox} is applicable.
We use Lemma~\ref{eq:tsall_mix_gap_bound} with~\eqref{eq:gen_sq_mix_gaps} to see that
\begin{equation}
\bigg(\sum_{t = 1}^T \delta_t^{\ts} \bigg)^2
\leq  2 D_{H_{1/2}} \sum_{t = 1}^T p_{t, A_t}^{-1/2} (M- y_{t, A_t})^2 + (M-m)\sum_{t = 1}^T \delta_t^{\ts}\,.
\end{equation}
Again, using the fact that for all $a, b, x \geq 0$, the inequality  $x^2 \leq a + bx$ implies $x \leq \sqrt a + b$ :
\begin{equation}\label{eq:mix_sum_tsal}
\sum_{t = 1}^T \delta_t^{\ts} \leq \sqrt{2D_{H_{1/2}} \sum_{t = 1}^T p_{t, A_t}^{-1/2} (M - y_{t, A_t})^2} + (M-m)
\end{equation}
By~\eqref{eq:gen_ftrl_mixsum}, by taking expectations, and by Jensen's inequality:
\begin{equation}
R_T(y_{1:T})
\leq 2 \E \big[ \sum_{t = 1}^T \delta_t^{\ts} \big]
\leq 2\sqrt{2D_{H_{1/2}} \sum_{t = 1}^T \E \Bigl[ p_{t, A_t}^{-1/2} (M- y_{t, A_t})^2 \Bigr]} + 2(M-m)\,.
\end{equation}
We conclude by observing that for all $t$, by definition of the payoff estimates,
\begin{align*}
\E \! \left[ p_{t, A_t}^{-1/2}  \big(M- y_{t, A_t}\big)^2\right]
= \E \! \left[ \sum_{a = 1}^K p_{t, a} \, p_{t, a}^{-1/2}  \big(M- y_{t, a}\big)^2\right]
& \leq (M-m)^2 \, \E \!\left[ \sum_{a =1}^K \sqrt{p_{a, t}} \right] \\
& \leq (M-m)^2 \sqrt{K}\,,
\end{align*}
where the last inequality follows from the concavity of the square root.
The final claim is obtained by bounding the diameter $D_{H_{1/2}}$ by $2 \sqrt K$.
\end{proof}

We conclude this section by providing a proof of Lemma~\ref{eq:tsall_mix_gap_bound}. \medskip

\bpc{of Lemma~\ref{eq:tsall_mix_gap_bound}}
The fact that $\delta^{\ts}_t \geq 0$ holds actually for all regularizers
and can be seen from the definition \eqref{eq:def_gen_mix_gap} with $p = p_t$. For the inequality $\delta^{\ts}_t \leq M-m$,
we start with elementary manipulations of the definition of the mixability gap~\eqref{eq:def_gen_mix_gap}.
Denoting by $\vec M$ the vector with coordinates $(M, \dots, M)$
and noting that $\dotp{p_t-q}{\vec M} = 0$ for all $q \in \mathcal S$, we have
\begin{equation}
\label{eq:Lm3rewrdeltats}
\delta^{\ts}_t
= \max_{q \in \mathcal S} \left\{ \dotp{p_t - q}{- \what y_t } - \frac{B_{H_{1/2}}(q, p_t)}{\eta_t} \right\}
=\max_{q \in \mathcal S} \left\{ \dotp{p_t - q}{\vec M- \what y_t } - \frac{B_{H_{1/2}}(q, p_t)}{\eta_t} \right\}.
\end{equation}
Since all the coordinates of $\vec M - \what y_{t}$ are non-negative and by non-negativity of the Bregman divergence, this implies that
\[
\delta^{\ts}_t
\leq \dotp{p_t}{\vec M - \what y_t}
= M - y_{A_t, t} \leq M -m\,.
\]

We now prove the second inequality; we may assume that $\eta_t < + \infty$, as the bound holds trivially otherwise.
By Proposition~\ref{prop:opt_condition} (and by calculations similar to the ones performed in the proof of
Reminder~\ref{rem:ftrl}) the maximum in the rewriting~\eqref{eq:Lm3rewrdeltats}
of $\delta^{\ts}_t$ is achieved on the interior of the domain of $H_{1/2}$, which equals~$(0, + \infty)^K$,
thus in the interior of $\cS$. We therefore only need to prove that
\begin{equation}
\label{lm:AppA4Lemma3-needtoprove}
\forall q \in \interior(\cS), \qquad
\dotp{p_t - q}{\vec M- \what y_t } - \frac{B_{H_{1/2}}(q, p_t)}{\eta_t} \leq
\eta_t \, p_{t, A_t}^{-1/2} (M- y_{t, A_t})^2\,.
\end{equation}
We fix such a $q \in \interior(\cS)$, i.e., such that $q_a > 0$ for all $a$. We consider two cases.
First, if $q_{A_t} \geq p_{t, A_t}$, then, given the observations made after~\eqref{eq:loss_estimates_known_M},
\[
\dotp{p_t - q}{ \vec M- \what y_t }  - \frac{B_{H_{1/2}}(q, p_t)}{\eta_t}
= \underbrace{\left(  \frac{M - y_{t, A_t}}{p_{t, A_t}}\right)}_{\geq 0} \underbrace{\big( p_{t, A_t} - q_{A_t}  \big)}_{\leq 0}
- \frac{B_{H_{1/2}}(q, p_t)}{\eta_t} \leq 0\,.
\]

Otherwise, when $q_{A_t} < p_{t, A_t}$, a standard way of bounding the mixability gap, detailed below,
indicates that
\begin{equation}
\label{eq:lm3--stdbdmixgap}
\dotp{p_t - q}{ M- \what y_t } - \frac{B_{H_{1/2}}(q, p_t)}{\eta_t}
\leq \frac{\eta_t}{2}  \dotpB{\vec M- \what y_t}{ \nabla^2 H_{1/2}(z)^{-1} \, \bigl( \vec M- \what y_t \bigr)}\,,
\end{equation}
where $z$ is some probability distribution of the open segment $\segm(q, p_t)$ between $q$ and $p_t$, and
where $\nabla^2 H_{1/2}(z)^{-1}$ denotes the inverse of the positive definite Hessian of
$H_{1/2}$ at $z$. Since at $w \in (0, + \infty)^K$, the function $H_{1/2}$ is indeed twice differentiable, with
\[
\nabla H_{1/2}(w) = \bigl( -w_a^{-1/2} \bigr)_{a \in[K]}
\qquad {\mbox{and}} \qquad
\nabla^2 H_{1/2}(w) = \mathrm {Diag}\big( w_a^{-3/2} /2\big)_{a \in [K]}\,,
\]
we have $\nabla^2 H_{1/2}(z)^{-1} = \mathrm {Diag}\big( 2 z_a^{3/2} \big)_{a \in [K]}$.
We substitute this value into~\eqref{eq:lm3--stdbdmixgap}
and recall that the vector $\vec M - \what y_t$ has null coordinates except for its $A_t$--th coordinate:
\[
\frac{\eta_t}{2}  \dotpB{\vec M- \what y_t}{ \nabla^2 H_{1/2}(z)^{-1} \, \bigl( \vec M- \what y_t \bigr)}
= \eta_t \,  z_{A_t}^{3/2} \bigl( M- \what y_{t, A_t} \bigr)^2\,.
\]
Finally, remember that $z$ lies in the open segment $\segm(q, p_t)$
and that we assumed  $q_{A_t} < p_{t, A_t}$; we thus also have $z_{A_t} < p_{t, A_t}$.
As a consequence, using the very definition of $\what y_{t, A_t}$,
\[
\eta_t \,  z_{A_t}^{3/2} \bigl( M- \what y_{t, A_t} \bigr)^2
\leq
\eta_t\,  p_{t, A_t}^{3/2} \bigl( M- \what y_{t, A_t} \bigr)^2
= \eta_t\,  p_{t, A_t}^{-1/2}(M- y_{t, A_t})^2\,.
\]
Therefore, in all cases, that is, whether $q_{A_t} \geq p_{t, A_t}$ or $q_{A_t} < p_{t, A_t}$,
the bound~\eqref{lm:AppA4Lemma3-needtoprove} is obtained. It only remains to prove
the standard inequality~\eqref{eq:lm3--stdbdmixgap}.

This inequality is essentially stated as Therorem~26.13 in \citet{lattimore_szepesvari_2020} but we provide a proof
for the sake of completeness.
As we assumed that $\eta_t < + \infty$, we have (as above, by Proposition~\ref{prop:opt_condition})
that $p_t$ lies in the interior of $\cS$. In particular, as both $p_t$ and $q$ are in the interior of $\cS$,
the function $H_{1/2}$ is $\mathcal C^2$ over the closed segment $\overline \segm(q, p_t)$
between $q$ and $p_t$. Therefore, by the mean-value theorem, there exists $z$ in the open segment $\segm(q, p_t)$ such that
\[
\underbrace{H_{1/2}(q) - H_{1/2}(p_t) - \dotp{\nabla H_{1/2}(p_t)}{q-p_t}}_{=B_{H_{1/2}}(q, p_t)}
= \frac{1}{2}\dotpB{q - p_t}{ \nabla^2 H_{1/2}(z) \, (q-p_t)}\,.
\]
It is useful to introduce the standard notation from convex analysis for the local norm
(which is indeed a norm because the Hessian is positive definite):
\[
\lVert q-p_t \rVert_{\nabla^2 H_{1/2}(z)}^2 \eqdef \dotpB{q - p_t}{ \nabla^2 H_{1/2}(z) \, (q-p_t)}\,.
\]
We therefore have so far the rewriting:
\[
- \frac{B_{H_{1/2}}(q, p_t)}{\eta_t} = - \frac{1}{2\eta_t} \dotpB{q - p_t}{ \nabla^2 H_{1/2}(z) \, (q-p_t)}\,.
\]
Now, by the Cauchy-Schwarz inequality,
\begin{align*}
\dotp{p_t - q}{ \vec M- \what y_t }
&= 	\dotpB{\nabla^2 H_{1/2}(z)^{1/2} \, (p_t - q)}{ \nabla^2 H_{1/2}(z)^{-1/2} \, \big(\vec M- \what y_t \big)} \\
&\leq \lVert p_t - q \rVert_{\nabla^2 H_{1/2}(z)} \,\, \lVert \vec M- \what y_t \rVert_{\nabla^2 H_{1/2}(z)^{-1}}\,.
\end{align*}
Combining the rewriting and the bound above, we get
\begin{align*}
\lefteqn{\dotp{p_t - q}{ M- \what y_t } - \frac{B_{H_{1/2}}(q, p_t)}{\eta_t}} \\
&\leq \lVert p_t - q \rVert_{\nabla^2 H_{1/2}(z)}  \lVert \vec M- \what y_t \rVert_{\nabla^2 H_{1/2}(z)^{-1}} - \frac{1}{2\eta_t} \lVert q-p_t \rVert_{\nabla^2 H_{1/2}(z)}^2 \\
&\leq \frac{\eta_t}{2}  \lVert \vec M- \what y_t \rVert_{\nabla^2 H_{1/2}(z)^{-1}}^2\,,
\end{align*}
where we used $ab - b^2 / 2\leq a^2 / 2$ to get the second inequality.
This is exactly~\eqref{eq:lm3--stdbdmixgap}.
\end{proof}

\end{document}